\documentclass[10pt,reqno]{amsart}
\usepackage{cmap}
\usepackage[a4paper,margin=25mm,bottom=25.5mm,headheight=14pt]{geometry}
\usepackage[T1]{fontenc}
\usepackage[utf8]{inputenc}
\usepackage{lmodern,amsmath,amssymb,amsthm,mathtools}
\relax 
\usepackage{booktabs,array,xcolor,etoolbox}
\definecolor{linkblue}{RGB}{0,45,120}
\usepackage[colorlinks=true,linkcolor=linkblue,citecolor=linkblue,urlcolor=linkblue,
bookmarksnumbered=true]{hyperref}

\allowdisplaybreaks[2]
\numberwithin{equation}{section}

\newtheoremstyle{annalsplain}
{8pt plus 2pt minus 2pt}
{8pt plus 2pt minus 2pt}
{\itshape}
{}
{\scshape}
{.}
{0.5em}
{}

\newtheoremstyle{annalsremark}
{7pt plus 2pt minus 2pt}
{7pt plus 2pt minus 2pt}
{\normalfont}
{}
{\scshape}
{.}
{0.5em}
{}

\theoremstyle{annalsplain}
\newtheorem{theorem}{Theorem}[section]
\newtheorem{proposition}[theorem]{Proposition}
\newtheorem{lemma}[theorem]{Lemma}
\newtheorem{corollary}[theorem]{Corollary}

\theoremstyle{annalsremark}
\newtheorem{remark}[theorem]{Remark}

\makeatletter
\renewcommand\subsection{\@startsection{subsection}{2}{\z@}%
	{.8\linespacing\@plus.4\linespacing}{.4\linespacing}%
	{\normalfont\bfseries}}
\renewcommand{\@settitle}{%
	\begin{center}
		\fontsize{10.2}{13.2}\fontseries{b}\selectfont
		\advance\topskip by 14pt\relax
		\uppercasenonmath\@title
		\@title\par
		\vspace{-2.25pt}
	\end{center}%
}
\patchcmd{\@setabstracta}{\skip@20\p@}{\skip@16.7\p@}{}{}
\patchcmd{\@setaddresses}{\footnotesize}{%
	\footnotesize
	\setlength{\baselineskip}{10.25pt}%
	\setlength{\parskip}{0pt}%
	\setlength{\bigskipamount}{7pt}%
	\setlength{\spaceskip}{2.5pt plus 1.25pt minus 1.1pt}%
}{}{}
\AddToHook{cmd/maketitle/after}{\enlargethispage{-6.5pt}}
\patchcmd{\ps@firstpage}{\hfil\thepage\hfil}{\hfil\raisebox{9.5pt}[0pt][0pt]{\thepage}\hfil}{}{}
\makeatother

\newcommand{\dd}{\,\mathrm{d}}
\newcommand{\Div}{\operatorname{div}}

\hypersetup{
	pdftitle={Nonlinear stability and time-decay rates of planar rarefaction waves for the full compressible Navier--Stokes--Fourier system in R3},
	pdfsubject={Nonlinear wave stability for the three-dimensional compressible Navier--Stokes--Fourier system},
	pdfauthor={},
	pdfkeywords={compressible Navier-Stokes-Fourier system, planar rarefaction waves, nonlinear stability, time-decay rates, relative entropy}
}

\title[Nonlinear stability of planar rarefaction waves in $\mathbb{R}^3$]{Nonlinear stability and time-decay rates of planar rarefaction waves\\for the full compressible Navier--Stokes--Fourier system in $\mathbb{R}^3$}

\author[R.-J. Duan]{Renjun Duan}
\address[R.-J. Duan]{Department of Mathematics, The Chinese University of Hong Kong, Shatin, Hong Kong, China}
\email{rjduan@math.cuhk.edu.hk}

\author[F.~Huang]{Feimin Huang}
\address[F.~Huang]{Academy of Mathematics and Systems Science, Chinese Academy of Sciences, Beijing 100190, China, and School of Mathematical Sciences, University of Chinese Academy of Sciences, Beijing 100049, China.}
\email{fhuang@amt.ac.cn}

\author[R.~Li]{Rui Li}
\address[R.~Li]{Department of Mathematics, The Chinese University of Hong Kong, Shatin, Hong Kong, China.}
\email{ruili001@cuhk.edu.hk}

\author[L.~Xu]{Lingda Xu}
\address[L.~Xu]{Department of Mathematics, The Chinese University of Hong Kong, Hong Kong, China.}
\email{lingdaxu@cuhk.edu.hk}

\subjclass[2020]{Primary 35Q30, 35L65; Secondary 76N10, 35B40.}
\keywords{Compressible Navier--Stokes--Fourier system, planar rarefaction waves, nonlinear stability, algebraic time-decay rates, relative entropy, negative Sobolev spaces.}

\begin{document}
	
	\begin{abstract}
		We prove the global-in-time nonlinear stability and establish explicit algebraic time-decay rates of planar rarefaction waves for the three-dimensional full compressible Navier--Stokes--Fourier system of an ideal polytropic gas in $\mathbb{R}^3$. Given a smooth initial rarefaction profile of small wave strength \(\delta>0\), we first construct the corresponding global one-dimensional Navier–Stokes–Fourier solution \(U_b(t,x_1)\) and show that it converges toward the self-similar Euler rarefaction fan in \(L^\infty(\mathbb{R})\) at rate \((1+t)^{-12/25}\). For the three-dimensional Cauchy problem under small weighted $H^3(\mathbb{R}^3)$ perturbations with zero integral of mass, momentum, and total energy, we prove that there exists a unique global strong solution $U(t,x)$ converging toward $U_b(t,x_1)$ in $L^\infty(\mathbb{R}^3)$ at rate $(1+t)^{-11/25}$, with $\nabla(U-U_b)$ decaying in $H^2(\mathbb{R}^3)$ at rate $(1+t)^{-11/25}$.
		
	\end{abstract}
	
	\maketitle
	
	\section{Introduction}
	\label{sec:main-results}
	\label{r5main:problem}
	
	\subsection{The physical system and background}
	
	The motion of a compressible, viscous, heat-conducting Newtonian fluid in $\mathbb{R}^3$ is governed by the three-dimensional full compressible Navier--Stokes--Fourier system:
	\begin{equation}
		\left\{
		\begin{aligned}
			&\rho_t+\Div(\rho u)=0,\\
			&(\rho u)_t+\Div(\rho u\otimes u)+\nabla p=\Div\tau,\\
			&E_t+\Div\bigl((E+p)u\bigr)=\Div(\tau u+\kappa\nabla\theta),
		\end{aligned}
		\right.
		\label{r5main:NSF}
	\end{equation}
	for $t>0$ and $x=(x_1,x_\perp)\in\mathbb{R}\times\mathbb{R}^2$. Here $\rho>0$, $u=(u_1,u_2,u_3)\in\mathbb{R}^3$, and $\theta>0$ denote the fluid density, velocity, and absolute temperature, respectively. For an ideal polytropic gas, the pressure $p$, total energy density $E$, specific heats $c_v,c_p$, and viscous stress tensor $\tau$ are given by
	\[
	\left\{
	\begin{aligned}
		&p=R_g\rho\theta,\qquad
		E=\rho\left(c_v\theta+\frac{|u|^2}{2}\right),\qquad
		c_v=\frac{R_g}{\gamma-1},\qquad c_p=\gamma c_v,\\
		&\tau=\mu(\nabla u+\nabla u^\top)+\lambda_v(\Div u)I_3.
	\end{aligned}
	\right.
	\]
	Throughout the paper, the adiabatic exponent $\gamma$, specific gas constant $R_g$, viscosity coefficients $\mu,\lambda_v$, and thermal conductivity $\kappa$ are fixed constants satisfying
	\[
	\gamma>1,\qquad R_g>0,\qquad \mu>0,\qquad \kappa>0,\qquad 2\mu+3\lambda_v\ge0,
	\]
	so that $\nu=2\mu+\lambda_v>0$. We write $U=(\rho,u,\theta)$ for the primitive variables and $\mathbf{U}=(\rho,\rho u,E)$ for the conserved variables.
	
	A classical problem in compressible fluid dynamics is the large-time asymptotic stability of basic Riemann wave patterns---shock waves, rarefaction waves, and contact discontinuities---connecting distinct end states $U_\pm$ as $x_1\to\pm\infty$. In one space dimension, the nonlinear stability theory of viscous conservation laws has been developed since the works of Il'in and Oleinik \cite{IO60}, Goodman \cite{Goodman86}, Matsumura and Nishihara \cite{MN85,MN86,MN92}, and Liu \cite{Liu85}. For one-dimensional rarefaction waves, Matsumura and Nishihara \cite{MN86,MN92} constructed a smooth approximate profile from the Burgers equation and used the positivity of the velocity slope to prove global $L^2(\mathbb{R})$ stability. This method was extended to the full Navier--Stokes--Fourier system by Liu and Xin \cite{LX88} and Xin \cite{Xin88}, to large rarefaction waves by Nishihara, Yang, and Zhao \cite{NYZ04}, and to pointwise Green's function bounds and composite waves by Liu and Zeng \cite{LZ97}, Huang, Matsumura, and Xin \cite{HMX06}, Huang, Xin, and Yang \cite{HXY08}, and Huang, Li, and Matsumura \cite{HLM10}; see also the recent work of Li, Song, Wang, and Xu \cite{LSWX26} on limiting pointwise decay for the one-dimensional isentropic system. However, because rarefaction waves expand in time and do not admit an $L^2(\mathbb{R})$ anti-derivative formulation, classical $H^2(\mathbb{R})$ energy methods gave only qualitative convergence as $t\to\infty$ without any explicit decay rate; obtaining algebraic time-decay rates toward rarefaction waves for the compressible Navier--Stokes equations has remained a long-standing open problem, and the first explicit decay rate in the $H^2(\mathbb{R})$ framework for the one-dimensional isentropic system was obtained only very recently by He, Huang, Liu, and Xu \cite{HHLX26}.
	
	In several space dimensions, wave propagation in the transverse directions and the coupling among velocity components introduce new difficulties. For scalar viscous conservation laws in $\mathbb{R}^n$, Xin \cite{Xin90} proved the asymptotic stability of planar rarefaction waves, and Huang and Xu \cite{HX24} established explicit algebraic decay rates; Hokari and Matsumura \cite{HM97} treated a two-dimensional barotropic model with artificial viscosity. For the physical compressible Navier--Stokes equations, Li and Wang \cite{LW18} proved the nonlinear stability of planar rarefaction waves for the two-dimensional isentropic system in a periodic channel $\mathbb{R}\times\mathbb{T}$, Li, Wang, and Wang \cite{LWW18} proved global stability for the three-dimensional full Navier--Stokes--Fourier system on $\mathbb{R}\times\mathbb{T}^2$, and Huang, Xu, and Yuan \cite{HXY22} proved the asymptotic stability of planar rarefaction waves under general three-dimensional periodic perturbations (see also Hou, Liu, Wang, and Xu \cite{HLSW24} for the vanishing viscosity limit to planar rarefaction waves with vacuum for the three-dimensional full system).
	
	Passing from the periodic channel $\mathbb{R}\times\mathbb{T}^2$ to the whole space $\mathbb{R}^3$ changes the analytical structure of the problem, and prior to the present work neither global stability on $\mathbb{R}^3$ nor explicit time-decay rates of planar rarefaction waves had been obtained for the three-dimensional full system \eqref{r5main:NSF}. On $\mathbb{R}\times\mathbb{T}^2$, the cross-section $\mathbb{T}^2$ is compact, so the residual of a smooth Euler rarefaction profile belongs to $L^2(\mathbb{R}\times\mathbb{T}^2)$ and can be controlled by the relative entropy method \cite{LWW18,HXY22}; moreover, nonzero transverse Fourier modes satisfy a Poincar\'e inequality, while the zero mode reduces to a one-dimensional system that can be integrated in $x_1$. On $\mathbb{R}^3$, any nonzero residual depending only on $x_1$ has infinite $L^2(\mathbb{R}^3)$ norm, the three-dimensional divergence $\Div=\partial_1+\nabla_\perp\cdot$ prevents one-dimensional integration, and the transverse frequency $\zeta\in\mathbb{R}^2$ can be arbitrarily small.
	
	In this paper, we perturb around the exact one-dimensional viscous rarefaction wave $U_b(t,x_1)$ of \eqref{r5main:NSF}, establish its explicit algebraic convergence rate toward the Euler rarefaction wave by adopting the first-order cross-family derivative correction method developed in our works on compressible Couette flow \cite{HLX25,LWXZ26} and one-dimensional pointwise decay \cite{LSWX26}, and combine the $H^3(\mathbb{R}^3)$ energy with a weighted $\dot{H}^{-1}(\mathbb{R}^3)$ energy based on the characteristic decomposition of $\partial_1 U_b$. This yields both global nonlinear stability on $\mathbb{R}^3$ and explicit algebraic decay rates for the planar viscous wave and its three-dimensional perturbations.
	
	\subsection{Formulation of the planar wave and main result}
	
	By a Galilean transformation, we normalize the left far-field state to have zero velocity, $U_-=(\rho_-,0,\theta_-)$ with $\rho_->0$ and $\theta_->0$. The left sound speed and the characteristic speeds of the forward 3-rarefaction wave are given by
	\[
	c_-=(\gamma R_g\theta_-)^{1/2},\qquad W_-=c_-,\qquad W_+=W_-+\delta,
	\]
	where $\delta>0$ is the wave strength. Fix a smoothing parameter $e>0$, and let $w^E(t,x_1)$ be the smooth solution to the inviscid Burgers Cauchy problem
	\begin{equation}
		\left\{
		\begin{aligned}
			&w^E_t+w^Ew^E_{x_1}=0,\\
			&w^E(0,x_1)=\frac{W_++W_-}{2}+\frac{\delta k_2}{2}\int_0^{e x_1}(1+z^2)^{-2}\dd z,
			\qquad k_2=\frac{4}{\pi}.
		\end{aligned}
		\right.
		\label{main:burgers-data}
	\end{equation}
	Setting $W^E(t,x_1)=w^E(t+1,x_1)$ and $e_1=(1,0,0)$, we define the smooth Euler rarefaction profile $\overline{U}_E=(\bar{\rho},\bar{u}e_1,\bar{\theta})$ by the 3-rarefaction Riemann invariant relations
	\[
	\left\{
	\begin{aligned}
		&\bar{c}=c_-+\frac{\gamma-1}{\gamma+1}(W^E-W_-),\qquad
		\bar{u}=\frac{2}{\gamma+1}(W^E-W_-),\\
		&\bar{\rho}=\rho_-\left(\frac{\bar{c}}{c_-}\right)^{2/(\gamma-1)},\qquad
		\bar{\theta}=\theta_-\left(\frac{\bar{c}}{c_-}\right)^2.
	\end{aligned}
	\right.
	\]
	We denote by
	\[
	U_b(t,x_1)=\bigl(\rho_b(t,x_1),u_b(t,x_1)e_1,\theta_b(t,x_1)\bigr)
	\]
	the one-dimensional planar solution to \eqref{r5main:NSF} with initial datum
	\begin{equation}
		U_b(0,x_1)=\overline{U}_E(0,x_1).
		\label{main:planar-data}
	\end{equation}
	The centered Euler rarefaction fan $U^R(t,x_1)$ is defined by the same state formulas with $W^E(t,x_1)$ replaced by $W^R(x_1/t)$, where
	\[
	W^R(\xi)=
	\begin{cases}
		W_-,&\xi\le W_-,\\
		\xi,&W_-<\xi<W_+,\\
		W_+,&\xi\ge W_+.
	\end{cases}
	\]
	
	Our main theorem is stated as follows.
	
	\begin{theorem}
		\label{thm:main}
		Fix a weight exponent $m_0>5/2$, the physical constants $(\gamma,R_g,\mu,\lambda_v,\kappa)$, the left state $U_-=(\rho_-,0,\theta_-)$. There exists $\delta_*>0$ such that for every $0<\delta\le\delta_*$, the planar Cauchy problem for \eqref{r5main:NSF} with initial datum \eqref{main:planar-data} admits a unique global smooth solution $U_b(t,x_1)$ taking values in a compact set of positive states.
		
		For each such \(\delta\), there exist constants $\varepsilon_*>0$ and $C_\delta<\infty$ with the following property.
		Let $f_0(x)$ be an initial perturbation of the conserved variables satisfying
		\begin{equation}
			\int_{\mathbb{R}^3}f_0(x)\dd x=0,\qquad
			\varepsilon=\sum_{|\alpha|\le3}
			\|\langle x\rangle^{m_0}\partial_x^\alpha f_0\|_{L^2(\mathbb{R}^3)}
			\le\varepsilon_*,
			\label{main:data-class}
		\end{equation}
		where $\alpha=(\alpha_1,\alpha_2,\alpha_3)\in\mathbb{N}_0^3$ and $\langle x\rangle=(1+|x|^2)^{1/2}$. Then the three-dimensional Cauchy problem for \eqref{r5main:NSF} with initial datum
		\[
		\mathbf{U}(0,x)=\overline{\mathbf{U}}_E(0,x_1)+f_0(x)
		\]
		admits a unique global strong solution $U=(\rho,u,\theta)$ with positive density and temperature. Its perturbation $Z=U-U_b=(Z_\rho,Z_u,Z_\theta)$ relative to the one-dimensional solution \(U_b(t,x_1)\) satisfies
		\[
		\left\{
		\begin{aligned}
			&Z\in C\bigl([0,\infty);H^3(\mathbb{R}^3)\bigr),\qquad
			\nabla Z_\rho\in L^2_{\mathrm{loc}}\bigl([0,\infty);H^2(\mathbb{R}^3)\bigr),\\
			&\nabla(Z_u,Z_\theta)\in L^2_{\mathrm{loc}}\bigl([0,\infty);H^3(\mathbb{R}^3)\bigr),
		\end{aligned}
		\right.
		\]
		and the perturbation obeys the \(L^\infty\) decay estimate
		\begin{equation}
			\|U(t)-U_b(t)\|_{L^\infty(\mathbb{R}^3)}
			\le C_\delta\varepsilon (1+t)^{-11/25},\qquad t\ge0.
			\label{main:Linfty}
		\end{equation}
		Moreover, \(U\) converges to the smooth Euler rarefaction profile \(\overline U_E\) and to the self-similar Euler rarefaction fan \(U^R\) with the respective bounds
		\begin{align}
			\|U(t)-\overline{U}_E(t)\|_{L^\infty(\mathbb{R}^3)}
			&\le C_\delta\bigl(\varepsilon (1+t)^{-11/25}+(1+t)^{-12/25}\bigr),
			&&t\ge0,\label{main:limit}\\
			\|U(t)-U^R(t)\|_{L^\infty(\mathbb{R}^3)}
			&\le C_\delta\bigl(\varepsilon (1+t)^{-11/25}+(1+t)^{-12/25}\bigr),
			&&t\ge1,\label{main:selfsimilar-limit}
		\end{align}
		where $C_\delta$ is independent of $t$ and $f_0$, and depends only on the fixed physical constants and $\delta$.
	\end{theorem}
	
	\begin{remark}
		The condition $\int_{\mathbb{R}^3}f_0(x)\dd x=0$ in \eqref{main:data-class} means that the initial perturbation has zero total mass, momentum, and energy. The derivative order $N=22$ used in Sections~\ref{sec:elementaryplanarinputs}--\ref{sec:allfamilyLone} is needed only for the smooth one-dimensional solution $U_b(t,x_1)$; the three-dimensional perturbation requires only $H^3(\mathbb{R}^3)$ regularity.
	\end{remark}
	
	\subsection{Main difficulties}
	
	\subsubsection{The viscous background, 1D rarefaction decay, and sign reversal in $\dot{H}^{-1}(\mathbb{R}^3)$}
	In one space dimension and in the periodic channel $\mathbb{R}\times\mathbb{T}^2$ \cite{MN86,LX88,LWW18,HXY22}, one can perturb directly around the smooth Euler profile $\overline{U}_E$, keep its viscous residual as a source term, and control it using the positive slope $\partial_1\bar{u}>0$ together with derivative dissipation. On $\mathbb{R}^3$, a nonzero planar residual does not depend on $x_\perp\in\mathbb{R}^2$ and therefore does not belong to $L^2(\mathbb{R}^3)$. Taking the exact viscous solution $U_b(t,x_1)$ as the background eliminates this source term, but the difference $\partial_1(U_b-\overline{U}_E)$ has no definite sign and produces cross terms of the form
	\[
	C\int_{\mathbb{R}^3}|\partial_1(U_b-\overline{U}_E)|\,|Z|^2\dd x,
	\qquad Z=U-U_b,
	\]
	whereas the physical dissipation controls only spatial derivatives of $Z$.
	
	Controlling these cross terms requires an explicit time-decay rate for the one-dimensional viscous rarefaction wave $U_b-\overline{U}_E$ and its derivatives. However, proving explicit algebraic decay rates toward rarefaction waves for the compressible Navier--Stokes system is itself a long-standing open problem: unlike shock waves or contact discontinuities, rarefaction waves do not admit an $L^2(\mathbb{R})$ anti-derivative formulation, so classical energy methods \cite{MN86,MN92,LX88} yield only qualitative convergence $\|U_b(t)-\overline{U}_E(t)\|_{L^\infty(\mathbb{R})}\to 0$ as $t\to\infty$ without any rate. Only very recently did He, Huang, Liu, and Xu \cite{HHLX26} establish the first explicit decay rate toward rarefaction waves in the $H^2(\mathbb{R})$ energy framework for the one-dimensional isentropic system. For the full Navier--Stokes--Fourier system \eqref{r5main:NSF}, where viscous and thermal dissipation in the forward rarefaction wave excite both a backward acoustic wave and an entropy wave, obtaining explicit decay rates for $U_b-\overline{U}_E$ requires a new one-dimensional analysis.
	
	To estimate low transverse frequencies on $\mathbb{R}^3$, Section~\ref{sec:common-energy} introduces the potential $h$ solving the elliptic equation
	\[
	\mathcal{A}_t h=\widetilde{f}=f+K_l^*\partial_l f,\qquad
	E_{-1}(t)=\frac{1}{2}\mathcal{A}_t[h,h]
	=\frac{1}{2}\int_{\mathbb{R}^3}\mathcal{W}\mathcal{J}\dd x,
	\]
	where $f=\mathbf{U}-\mathbf{U}_b$, $K_l^*$ is the constant Kawashima matrix, and $\mathcal{J}$ is defined in \eqref{r5a:jet-form}. Passing from $L^2(\mathbb{R}^3)$ to $\dot{H}^{-1}(\mathbb{R}^3)$ via $\mathcal{A}_t^{-1}$ reverses the sign of the principal Euler transport term: in terms of $s_\pm$ in \eqref{r5a:acoustic-projections}, its leading acoustic part in $E_{-1}'(t)$ is
	\begin{equation}\label{intro:sign-reversal}
		\beta\int_{\mathbb{R}^3}\mathcal{W}
		\bigl(\mu_+s_+^2+\mu_-s_-^2\bigr)\dd x,
		\qquad
		\mu_\pm=\frac{1}{2}\left(\partial_1u_b
		\pm\frac{\partial_1p_b}{\rho_bc_E}\right),
		\qquad \beta=\frac{\gamma+1}{4}.
	\end{equation}
	Thus a positive slope $\mu_\pm>0$, which gives damping in $L^2$, produces growth in $\dot{H}^{-1}(\mathbb{R}^3)$.
	
	\subsubsection{Non-integrable cross terms and the allowable growth of the $\dot{H}^{-1}(\mathbb{R}^3)$ energy}
	Write $T=1+t$ and $g=\mathbf{U}_b-\mathbf{V}(W)$, where $W$ is the viscous Burgers profile. Even after we establish the explicit one-dimensional decay rates in Section~\ref{sec:planar-decay-all},
	\begin{equation}\label{intro:ux-bad}
		\|\partial_1V(W)(t)\|_{L^\infty}\le CT^{-1},\qquad
		\|g(t)\|_{L^\infty}\le CT^{-12/25},\qquad
		\|\partial_1g(t)\|_{L^\infty}\le CT^{-24/25},
	\end{equation}
	the derivative decay rate $T^{-24/25}$ is still not integrable in time. Consequently, even with \eqref{intro:ux-bad} in hand, the cross terms cannot be treated by a direct application of Gronwall's inequality, which would produce a factor $\exp(CT^{1/25})$ and prevent any algebraic decay.
	
	Instead, we establish a combined energy inequality for $\mathcal{F}=E_{-1}+L_*E_3$ of the form
	\[
	\mathcal{F}'+c\mathcal{D}\le k_*\mathcal{F},\qquad
	\int_0^t k_*(\tau)\dd\tau
	\le C_\delta+A_*\log T+B_*\log\log(2+t).
	\]
	In the two-dimensional Fourier-splitting argument of Section~\ref{sec:main-proof}, convergence of the dyadic sum for the nonlinear source,
	\[
	\sum_{j\ge0}2^{j(A_*/2-1)}(1+j)^C<\infty,
	\]
	requires $A_*<2$. Along the isentropic forward rarefaction curve, $2\beta\mu_+=\partial_1\lambda_+^E$ has leading upper bound $T^{-1}$, which already contributes $1$ to $A_*$. If the {coefficients $2\beta\mu_-$ and $2\beta\mu_+$} were bounded separately by $T^{-1}$ and added, the leading coefficient would be $2/T$, violating $A_*<2$. Hence the slope error $\partial_1 g$, the backward slope $\mu_-$, and the entropy derivative $S_1=\partial_1 S$ must all be estimated without adding another full unit to $A_*$.
	
	\subsubsection{The entropy gradient and coercivity}
	\label{intro:entropy-obstruction}
	The background entropy gradient $S_1=\partial_1S$ enters $E_{-1}'(t)$ through the integral
	\begin{equation}\label{intro:contact-vorticity}
		\int_{\mathbb{R}^3}\mathcal{W}w\theta_bS_1\,\Omega
		\cdot(e_1\times\nabla h_E)\dd x,
		\qquad \Omega=\operatorname{curl}(h_m+u_bh_Ee_1).
	\end{equation}
	{For the entropy production component $H_P$ from \eqref{ev:decomposition}, whose derivative enters $S_1$ in \eqref{intro:contact-vorticity}, $\|(H_P)_y\|_{L^\infty}\le C\delta T^{-1}\log(2+t)+k_P(t)$ with $k_P$ from \eqref{ev:HPy} satisfying $k_P\in L^1(0,\infty)$} near the trajectory $x_c(t)$. To avoid a factor $T^{-1}\log(2+t)$ in Gronwall's inequality, we choose the weight $w$ to be of order $[\log(2+t)]^{-2}$ near $x_c(t)$, so that $\sqrt{w}$ cancels $\log(2+t)$. However, when $w\approx[\log(2+t)]^{-2}\ll 1$, the coefficient $K=w-1$ is close to $-1$, and the integrand $|\nabla h_m|^2+K|\Omega|^2$ in \eqref{r5a:jet-form} is no longer pointwise nonnegative. Coercivity of $\mathcal{A}_t$ must therefore be proved using the divergence--curl structure of the integral.
	
	\subsection{Main ideas and proof outline}
	
	In Lagrangian mass coordinates $y=Y(t,x_1)$ with $Y_{x_1}=\rho_b$, let $r=u_b-I(p_b)$, $s=u_b+I(p_b)$, and $\sigma=S-S_-$ denote the backward acoustic, forward acoustic, and relative entropy variables of the planar wave. By \eqref{v3:state-recovery}, the background derivatives $\mu_\pm$ and $S_1=\partial_1 S$ in \eqref{r5a:explicit-acoustic-polynomial} can be written as
	\begin{equation}\label{intro:mu-forward}
		\left\{
		\begin{aligned}
			&\mu_+=\frac{\rho_b}{4}
			\left[\bigl(1+e^{\sigma/(2c_p)}\bigr)s_y
			+\bigl(1-e^{\sigma/(2c_p)}\bigr)r_y\right],\\
			&\mu_-=\frac{\rho_b}{4}
			\left[\bigl(1+e^{\sigma/(2c_p)}\bigr)r_y
			+\bigl(1-e^{\sigma/(2c_p)}\bigr)s_y\right],\\
			&S_1=\rho_b\sigma_y,
		\end{aligned}
		\right.
	\end{equation}
	evaluated at $y=Y(t,x_1)$. Thus the leading acoustic slopes are $\rho_bs_y/2$ and $\rho_br_y/2$ up to $O(\rho_b|\sigma|(|r_y|+|s_y|))$ terms, while the entropy wave has zero speed in mass coordinates ($Y_t+u_bY_{x_1}=0$). The proof proceeds in five main steps.
	
	\subsubsection{Decay of the one-dimensional viscous wave}
	In Section~\ref{sec:planar-decay-all}, we construct the planar solution $U_b$ and prove the algebraic decay of $g=\mathbf{U}_b-\mathbf{V}(W)$. Relative entropy and Kawashima energy estimates give uniform bounds in $H^N(\mathbb{R})$ ($N=22$) and derivative decay estimates conditional on the decay of $\|g(t)\|_{L^2}$, while weighted relative entropy bounds outside $\{|x|\le Ct\}$ and interpolation on dyadic annuli give time-integrable $L^1(\mathbb{R})$ bounds for $\partial_x^3 g$ and $\partial_x^4 g$.
	
	To estimate $\|g(t)\|_{L^1(\mathbb{R})}$, we project $g$ onto the three characteristic wave families $q=L(W)g=(q_-,q_0,q_+)^\top$. Because viscosity and thermal conductivity couple the three characteristic modes, the equation for $q_i$ contains off-diagonal second-order viscous terms $\sum_{j\ne i}D_{ij}(W)\partial_x^2 q_j$ that cannot be absorbed in $L^1(\mathbb{R})$. To overcome this difficulty, we adopt the first-order cross-family derivative correction method developed in our works on compressible Couette flow by Huang, Li, and Xu \cite{HLX25} and Li, Wang, Xu, and Zhang \cite{LWXZ26}, together with our recent work by Li, Song, Wang, and Xu \cite{LSWX26} (cf.\ (1.8) in \cite{LSWX26}): using the separation of wave speeds $\lambda_i-\lambda_j\ne 0$ ($i\ne j$), we add a first-order spatial derivative correction of the opposite characteristic families,
	\[
	G_i=q_i+\sum_{j\ne i}K_{ij}(W)\partial_x q_j,\qquad
	K_{ij}(W)=-\frac{D_{ij}(W)}{\lambda_i(W)-\lambda_j(W)}\quad (i\ne j),
	\]
	followed by the profile derivative correction $\widetilde{G}_i=G_i+\chi_i(W)W_x$ ($i\ne +$) and the quadratic normal-form correction \eqref{af:normalform}. Since $(\partial_t+\lambda_i\partial_x)q_j=(\lambda_i-\lambda_j)\partial_x q_j+\cdots$, the transport commutator $[\Lambda,K]\partial_x^2 q$ {cancels the off-diagonal viscous second derivatives of the linearized system}, reducing the system to the diagonally diffused wave equations \eqref{af:h-system} (the same derivative correction idea is also used in Lagrangian mass coordinates in Lemma~\ref{inc:reduction} and Proposition~\ref{prop:outgoing-reduction} below to eliminate cross-family viscous second derivatives in the backward and forward slope equations). Exponential weights in $W$ then control the linear coupling with $W_x$ for the backward and entropy modes, and a Glimm--Bony interaction functional controls products of waves from different families. This gives
	\[
	\|g(t)\|_{L^1(\mathbb{R})}
	\le C\bigl(\varepsilon_\delta+\delta\log(2+t)\bigr),
	\qquad \lim_{\delta\rightarrow 0}\varepsilon_\delta=0.
	\]
	Combining this $L^1$ bound with the one-dimensional Nash inequality closes the bootstrap with $a=12/25$:
	\[
	\|g(t)\|_{L^2}\le CT^{-a/2},\qquad
	\|g_x(t)\|_{L^2}\le CT^{-3a/2},\qquad
	\|g(t)\|_{L^\infty}\le CT^{-a}.
	\]
	In particular,
	\[
	K_b(t)=\sum_{j=1}^4\|\partial_1^jU_b(t)\|_{L^\infty}^2
	\le C_\delta T^{-48/25}.
	\]
	
	\subsubsection{The forward rarefaction slope}
	\label{intro:forward-estimate}
	On the exterior region where the pressure exponent in $\mathcal{W}$ equals a constant $m>0$, the forward acoustic part of \eqref{r5a:explicit-acoustic-polynomial} is
	\[
	\int_{\mathbb{R}^3}\zeta_e\mathcal{W}\mu_+
	\left(\mathcal{Q}_+-\frac{m}{2}\mathcal{F}_+^2\right)\dd x,
	\]
	where $\zeta_e$ is the exterior cutoff in \eqref{r5a:retained-D}, and $s_+$ and $\mathcal{F}_+^2$ are defined in \eqref{r5a:acoustic-projections}. The factor $(p_b/p_-)^{-m/\gamma}$ in \eqref{r5a:common-weights} produces the term $-m\mu_+\mathcal{F}_+^2/2$; for $m$ large enough, Lemma~\ref{lem:3d-acoustic-completion} gives
	\[
	\mathcal{Q}_+-\frac{m}{2}\mathcal{F}_+^2
	\le(\beta+\epsilon)s_+^2-c_\epsilon\mathcal{F}_+^2.
	\]
	We apply this inequality to a strictly increasing scalar reference profile $\bar{s}$, without assuming monotonicity of $s_y$.
	
	To construct $\bar{s}$, Proposition~\ref{prop:outgoing-reduction} introduces an approximate entropy pair $(\widetilde{E},J_0)$, a shifted coordinate $z=Z_+(t,y)$, and a rescaled variable $w(t,Z_+(t,y))=h_0(s(t,y))-C_+(t,y)$ such that the pure $s_y^2$ source term cancels on $r=\sigma=0$. Because of this cancellation, the source term in the equation for $q_+=w-\Omega$, where $\Omega$ solves the scalar viscous equation \eqref{eq:out-reference-new}, is integrable in $L^1(\mathbb{R})$ with weight $T^{1/20}$ (Lemma~\ref{lem:outgoing-source}). The Nash estimate in Lemma~\ref{lem:outgoing-nash} then gives a small logarithmic bound for $\int_0^t\|q_+(\tau)\|_{L^\infty}^2\dd\tau$. Setting $\bar{s}(t,y)=h_0^{-1}(\Omega(t,y))$, Proposition~\ref{prop:outgoing-amplitude} yields
	\[
	\left\{
	\begin{aligned}
		&0\le\partial_1\lambda_+^E\bigl(\bar{s}(t,Y)\bigr)
		\le\frac{1+\varepsilon_++o_\delta(1)}{T}+k_+(t),
		\qquad k_+\in L^1(0,\infty),\\
		&\int_0^t\|s(\tau)-\bar{s}(\tau)\|_{L^\infty}^2\dd\tau
		\le C_\delta+\varepsilon_\delta\log(1+t),
		\qquad \lim_{\delta\rightarrow 0}\varepsilon_\delta=0.
	\end{aligned}
	\right.
	\]
	In the exterior energy estimate, integrating $\partial_1(s-\bar{s})$ by parts once puts the derivative onto $\nabla h$, and Young's inequality bounds the result by
	\[
	\epsilon\mathcal{H}
	+C_\epsilon\bigl(\|s-\bar{s}\|_{L^\infty}^2+\mathcal{K}(t)\bigr)E_{-1},
	\qquad \mathcal{K}\in L^1(0,\infty),
	\]
	which adds at most $C_\epsilon\varepsilon_\delta\log T$ to the integrated energy growth.
	
	\subsubsection{The backward acoustic slope and monotone phase}
	\label{intro:backward-estimate}
	Although $r$ vanishes on the profile $V(W)$ and $r_y$ changes sign, the term involving $s_-^2$ requires only an upper bound on $r_y$. Using Lemma~\ref{inc:reduction}, a first-order correction $\tilde{r}=r-kU_y$ (with $U=(v,u,\sigma)^\top$) and normalization $w_-=e^{h_-}\tilde{r}_y$ reduce the backward slope equation to the Riccati equation
	\[
	(w_-)_t+\bigl(b_--(d_-)_y+2d_-(h_-)_y\bigr)(w_-)_y
	=d_-(w_-)_{yy}-\kappa_0(p)w_-^2+\epsilon_1w_-+e^{h_-}R_y,
	\]
	where $\int_0^\infty T\|R_y(t)\|_{L^\infty}\dd t<\infty$ by the central acoustic and exterior entropy estimates. 
	The maximum principle for $w_-^+$ and the relation \(\alpha=a_b(s)\tilde r\) then give the one-sided upper bound \eqref{inc:regular-slope}:
	\[
	\kappa_*\alpha_y\le\frac{1+\varepsilon_s+o_\delta(1)}{T}+k_s(t),
	\qquad k_s\in L^1(0,\infty).
	\]
	Since \(r_y=a_b(s)^{-1}\alpha_y+O(T^{-\frac{36}{25}})\), this estimate also gives a one-sided upper bound on \(r_y\) and hence on the coefficient of \(s_-^2\) in the three-dimensional energy estimate.
	
	Moreover, by \eqref{r5a:acoustic-projections}, the projections $s_-$ and $s_+$ satisfy
	\[
	s_-^2+s_+^2=\hat{a}_1^2+v_{11}^2\le\mathcal{J}
	\qquad\text{on the exterior region where }K\ge0.
	\]
	Therefore, the upper bounds $(1+\varepsilon_s)/T$ {for the coefficient of $s_-^2$ and $(1+\varepsilon_+)/T$ for that of $s_+^2$} do not add to $2/T$; their sum is bounded by the larger coefficient times $s_-^2+s_+^2\le\mathcal{J}$. Together with $E_e(t)\le(1+\epsilon_0)E_{-1}(t)$ from \eqref{r5a:fan-additivity}, this gives a single leading growth coefficient:
	\[
	\frac{1+\varepsilon_s+C\varepsilon_++C\epsilon+o_\delta(1)}{T}
	\int_{\mathbb{R}^3}\frac{\zeta_e\mathcal{W}}{2}(s_-^2+s_+^2)\dd x
	\le\frac{(1+\epsilon_0)(1+\varepsilon_s+C\varepsilon_++C\epsilon+o_\delta(1))}{T}E_{-1}(t).
	\]
	
	To control the remaining mixed backward terms involving $|\partial_1\alpha_E|\mathcal{F}_-$ (with $\alpha_E$ obtained from the mollified, truncated \(\alpha\) in \eqref{111}), Section~\ref{sec:actual-incoming-phase} constructs a nonnegative comparison solution $m_-(t,y)\ge|\alpha(t,y)|$ with $\|m_-(t)\|_{L^1}\le CT^{1/25}$ and defines the phase slope
	\[
	\psi=a_1(H_\ell*m_-)>0,\qquad
	a_1=\eta_0T^{-1/4},\qquad \ell=T^{2/5}.
	\]
	After flattening in the central region and changing to Eulerian coordinates, $\psi_E=\partial_1\Psi_E\ge0$ is positive on $\operatorname{supp}\alpha_E$ and satisfies (Proposition~\ref{prop:positive-phase})
	\[
	I_E(t)=\int_{\operatorname{supp}\alpha_E}
	\frac{|\partial_1\alpha_E|^2}{\psi_E}\dd x_1
	\le C_\eta T^{-51/100},\qquad
	\int_0^\infty[\log(2+t)]^M I_E(t)^2\dd t<\infty
	\quad(M\ge0).
	\]
	Because the backward wave moves {with speed $(u_b-c_E)-u_b=-c_E$ relative to the fluid}, the factor $e^{\lambda_h\Psi_E}$ in $\mathcal{W}$ produces the damping term
	\[
	-\frac{\lambda_h}{2}\int_{\mathbb{R}^3}\zeta_e\mathcal{W}c_E\psi_E\mathcal{F}_-^2\dd x
	\le-c_T\mathcal{T}_\rho,
	\]
	which controls the mixed backward terms by $\epsilon(\mathcal{H}+\mathcal{T}_\rho)+C_\epsilon(I_E^2+\mathcal{K})E_{-1}$ in Lemma~\ref{lem:3d-signed-incoming}.
	
	\subsubsection{Entropy decomposition and elliptic coercivity}
	\label{intro:entropy-estimate}
	In Proposition~\ref{prop:entropy-decomposition}, we split the relative entropy as $\sigma=H_P+\sigma_0$. Here $H_P\ge0$ solves the heat equation forced by the Burgers dissipation $P=K(W)W_x^2$; {its derivative is bounded by $C\delta T^{-1}\log(2+t)+k_P(t)$} for $|Y|\le T^{13/25}$ and is integrable in time outside that region. Taking $w\approx[\log(2+t)]^{-2}$ in the central region gives
	\[
	\|\sqrt{w}\,\partial_1H_P(t)\|_{L^\infty}
	\le\frac{C\delta}{T}+k_C(t),\qquad k_C\in L^1(0,\infty),
	\]
	while $\sigma_0$ satisfies $\int_0^t\|\sigma_0(\tau)\|_{L^\infty}^2\dd\tau\le C_\delta+\varepsilon_\delta\log T$ and is estimated by integration by parts in \eqref{r5a:contact-X-allocation}.
	
	Let $x_b(t,\cdot)=Y(t,\cdot)^{-1}$ be the inverse mass-coordinate map and set
	\begin{equation}\label{intro:contact-trajectory}
		x_c(t)=x_b(t,0),\qquad
		Y\bigl(t,x_c(t)\bigr)=0,\qquad
		x_c'(t)=u_b\bigl(t,x_c(t)\bigr),\qquad x_c(0)=0.
	\end{equation}
	The region where $K=w-1<0$ is contained in
	\[
	I(t)=\bigl\{|x_1-x_c(t)|\le L(t)\bigr\},\qquad
	L(t)\lesssim T^{13/25}[\log(2+t)]^{1/p_c}.
	\]
	We take the weight $\mathcal{W}$ to be constant on the wider interval
	\[
	J_0(t)=\bigl\{|x_1-x_c(t)|\le\ell_0(t)\bigr\},\qquad
	\ell_0(t)\approx T^{53/100},\qquad \lim_{t\to\infty}\frac{L(t)}{\ell_0(t)}=0.
	\]
	In Lemma~\ref{lem:3d-coercivity}, we combine the identity $\int|\nabla v|^2=\int(|\Div v|^2+|\operatorname{curl}v|^2)$ where the weight is constant with an off-diagonal bound of order $(L/\ell_0)^{1/2}\ll w_{\min}$ for the operator $\mathbf{1}_{J_0^c}\nabla\operatorname{curl}(-\Delta)^{-1}\mathbf{1}_I$, which yields the coercivity estimate \eqref{r5a:strong-first}.
	
	\subsubsection{Three-dimensional decay}
	Combining the above estimates gives Proposition~\ref{prop:3d-common-growth} with $A_*\le 203/200<2$. Two-dimensional Fourier splitting in $x_\perp$ then gives the $H^3(\mathbb{R}^3)$ bound
	\[
	\|Z(t)\|_{H^3(\mathbb{R}^3)}\le B_0(t)=C_\delta\varepsilon T^{1/50}.
	\]
	Inserting this bound into the derivative energy inequality for $E_1\approx\|\nabla Z\|_{H^2(\mathbb{R}^3)}^2$ and using the interpolation $E_1^2\le CB_0^2D_1$ {gives the following inequality for $\mathcal{Y}=E_1/B_0^2$ when $\varepsilon>0$ (for $\varepsilon=0$, uniqueness gives $U=U_b$):}
	\[
	\mathcal{Y}'+\frac{1}{25T}\mathcal{Y}+c\mathcal{Y}^2\le CK_b(t).
	\]
	{For $t\ge1$, comparison on $[t/2,t]$ together with} $\int_{t/2}^tK_b(\tau)\dd\tau\le C_\delta T^{-23/25}$ gives $\mathcal{Y}(t)\le C_\delta T^{-23/25}$, and hence
	\[
	\|Z(t)\|_{L^\infty(\mathbb{R}^3)}+\|\nabla Z(t)\|_{H^2(\mathbb{R}^3)}
	\le C_\delta\varepsilon T^{-11/25}.
	\]
	
	The rest of the paper is organized as follows. Section~\ref{sec:planar-decay-all} constructs the planar wave $U_b$ and proves its decay in $L^2(\mathbb{R})$ and $L^\infty(\mathbb{R})$. Section~\ref{sec:entropy-central-all} proves the entropy decomposition and the central acoustic estimates. Sections~\ref{sec:backward-all} and \ref{sec:forward-all} establish the backward and forward acoustic estimates. Section~\ref{sec:three-dimensional-all} proves the coercivity of the potential energy and the combined energy inequality. Section~\ref{sec:global-all} completes the local continuation and three-dimensional decay proof.
	
	\subsection{Notation and thermodynamic coordinates}
	
	Throughout the paper, we set
	\[
	T=1+t,\qquad J=\log(2+t),\qquad a=\frac{12}{25}.
	\]
	Unless otherwise stated, $L^p$ and $H^k$ norms are taken over $\mathbb{R}$ in the one-dimensional analysis (Sections~\ref{sec:elementaryplanarinputs}--\ref{sec:outgoing-energy-replacement}) and over $\mathbb{R}^3$ in the three-dimensional analysis (Sections~\ref{sec:r5a-full-interface}--\ref{sec:main-proof}). We write $\partial_1=\partial_{x_1}$, and in one-dimensional calculations we write $x$ for $x_1$. {In the intermediate estimates, constants $C>0$ may depend on the fixed physical parameters, $U_-$, $e$, $m_0$, and the derivative order}, whereas $C_\delta$ may also depend on $\delta$. The notation $o_\delta(1)$ stands for a quantity tending to zero as $\delta\rightarrow 0$ after all fixed parameters have been chosen.
	
	In the planar analysis, $W(t,x)$ denotes the solution of the viscous Burgers equation \eqref{planar:viscous-reference}, and $W^E(t,x)$ denotes the inviscid Burgers profile. The rarefaction profile map is written as $V(W)$ in primitive variables and $\mathbf{V}(W)$ in conserved variables, so that the planar perturbation $g$ and three-dimensional perturbation $f$ are
	\[
	g=\mathbf{U}_b-\mathbf{V}(W),\qquad f=\mathbf{U}-\mathbf{U}_b.
	\]
	
	To write the planar equations in Lagrangian mass coordinates, we set $v_-=\rho_-^{-1}$ and $p_-=R_g\rho_-\theta_-$, and define
	\begin{equation}\label{v3:I-explicit}
		I(p)=\frac{2c_-}{\gamma-1}\left[\left(\frac{p}{p_-}\right)^{(\gamma-1)/(2\gamma)}-1\right],\qquad p>0.
	\end{equation}
	For a planar state with specific volume $v=1/\rho$, velocity $u$, and specific entropy $S=c_v\log\theta+R_g\log v$, we define
	\begin{equation}\label{v3:acoustic-variables}
		r=u-I(p),\qquad s=u+I(p),\qquad \sigma=S-S_-.
	\end{equation}
	Differentiating \eqref{v3:I-explicit} and using $p_-/c_-=\sqrt{p_-/(\gamma v_-)}$ gives
	\begin{equation}\label{v3:coefficient-functions}
		c_0(p)=\frac{1}{I'(p)}
		=\sqrt{\frac{\gamma p_-}{v_-}}\left(\frac{p}{p_-}\right)^{(\gamma+1)/(2\gamma)},
		\qquad
		v_0(p)=\frac{\gamma p}{c_0(p)^2}
		=v_-\left(\frac{p_-}{p}\right)^{1/\gamma}.
	\end{equation}
	Substituting $\theta=pv/R_g$ into the entropy formula yields
	\[
	\sigma=S-S_-=c_v\log\left(\frac{p}{p_-}\right)+c_p\log\left(\frac{v}{v_-}\right),
	\]
	which gives $v$, the mass-coordinate sound speed $c_L$, and the Eulerian sound speed $c_E$:
	\begin{equation}\label{v3:state-recovery}
		v=v_0(p)e^{\sigma/c_p},\qquad
		c_L=\sqrt{\frac{\gamma p}{v}}=c_0(p)e^{-\sigma/(2c_p)},\qquad
		c_E=vc_L=\sqrt{\gamma R_g\theta}.
	\end{equation}
	The state variables are recovered from $(r,s,\sigma)$ through
	\[
	u=\frac{r+s}{2},\qquad p=I^{-1}\left(\frac{s-r}{2}\right),\qquad
	v=v_0(p)e^{\sigma/c_p},\qquad \theta=\frac{pv}{R_g}.
	\]
	We set $c_*=c_0(p_-)$.
	
	Writing $V(W)=(\rho_v,u_v,\theta_v)$ with $v_v=\rho_v^{-1}$ and $p_v=R_g\rho_v\theta_v$, direct substitution of the state formulas from Section~\ref{sec:main-results} gives
	\begin{equation}\label{v3:profile-identities}
		\theta_vv_v^{\gamma-1}=\theta_-v_-^{\gamma-1},\qquad
		I(p_v)=u_v,\qquad r\bigl(V(W)\bigr)=\sigma\bigl(V(W)\bigr)=0.
	\end{equation}
	Thus both $r$ and $\sigma$ vanish identically on $V(W)$. {Since $a_0=(c_0(p)/2)(e^{-\sigma/c_p}-1)$ in \eqref{v3:coupling-bound},} the coefficient $a_0$ satisfies $|a_0|\le C|\sigma|$, which is used in \eqref{ce:7} and in Section~\ref{sec:global-slope-new}.
	
	\section{Construction and decay of the planar wave}
	\label{sec:planar-decay-all}
	
	\subsection{Planar energy and derivative estimates}
	\label{sec:elementaryplanarinputs}
	
	In this section, we construct the one-dimensional planar rarefaction wave $U_b(t,x_1)$ and prove uniform Sobolev bounds and conditional time-decay estimates on any finite time interval. Suppressing the vanishing transverse velocity components in one dimension, we write
	\[
	\mathbf{U}=(\rho,\rho u,E)^\top,\qquad
	\theta=\frac{E/\rho-u^2/2}{c_v},\qquad p=R_g\rho\theta,
	\]
	so that the one-dimensional Navier--Stokes--Fourier system reads
	\begin{equation}\label{planar:conserved-equation}
		\mathbf{U}_t+F(\mathbf{U})_x=\bigl(B(\mathbf{U})\mathbf{U}_x\bigr)_x,
	\end{equation}
	where
	\[
	F(\mathbf{U})=
	\begin{pmatrix}\rho u\\\rho u^2+p\\u(E+p)\end{pmatrix},
	\]
	\[
	B(\mathbf{U})=\begin{pmatrix}
		0&0&0\\
		-\nu u/\rho&\nu/\rho&0\\
		-\nu u^2/\rho+\kappa(u^2/2-c_v\theta)/(\rho c_v)
		&\nu u/\rho-\kappa u/(\rho c_v)&\kappa/(\rho c_v)
	\end{pmatrix}.
	\]
	Let $\mathbf{V}(W)$ be the conserved approximate rarefaction profile defined in Section~\ref{sec:main-results}, with primitive components $(\rho_v,u_v,\theta_v)=V(W)$. The viscous Burgers profile $W(t,x)$ solves the Cauchy problem
	\begin{equation}\label{planar:viscous-reference}
		\left\{
		\begin{aligned}
			&W_t+WW_x=d_*W_{xx},\\
			&W(0,x)=W^E(0,x),\qquad d_*=\frac{\nu+(\gamma-1)\kappa/c_p}{2\rho_-}>0.
		\end{aligned}
		\right.
	\end{equation}
	Throughout Section~\ref{sec:planar-decay-all}, $x$ denotes $x_1$ and all norms are taken over $\mathbb{R}$. We write the conserved and primitive planar perturbations as
	\[
	g=\mathbf{U}_b-\mathbf{V}(W),\qquad
	z=(\phi,\psi,\zeta)=(\rho_b-\rho_v,u_b-u_v,\theta_b-\theta_v).
	\]
	In the energy estimates below, we write $(\rho,u,\theta)=U_b$ and use the subscript $k$ in $g_k,z_k,\phi_k,\psi_k,\zeta_k$ to denote $\partial_x^k$. {On the compact set of positive states, the product estimates and the profile derivative bounds in Lemma~\ref{lem:planar-burgers} give equivalent Sobolev norms for the small perturbations $g$ and $z$.} We fix $N=22$ and $a=12/25$.
	
	\subsubsection{The viscous Burgers profile}
	
	\begin{lemma}
		\label{lem:planar-burgers}
		The solution $W(t,x)$ of \eqref{planar:viscous-reference} is strictly increasing in $x$ and satisfies, for all $t\ge 0$,
		\[
		\|W_x(t)\|_{L^1}=\delta,\qquad
		\|W_x(t)\|_{L^\infty}+\|W_{xx}(t)\|_{L^1}\le C\min\{\delta,T^{-1}\},
		\]
		together with the higher-derivative bounds
		\[
		\|\partial_x^jW(t)\|_{L^2}^2\le C_jT^{-j-1/2}\quad(j\ge2),
		\qquad
		\|\partial_x^jW(t)\|_{L^\infty}\le C_jT^{-(j+1)/2}\quad(j\ge1).
		\]
		Here $C,C_j$ are uniform for small $\delta>0$ at fixed $e>0$, and every fixed-order spatial derivative of $W$ tends uniformly to zero as $\delta\rightarrow0$. In particular,
		\[
		\int_0^\infty\|\mathbf{V}(W)_{xx}(t)\|_{L^2}\dd t\le C\delta^{1/5}.
		\]
	\end{lemma}
	
	\begin{proof}
		Let $w_0$ denote the initial datum in \eqref{main:burgers-data} and set $\bar{W}=(W_-+W_+)/2$. By definition, $W(0,x)=W^E(0,x)=w^E(1,x)$, which is given along characteristics by
		\[
		x=\xi+w_0(\xi),\qquad
		W(0,x)=w_0(\xi),\qquad
		W_x(0,x)=\frac{w_0'(\xi)}{1+w_0'(\xi)}.
		\]
		The shifted coordinate $x-\bar{W}=\xi+(w_0(\xi)-\bar{W})$ is an odd, strictly increasing function of $\xi$. Because $w_0'(\xi)$ is even, positive, and strictly decreasing in $|\xi|$, the initial slope $W_x(0,x)$ is positive, even, and has a single maximum at $x=\bar{W}$, and all of its fixed-order derivatives are $O(\delta)$. The viscous Burgers equation preserves the positivity of $W_x$ and the odd symmetry of $W(t,x)-\bar{W}$ about $x_B(t)=\bar{W}(1+t)$. Applying the maximum principle to
		\[
		W_{xt}+WW_{xx}=d_*W_{xxx}-(W_x)^2
		\]
		gives
		\[
		\|W_x(t)\|_{L^\infty}\le
		\frac{\|W_x(0)\|_{L^\infty}}{1+t\|W_x(0)\|_{L^\infty}}
		\le C\min\{\delta,T^{-1}\},\qquad \|W_x(t)\|_{L^1}=\delta.
		\]
		In the moving coordinate $\xi=x-x_B(t)$, $W_{xx}$ is odd in $\xi$, vanishes at $\xi=0$, and is nonpositive at $t=0$ on $\xi>0$. By the maximum principle on $(0,\infty)$, $W_{xx}(t,\cdot)$ remains nonpositive for $\xi>0$, so $W_x(t,\cdot)$ has a single maximum at $\xi=0$ and satisfies
		\[
		\|W_{xx}(t)\|_{L^1}=2\|W_x(t)\|_{L^\infty}\le C\min\{\delta,T^{-1}\}.
		\]
		
		Next, $W_{xx}$ and its backward-in-time dual variable $v(s,x)$ satisfy
		\[
		\left\{
		\begin{aligned}
			&W_{xxt}+(WW_{xx})_x=d_*W_{xxxx}-2W_xW_{xx},\\
			&{-v_s=d_*v_{xx}+Wv_x-2W_xv}.
		\end{aligned}
		\right.
		\]
		Since $W_x\ge 0$, the forward equation is $L^1$-contractive and has a nonpositive zeroth-order term in its $L^2$ energy identity; likewise, the {dual equation in reverse time has nonpositive $L^1$ derivative} $-3\int_{\mathbb{R}} W_xv\dd x$ for $v\ge 0$ and the same nonpositive zeroth-order term in its $L^2$ energy. Combining the $L^2$ energy inequality with the one-dimensional Nash inequality gives an $L^1\to L^2$ bound with factor $(t-s)^{-1/4}$ for the forward equation, and duality gives an $L^2\to L^\infty$ bound with the same factor $(t-s)^{-1/4}$. Applying these bounds on $[t/2,t]$ yields
		\begin{equation}\label{af:Burgers-h}
			\|W_{xx}(t)\|_{L^2}\le CT^{-5/4},\qquad
			\|W_{xx}(t)\|_{L^\infty}\le CT^{-3/2}.
		\end{equation}
		
		To bound the higher derivatives, set $B_k(t)=\|\partial_x^kW_{xx}(t)\|_{L^2}^2$ for $k\ge 1$. Differentiating $W_{xxt}+WW_{xxx}=d_*W_{xxxx}-3W_xW_{xx}$ $k$ times, using $W_x\ge 0$, and applying the Leibniz rule gives, for $k\ge 1$,
		\[
		B_k'+cB_{k+1}\le CT^{-1}B_k+
		C\|W_{xx}\|_{L^\infty}B_{k-1}^{1/2}B_k^{1/2}.
		\]
		Indeed, every nonlinear commutator beyond the $W_x$ term is a sum of products $\partial_x^jW_{xx}\,\partial_x^{k-1-j}W_{xx}$ ($0\le j\le k-1$), whose $L^2$ norm is bounded by Gagliardo--Nirenberg interpolation:
		\[
		\bigl\|\partial_x^jW_{xx}\,\partial_x^{k-1-j}W_{xx}\bigr\|_{L^2}
		\le C\|W_{xx}\|_{L^\infty}\|\partial_x^{k-1}W_{xx}\|_{L^2}.
		\]
		Assuming inductively that $B_{k-1}(t)\le C_{k-1}T^{-k-3/2}$, Young's inequality and \eqref{af:Burgers-h} give
		\[
		B_k'+cB_{k+1}\le CT^{-1}B_k+CT^{-k-7/2}.
		\]
		By interpolation between adjacent derivatives,
		\[
		B_{k+1}\ge \frac{B_k^2}{B_{k-1}}\ge cT^{k+3/2}B_k^2.
		\]
		Comparing $B_k(t)$ with $H_kT^{-k-5/2}$ for large $H_k>0$ closes the induction: at any first crossing time, the quadratic term $-cH_k^2T^{-k-7/2}$, the time derivative of $H_kT^{-k-5/2}$, and the right-hand side all have power $T^{-k-7/2}$, and the quadratic term dominates. Induction and interpolation then yield
		\begin{equation}\label{af:Burgers-all}
			\|\partial_x^jW(t)\|_{L^2}^2\le C_jT^{-j-1/2}\quad(j\ge2),
			\qquad
			\|\partial_x^jW(t)\|_{L^\infty}\le C_jT^{-(j+1)/2}\quad(j\ge1).
		\end{equation}
		
		In addition, all of these norms tend to zero as $\delta\rightarrow 0$. For $k=0$, the $L^2$ norms of $W_x$ and $W_{xx}$ are nonincreasing in time. {For $k\ge1$, finite-time continuous dependence gives $B_k=o_\delta(1)$ uniformly on each fixed time interval, while \eqref{af:Burgers-all} gives $B_k(t)\le C_kT^{-k-5/2}$ uniformly on its complement. Taking first a sufficiently long interval and then $\delta$ small gives $\sup_{t\ge0}B_k(t)=o_\delta(1)$.} By the chain and product rules, $\mathbf{V}(W)$ satisfies
		\[
		\|\mathbf{V}(W)_{xx}(t)\|_{L^2}\le C\min\{\delta,T^{-5/4}\},\qquad
		\int_0^\infty\|\mathbf{V}(W)_{xx}(t)\|_{L^2}\dd t\le C\delta^{1/5},
		\]
		while $(W_x)^2$ decays faster in $L^2$. Applying the exterior interpolation of Section~\ref{sec:exterior-L1} to $W$ gives the higher-derivative $L^1$ estimates used in \eqref{af:puresource}.
	\end{proof}
	
	\subsubsection{Energy bounds and conditional decay}
	
	\begin{proposition}
		\label{prop:planar-energy}
		For sufficiently small $\delta>0$, the one-dimensional planar solution $U_b$ with initial datum \eqref{main:planar-data} exists globally in time and remains in a compact subset of positive states. There exists $A_{\delta,N}$ with $A_{\delta,N}\to 0$ as $\delta\rightarrow0$ such that
		\[
		\sup_{t\ge0}\|g(t)\|_{H^N}^2+
		\int_0^\infty\left(
		\|\phi_x(t)\|_{H^{N-1}}^2+
		\|\bigl(\psi_x(t),\zeta_x(t)\bigr)\|_{H^N}^2+
		\int_{\mathbb{R}}W_x|g|^2\dd x\right)\dd t
		\le A_{\delta,N}.
		\]
		Furthermore, on any finite time interval $[0,t_*]$ on which $\|g(t)\|_{L^2}^2\le T^{-a}$ holds with $a=12/25$, the derivatives of $g$ satisfy
		\begin{equation}\label{af:bootstrap}
			\|\partial_x^kg(t)\|_{L^2}\le C_kT^{-(k+1/2)a}
			\quad(1\le k\le21),\qquad
			\|\partial_x^kg(t)\|_{L^\infty}\le C_kT^{-(k+1)a}
			\quad(0\le k\le20),
		\end{equation}
		where the constants $C_k$ are independent of $t_*$ and $\delta$.
	\end{proposition}
	
	\begin{proof}
		We divide the proof into four steps.
		
		\medskip
		\noindent\textit{Step 1: Relative entropy and low-order energy.}
		Since $V(W)=(\rho_v,u_v,\theta_v)$ solves the one-dimensional Euler equations up to the diffusion terms from $d_*W_{xx}$, the mass, momentum, and temperature perturbation equations contain the profile terms
		\[
		-d_*\rho_v'(W)W_{xx},\qquad -d_*u_v'(W)W_{xx},\qquad -d_*\theta_v'(W)W_{xx},
		\]
		together with the viscous profile terms
		\[
		\frac{\nu}{\rho_v}u_{v,xx},\qquad
		\frac{\kappa}{\rho_v c_v}\theta_{v,xx},\qquad
		\frac{\nu}{\rho_v c_v}u_{v,x}^2.
		\]
		Replacing $\rho_v$ in these denominators by $\rho$ {introduces products of $g$ with second derivatives of the profile and the term $-\nu\phi u_{v,x}^2/(c_v\rho\rho_v)$}. Hence the profile residual $\mathcal{R}_0$ satisfies, by Lemma~\ref{lem:planar-burgers},
		\[
		\|\mathcal{R}_0(t)\|_{L^2}\le CT^{-5/4},\qquad
		\|\mathcal{R}_0(t)\|_{H^2}^2\le CT^{-5/2},
		\]
		while all terms containing $W_x^2$ decay faster in $L^2(\mathbb{R})$.
		
		Define the relative entropy functional
		\[
		E_0=\int_{\mathbb{R}}\left(R_g\rho\theta_v\Phi_0\left(\frac{\rho_v}{\rho}\right)
		+c_v\rho\theta_v\Phi_0\left(\frac{\theta}{\theta_v}\right)+\frac{1}{2}\rho\psi^2\right)\dd x,
		\qquad\Phi_0(s)=s-1-\log s,
		\]
		along with the first- and second-order derivative energies and the Kawashima cross terms:
		\[
		E_k=\frac{1}{2}\int_{\mathbb{R}}\left(\frac{R_g}{\rho}\phi_k^2
		+\frac{\rho}{\theta}\psi_k^2
		+\frac{\rho c_v}{\theta^2}\zeta_k^2\right)\dd x,\qquad
		C_k=\int_{\mathbb{R}}\psi_{k-1}\phi_k\dd x,\qquad k=1,2.
		\]
		For small fixed constants $0<\varepsilon_2\ll\varepsilon_1\ll 1$, the combined energy
		\[
		\mathcal{E}=E_0+\varepsilon_1(E_1+E_2)+\varepsilon_2(C_1+C_2)
		\]
		is uniformly equivalent to $\|g\|_{H^2}^2$, and its dissipation controls
		\[
		\mathcal{D}=\|\phi_x\|_{H^1}^2+\|(\psi_x,\zeta_x)\|_{H^2}^2,
		\qquad P_0=\int_{\mathbb{R}}W_x|g|^2\dd x.
		\]
		To check the positivity of the $W_x$ term in ${-\frac{\dd}{\dd t}E_0}$, set
		\[
		X_1=\sqrt{\frac{R_g\theta_v}{\rho_v}}\,\phi,\qquad
		X_2=\sqrt{\rho_v}\,\psi,\qquad
		X_3=\sqrt{\frac{\rho_v c_v}{\theta_v}}\,\zeta.
		\]
		Dividing the quadratic form in ${-\frac{\dd}{\dd t}E_0}$ by $u_{v,x}>0$ gives
		\[
		\begin{aligned}
			&X_2^2+\frac{\gamma-1}{2}(X_1^2+X_3^2)
			-\frac{\gamma-1}{\sqrt{\gamma}}X_1X_2
			+\sqrt{\frac{\gamma-1}{\gamma}}X_2X_3\\
			&\quad=\frac{\gamma-1}{2}\left(X_1-\frac{X_2}{\sqrt{\gamma}}\right)^2
			+\frac{\gamma-1}{2}\left(X_3+\frac{X_2}{\sqrt{(\gamma-1)\gamma}}\right)^2
			+\frac{1}{2}X_2^2,
		\end{aligned}
		\]
		and the cubic Taylor remainder $O(|g|^3)W_x$ is absorbed by $P_0$ for small $\|g\|_{L^\infty}$. {For the second density derivative, Young's inequality gives
			\[
			\left|d_*\int_{\mathbb R}\frac{R_g}{\rho}\phi_{xx}\,
			\partial_x^2\bigl(\rho_v'(W)W_{xx}\bigr)\dd x\right|
			\le \epsilon\|\phi_{xx}\|_{L^2}^2+
			C_\epsilon\|\mathcal R_0\|_{H^2}^2 .
			\]
			In the remaining residual pairings we integrate by parts only when the resulting derivative is controlled by $\mathcal D$. }Combining the derivative energy estimates with $C_1,C_2$ and absorbing cubic terms via $\|g\|_{H^N}\ll 1$, we obtain
		\begin{equation}\label{af:physical-energy}
			\mathcal{E}'+c(\mathcal{D}+P_0)
			\le k_0(t)\mathcal{E}+
			C\|\mathcal{R}_0(t)\|_{L^2}\|g(t)\|_{L^2}+CQ_0(t),
		\end{equation}
		where
		\[
		{Q_0(t)\le CT\|\mathcal{R}_0(t)\|_{H^2}^2\le CT^{-3/2}},\qquad
		k_0(t)\le C\bigl(\|W_{xx}(t)\|_{L^\infty}+\|W_x(t)\|_{L^\infty}^2+\|g_x(t)\|_{L^\infty}^2\bigr).
		\]
		
		\medskip
		\noindent\textit{Step 2: Global existence and uniform smallness in $H^N(\mathbb{R})$.}
		By the one-dimensional Gagliardo--Nirenberg inequality,
		\[
		\|g_x\|_{L^\infty}^2\mathcal{E}
		{\le C\mathcal E\bigl(\|z_x\|_{L^2}\|z_{xx}\|_{L^2}
			+\|W_x\|_{L^\infty}^2\|z\|_{L^\infty}^2\bigr)
			\le C\mathcal E\mathcal D+C\|W_x\|_{L^\infty}^2\mathcal E,}
		\]
		{by the chain rule for $g$ and the smallness of $z$. The first term is absorbed by $c\mathcal D$ when $\sup_{0\le s\le t}\mathcal E(s)$ is small.} Since $g(0)=0$ and
		\[
		\int_0^\infty\bigl(\|W_{xx}\|_{L^\infty}+\|W_x\|_{L^\infty}^2\bigr)\dd t\le C,\qquad
		\int_0^\infty\bigl(\|\mathcal{R}_0(t)\|_{L^2}+Q_0(t)\bigr)\dd t=o_\delta(1),
		\]
		Gronwall's inequality gives
		\begin{equation}\label{af:direct-classical}
			\sup_{t\ge 0}\|g(t)\|_{H^2}^2+
			\int_0^\infty\bigl(\mathcal{D}(t)+P_0(t)\bigr)\dd t\le A_\delta,
			\qquad \lim_{\delta\rightarrow 0}A_\delta=0.
		\end{equation}
		
		To extend this bound to order $N=22$, we multiply the momentum and temperature equations by $\rho$ and $\rho c_v$ so that the second-order viscosity and heat-conduction coefficients are constant. Pairing the $k$-th derivatives with
		\[
		\left(\frac{R_g\theta}{\rho}\phi_k,\ \psi_k,\ \frac{\zeta_k}{\theta}\right),
		\]
		adding $\eta\int_{\mathbb{R}}\psi_{k-1}\phi_k\dd x$, and estimating commutators by the Moser inequality
		\[
		\bigl\|[\partial_x^k,a]b\bigr\|_{L^2}
		\le C_k\bigl(\|a_x\|_{L^\infty}\|\partial_x^{k-1}b\|_{L^2}
		+\|\partial_x^ka\|_{L^2}\|b\|_{L^\infty}\bigr),
		\]
		{we integrate the viscous residuals by parts once and bound the top-order density residual directly by $\epsilon\|\phi_N\|_{L^2}^2+C_\epsilon\|\partial_x^N(\rho_v^\prime(W)W_{xx})\|_{L^2}^2$. Lemma~\ref{lem:planar-burgers} then gives}
		\[
		\sup_{t\ge 0}\|g(t)\|_{H^N}^2+
		\int_0^\infty\left(\|\phi_x\|_{H^{N-1}}^2+
		\|(\psi_x,\zeta_x)\|_{H^N}^2+P_0\right)\dd t
		\le A_{\delta,N}\longrightarrow 0\quad(\delta\rightarrow 0).
		\]
		Combined with local existence, this yields the global smooth solution $U_b$ for $N=22$.
		
		\medskip
		\noindent\textit{Step 3: High-order derivative energy inequalities.}
		For $1\le k<N$, define
		\[
		\left\{
		\begin{aligned}
			&\mathcal{E}_k=\frac{1}{2}\sum_{j=k}^N\int_{\mathbb{R}}\left(
			\frac{R_g}{\rho}\phi_j^2+\frac{\rho}{\theta}\psi_j^2+
			\frac{\rho c_v}{\theta^2}\zeta_j^2\right)\dd x
			+\eta\sum_{j=k}^{N-1}\int_{\mathbb{R}}\psi_j\phi_{j+1}\dd x,\\
			&\mathcal{D}_k=\sum_{j=k+1}^N\|\phi_j\|_{L^2}^2+
			\sum_{j=k+1}^{N+1}\|(\psi_j,\zeta_j)\|_{L^2}^2,
		\end{aligned}
		\right.
		\]
		with $\eta>0$ fixed small. We show that
		\begin{equation}\label{af:derivative-energy}
			\mathcal{E}_k'+c\mathcal{D}_k
			\le C\|g_x\|_{L^\infty}\|g_k\|_{L^2}^2+CT^{-1}\mathcal{E}_k
			+CT^{-k-1}\|g\|_{L^2}^2+CT^{-k-3/2}.
		\end{equation}
		
		Write $z=(\phi,\psi,\zeta)$ and $S=\operatorname{diag}(R_g/\rho,\rho/\theta,\rho c_v/\theta^2)$. By Step 2 and Lemma~\ref{lem:planar-burgers}, for small $\delta>0$ we have
		\[
		\|z(t)\|_{H^N}+\max_{1\le j\le N+2}\|\partial_x^jV(t)\|_{L^\infty}\le\varepsilon_0,
		\qquad N=22,
		\]
		where $\varepsilon_0>0$ can be made arbitrarily small. The perturbation $z=(\phi,\psi,\zeta)$ satisfies
		\[
		\left\{
		\begin{aligned}
			&\phi_t+u\phi_x+\rho\psi_x
			=-u_{v,x}\phi-\rho_{v,x}\psi-d_*\rho_v'(W)W_{xx},\\
			&\psi_t+u\psi_x+\frac{R_g\theta}{\rho}\phi_x+R_g\zeta_x-\frac{\nu}{\rho}\psi_{xx}
			=-u_{v,x}\psi-R_g\left(\frac{\theta}{\rho}-\frac{\theta_v}{\rho_v}\right)\rho_{v,x}
			+\frac{\nu}{\rho}u_{v,xx}-d_*u_v'(W)W_{xx},\\
			&\zeta_t+u\zeta_x+(\gamma-1)\theta\psi_x-\frac{\kappa}{\rho c_v}\zeta_{xx}
			=-\theta_{v,x}\psi-(\gamma-1)u_{v,x}\zeta+\frac{\kappa}{\rho c_v}\theta_{v,xx}\\
			&\hspace{5.5cm}+\frac{\nu}{\rho c_v}(\psi_x+u_{v,x})^2-d_*\theta_v'(W)W_{xx}.
		\end{aligned}
		\right.
		\]
		
		For the nonlinear terms in $z$, set $b_1=\|z_x\|_{L^\infty}$ and $b_2=\|z_{xx}\|_{L^\infty}$. By the Leibniz rule, $\|\partial_x^jF(z)\|_{L^2}\le C_j\|z_j\|_{L^2}$, and Gagliardo--Nirenberg interpolation, we have for $v\in\{\psi,\zeta\}$:
		\begin{align*}
			&\bigl\|[\partial_x^j,F(z)]z_x\bigr\|_{L^2}
			\le C_j b_1\|z_j\|_{L^2},\\
			&\bigl\|\partial_x^j\bigl(F(z)_xz_x\bigr)\bigr\|_{L^2}
			\le C_j b_1\|z_{j+1}\|_{L^2},\\
			&\bigl\|[\partial_x^j,F(z)]v_{xx}\bigr\|_{L^2}
			\le C_j\bigl(b_1\|v_{j+1}\|_{L^2}+b_2\|z_j\|_{L^2}\bigr).
		\end{align*}
		For $k\le j<N$, we write $F(z)v_{xx}=\partial_x(F(z)v_x)-F(z)_xv_x$ and integrate by parts once, bounding the resulting inner products by
		\[
		Cb_1\|z_j\|_{L^2}\|z_{j+1}\|_{L^2}+Cb_1^2\|z_j\|_{L^2}^2.
		\]
		{Writing $\widetilde A$ for the hyperbolic coefficient matrix in the displayed primitive system, the viscous heating and hyperbolic commutators satisfy, with profile terms estimated below,}
		\[
		\bigl\|\partial_x^j\bigl(F(z)\psi_x^2\bigr)\bigr\|_{L^2}
		\le C\bigl(b_1\|\psi_{j+1}\|_{L^2}+b_1^2\|z_j\|_{L^2}\bigr),
		\]
		\[
		\bigl|\bigl\langle[\partial_x^j,S\widetilde{A}]z_x,z_j\bigr\rangle\bigr|
		{\le Cb_1\|z_j\|_{L^2}^2+\text{terms containing profile derivatives}.}
		\]
		At $j=N$, we do not integrate the viscous commutator by parts; instead, the third estimate above bounds its inner product directly by
		\[
		Cb_1\|z_N\|_{L^2}\|(\psi_{N+1},\zeta_{N+1})\|_{L^2}+Cb_2\|z_N\|_{L^2}^2.
		\]
		Since $k<N$, both $\|z_N\|_{L^2}^2$ and $\|(\psi_{N+1},\zeta_{N+1})\|_{L^2}^2$ lie in $\mathcal{D}_k$ with small coefficients $b_1,b_2\le\varepsilon_0$, so $\phi_{N+1}$ does not appear.
		
		For the cross terms $C_j=\int_{\mathbb{R}}\psi_j\phi_{j+1}\dd x$ ($k\le j\le N-1$), differentiating in time gives
		\[
		C_j'=-\int_{\mathbb{R}}\frac{R_g\theta}{\rho}\phi_{j+1}^2\dd x
		+\int_{\mathbb{R}}\rho\psi_{j+1}^2\dd x
		-R_g\int_{\mathbb{R}}\zeta_{j+1}\phi_{j+1}\dd x
		+\int_{\mathbb{R}}\frac{\nu}{\rho}\psi_{j+2}\phi_{j+1}\dd x
		+R_{C,j},
		\]
		where the viscous commutator in $R_{C,j}$ is bounded by
		\[
		C\bigl(b_1\|z_{j+1}\|_{L^2}+b_2\|z_j\|_{L^2}\bigr)\|\phi_{j+1}\|_{L^2}.
		\]
		Since Gagliardo--Nirenberg interpolation gives
		\[
		b_2^2\le C\|z_x\|_{L^\infty}\|z_{xxx}\|_{L^\infty}\le C\varepsilon_0 b_1,
		\]
		Young's inequality absorbs all terms with $j>k$ into $\mathcal{D}_k$ and leaves $Cb_1\|z_k\|_{L^2}^2$.
		
		Next, for the terms mixing $V(W)$ and $z$, Lemma~\ref{lem:planar-burgers} bounds any product of profile derivatives of total order $\ell\ge 1$ in $L^\infty(\mathbb{R})$ by $A_\ell\le C_\ell T^{-(\ell+1)/2}$. The resulting inner products are bounded by sums of
		\[
		A_\ell\|z_{j+1-\ell}\|_{L^2}\|z_j\|_{L^2},\qquad
		A_\ell\|z_{j+2-\ell}\|_{L^2}\|z_j\|_{L^2}.
		\]
		{At $j=N$, the order $N+1$ factor in the second product is $\|(\psi_{N+1},\zeta_{N+1})\|_{L^2}$. When $\ell=1$, these are bounded} by $\epsilon\mathcal{D}_k+CT^{-1}\mathcal{E}_k$. When $\ell\ge 2$ and $j<N$, interpolation between $z$ and $z_{j+1}$ and Young's inequality give
		\[
		C_\epsilon A_\ell^{2(j+1)/(\ell+1)}\|z\|_{L^2}^2
		+C_\epsilon A_\ell^{2(j+1)/\ell}\|z\|_{L^2}^2
		\le CT^{-j-1}\|z\|_{L^2}^2,
		\]
		while at $j=N$ interpolating with $z_N\in\mathcal{D}_k$ gives $CT^{-N}\|z\|_{L^2}^2$. Likewise, the profile terms in $C_j'$ ($k\le j\le N-1$, $\ell\ge 2$) satisfy
		\[
		A_\ell\|z_{j+2-\ell}\|_{L^2}\|z_{j+1}\|_{L^2}
		\le\epsilon\|z_{j+1}\|_{L^2}^2+C_\epsilon A_\ell^{2(j+1)/(\ell-1)}\|z\|_{L^2}^2
		\le\epsilon\mathcal{D}_k+CT^{-j-1}\|z\|_{L^2}^2.
		\]
		Summing over $k\le j\le N$ bounds these mixed terms by $\epsilon\mathcal{D}_k+CT^{-1}\mathcal{E}_k+CT^{-k-1}\|z\|_{L^2}^2$.
		
		Finally, the profile residual is a sum of terms $\partial_x(F_1(W)W_x)+F_2(W)W_x^2$. For $1\le j<N$, integrating by parts once and using Lemma~\ref{lem:planar-burgers} gives
		\[
		\bigl\|\partial_x^j\bigl(F_1(W)W_x\bigr)\bigr\|_{L^2}^2
		+\bigl\|\partial_x^{j-1}\bigl(F_2(W)W_x^2\bigr)\bigr\|_{L^2}^2
		\le CT^{-j-3/2},
		\]
		while at $j=N$ pairing directly with $z_N\in\mathcal{D}_k$ gives $\epsilon\|z_N\|_{L^2}^2+CT^{-N-5/2}$. Combining all three parts yields
		\begin{equation}\label{af:primitive-derivative-energy}
			\mathcal{E}_k'+c\mathcal{D}_k
			\le C\|z_x\|_{L^\infty}\|z_k\|_{L^2}^2
			+CT^{-1}\mathcal{E}_k
			+CT^{-k-1}\|z\|_{L^2}^2+CT^{-k-3/2}.
		\end{equation}
		Since
		\[
		\|z_x\|_{L^\infty}\le C\|g_x\|_{L^\infty}+CT^{-1}\|g\|_{L^\infty},\qquad
		\|z_k\|_{L^2}^2\le C\|g_k\|_{L^2}^2+CT^{-k-1}\|g\|_{L^2}^2,
		\]
		\eqref{af:primitive-derivative-energy} implies \eqref{af:derivative-energy}.
		
		\medskip
		\noindent\textit{Step 4: Closing the derivative decay bounds.}
		Assume that $\|g(t)\|_{L^2}^2\le T^{-a}$ on $[0,t_*]$. By interpolation and Young's inequality (using $\|z_x\|_{L^\infty}\|z_x\|_{L^2}^2\le C\mathcal{E}_1^{5/4}\mathcal{D}_1^{1/4}$ at $k=1$), the first term on the right-hand side of \eqref{af:primitive-derivative-energy} is bounded by
		\[
		{\epsilon\mathcal D_1+C_\epsilon\mathcal{E}_1^{5/3}\quad(k=1)},\qquad
		C\mathcal{E}_1^{1/4}\mathcal{E}_2^{5/4}\quad(k=2),\qquad
		C\mathcal{E}_1^{1/4}\mathcal{E}_2^{1/4}\mathcal{E}_k\quad(k\ge3),
		\]
		while interpolation between adjacent derivative orders gives
		\[
		\mathcal{D}_1\ge c\min\left\{\mathcal{E}_1,\frac{\mathcal{E}_1^2}{\|g\|_{L^2}^2}\right\},\qquad
		\mathcal{D}_k\ge c\min\left\{\mathcal{E}_k,\frac{\mathcal{E}_k^2}{\mathcal{E}_{k-1}}\right\}\quad(k\ge 2).
		\]
		Comparing $\mathcal{E}_k(t)$ inductively with $H_kT^{-(2k+1)a}$ for successively large $H_1,H_2,\dots,H_{N-1}$ gives
		\[
		\mathcal{E}_k(t)\le C_kT^{-(2k+1)a}\qquad(1\le k\le N-1,\ 0\le t\le t_*).
		\]
		Since $\|g_k\|_{L^2}^2\le C\|z_k\|_{L^2}^2+CT^{-k-1}\|z\|_{L^2}^2$ and $T^{-k-1-a}\le T^{-(2k+1)a}$, interpolation between adjacent orders gives \eqref{af:bootstrap}.
	\end{proof}
	
	\begin{remark}
		Taking $N=22$ provides derivative bounds up to order $21$ in $L^2(\mathbb{R})$ and order $20$ in $L^\infty(\mathbb{R})$, which suffices for the order-16 interpolation in Section~\ref{sec:exterior-L1}, the order-8 interpolation in Section~\ref{sec:central-energy-replacement}, and all third- and fourth-derivative estimates below. If $x=x_b(t,y)$ denotes the Lagrangian mass map with $(x_b)_y=v=\rho_b^{-1}$ and $\partial_y=v\partial_x$, expanding
		\[
		(v\partial_x)^jg=v^j\partial_x^jg+\sum_{i=1}^{j-1}c_{j,i}(v,\dots,\partial_x^{j-i}v)\partial_x^ig
		\]
		shows that the derivatives of $g$ in $y$ satisfy the same bounds \eqref{af:bootstrap}. For $U_b$ itself, $\|U_{b,y}\|_{L^2}\le CT^{-1/2}$ because of the profile slope $V_x$, while all higher derivatives of $U_b$ and all $L^\infty$ derivative norms obey \eqref{af:bootstrap}.
	\end{remark}
	
	\subsection{Exterior localization}
	\label{sec:exterior-L1}
	
	In this section, we combine the spatial decay of the initial datum \eqref{main:burgers-data} with Proposition~\ref{prop:planar-energy} to prove time-integrable $L^1(\mathbb{R})$ bounds for $\partial_x^3 g$ and $\partial_x^4 g$.
	
	\subsubsection{Exterior relative entropy}
	
	Fix either far-field state $U_\pm=(\rho_\pm,u_\pm,\theta_\pm)$, and let $\eta_\pm(U)$ be the relative entropy density with respect to $U_\pm$:
	\[
	\eta_\pm(U)=
	\frac{1}{2}\rho|u-u_\pm|^2
	+c_v\rho\theta_\pm\Phi_0\left(\frac{\theta}{\theta_\pm}\right)
	+R_g\rho\theta_\pm\Phi_0\left(\frac{\rho_\pm}{\rho}\right),
	\qquad \Phi_0(s)=s-1-\log s.
	\]
	On the compact set of positive states, $\eta_\pm(U)$ is comparable to $|U-U_\pm|^2$. Differentiating along \eqref{planar:conserved-equation} gives
	\[
	\partial_t\eta_\pm+\partial_x q_\pm+d_\pm=\partial_x j_\pm,
	\]
	where $d_\pm\ge 0$ is comparable to $|u_x|^2+|\theta_x|^2$, and the fluxes $q_\pm$ and $j_\pm$ satisfy
	\[
	|q_\pm|\le C_*\eta_\pm,\qquad
	|j_\pm|\le C_J\sqrt{\eta_\pm d_\pm}.
	\]
	
	Let $\omega:\mathbb{R}\to(0,\infty)$ be a smooth nondecreasing weight with $\omega'(z)/\omega(z)\le K$ for all $z\in\mathbb{R}$ and
	\[
	\omega(z)\approx e^{kz}\quad(z<0),\qquad
	\omega(z)\approx(1+z)^{2m}\quad(z>0).
	\]
	Multiplying the relative entropy identity for $U_+$ by $\omega(x-Ct)$, integrating by parts in $x$, and using the Cauchy--Schwarz inequality gives, for $C>C_*+C_J^2K/2$,
	\begin{align}
		&\int_{\mathbb{R}}\omega(x-Ct)\eta_+\bigl(U_b(t,x)\bigr)\dd x
		+\frac{1}{2}\int_0^t\!\int_{\mathbb{R}}\omega(x-Cs)d_+\dd x\dd s\notag\\
		&\quad+
		\left(C-C_*-\frac{1}{2}C_J^2K\right)
		\int_0^t\!\int_{\mathbb{R}}\omega'(x-Cs)\eta_+\dd x\dd s
		\le\int_{\mathbb{R}}\omega(x)\eta_+\bigl(U_b(0,x)\bigr)\dd x.
		\label{ex:weighted}
	\end{align}
	By \eqref{main:burgers-data}, $|U_b(0,x)-U_\pm|\le C\delta(1+|x|)^{-3}$ as $\pm x\to\infty$, so the right-hand side of \eqref{ex:weighted} is bounded by $C_m\delta^2$ for every fixed $m<5/2$. Together with the symmetric estimate as $x\to-\infty$, this yields
	\begin{equation}
		\|U_b(t)-U_\pm\|_{L^2(\{\pm x\ge Ct+R\})}
		\le C_m\delta(1+R)^{-m},\qquad R\ge 0,\quad m<\frac{5}{2}.
		\label{ex:state}
	\end{equation}
	Since $W(t,x)$ satisfies the same exterior estimate, $g=\mathbf{U}_b-\mathbf{V}(W)$ also satisfies \eqref{ex:state}, and the same bound holds in mass coordinates $y$ after enlarging $C$.
	
	\subsubsection{Annular Gagliardo--Nirenberg interpolation}
	
	\begin{lemma}
		\label{lem:exterior-L1}
		Fix $m=12/5$. Suppose that on $[0,t_*]$ the planar error $g$ satisfies \eqref{ex:state} together with
		\[
		\sup_{0\le t\le t_*}\|g(t)\|_{H^{16}}\le A_\delta\quad\left(\lim_{\delta\rightarrow 0}A_\delta=0\right),
		\qquad
		\|\partial_x^j g(t)\|_{L^2}\le C_jT^{-(j+1/2)a}\quad(1\le j\le 16),
		\]
		with $C_j$ independent of $t_*$ and $\delta$. Then for $j=3,4$ and every $0\le\epsilon<7a/2-3/2=9/50$,
		\begin{equation}
			\int_0^{t_*} T^\epsilon\|\partial_x^j g(t)\|_{L^1}\dd t=o_\delta(1),
			\label{ex:L1-high}
		\end{equation}
		where $\epsilon=1/20$ is admissible for both $j=3$ and $j=4$, and any fixed power of $\log(2+t)$ may be included in the integrand. Moreover,
		\begin{equation}
			\sup_{0\le t\le t_*}\bigl(\|g_x(t)\|_{L^1}+\|g_{xx}(t)\|_{L^1}\bigr)=o_\delta(1),
			\qquad
			\|g_x(t)\|_{L^1}+\|g_{xx}(t)\|_{L^1}\le CT^{-11/50}.
			\label{ex:variation}
		\end{equation}
	\end{lemma}
	
	\begin{proof}
		We split $\mathbb{R}$ into the interior region $\{|x|\le C_0T\}$ and dyadic exterior annuli. On $\{|x|\le C_0T\}$, the Cauchy--Schwarz inequality gives
		\[
		\|\partial_x^j g(t)\|_{L^1(\{|x|\le C_0T\})}
		\le CT^{1/2}\min\bigl\{A_\delta,C_jT^{-s_j}\bigr\},
		\qquad s_j=\left(j+\frac{1}{2}\right)a.
		\]
		Whenever $s_j>3/2+\epsilon$, integrating this minimum over $t\in[0,t_*]$ gives
		\[
		\int_0^{t_*}T^{\epsilon+1/2}\min\bigl\{A_\delta,C_jT^{-s_j}\bigr\}\dd t
		\le C_{\epsilon,j}A_\delta^{1-(3/2+\epsilon)/s_j}=o_\delta(1).
		\]
		Since $a=12/25$, we have {$s_3=42/25=3/2+9/50$} and $s_4=54/25$, which bounds the interior part of \eqref{ex:L1-high}.
		
		On $\{|x|\ge C_0T\}$, we cover the region by dyadic annuli $I_R=\{R\le |x|\le 2R\}$ with $R=2^k C_0T$ ($k\ge 0$), and let $I_R^*$ be the slightly enlarged annulus. Applying the Gagliardo--Nirenberg inequality on $I_R^*$ followed by the Cauchy--Schwarz inequality and \eqref{ex:state} with $m=12/5$ gives
		\begin{align}
			\|\partial_x^j g\|_{L^1(I_R)}
			&\le C R^{1/2}\|g\|_{L^2(I_R^*)}^{1-j/16}\|\partial_x^{16} g\|_{L^2(I_R^*)}^{j/16}
			+CR^{1/2-j}\|g\|_{L^2(I_R^*)}\notag\\
			&\le C\delta^{1-j/16}A_\delta^{j/16}R^{1/2-m(1-j/16)}
			+C\delta R^{1/2-j-m}.
			\label{ex:annulus}
		\end{align}
		For $j=3$ and $j=4$, the exponents of $R$ in the first term are $-29/20<-1$ and $-13/10<-1$. Summing \eqref{ex:annulus} over $R=2^kC_0T$ and integrating in time against $T^\epsilon$ with $\epsilon<9/50$ gives a convergent integral of size $O(\delta^{1-j/16}A_\delta^{j/16}+\delta)=o_\delta(1)$.
		
		For $j=1,2$, the exterior sum in \eqref{ex:annulus} is bounded by $o_\delta(1)T^{1/2-m(1-j/16)}$, while on $\{|x|\le C_0T\}$ the Cauchy--Schwarz inequality gives
		\[
		T^{1/2-3a/2}=T^{-11/50}\quad(j=1),\qquad
		T^{1/2-5a/2}=T^{-7/10}\quad(j=2),
		\]
		and $\sup_{t\ge 0}T^{1/2}\min\{A_\delta,CT^{-s_j}\}=o_\delta(1)$ since $s_1=18/25>1/2$. This proves \eqref{ex:variation}.
	\end{proof}
	
	\begin{remark}
		Using $\|\partial_x^{16}g(t)\|_{L^2}\le C_{16}T^{-(33/2)a}$ in \eqref{ex:annulus} in place of $A_\delta$ gains an extra factor $T^{-ja(1+1/32)}$ in the exterior sum, which yields
		\[
		\|\partial_x^3 U_b(t)\|_{L^1}\le CT^{1/2-7a/2}.
		\]
		Keeping $R$ and $T$ separate in \eqref{ex:annulus} also gives the exterior bound \eqref{new:source-tail} used below.
	\end{remark}
	
	\subsection{Characteristic interactions and decay}
	\label{sec:allfamilyLone}
	
	Using Proposition~\ref{prop:planar-energy} and Lemma~\ref{lem:exterior-L1}, we now prove the unconditional $L^1$ and $L^2$ decay estimates for the planar wave.
	
	\begin{proposition}
		\label{prop:planar-decay}
		The global planar solution $U_b$ of Proposition~\ref{prop:planar-energy} satisfies
		\[
		\|g(t)\|_{L^1}\le C\bigl(\varepsilon_\delta+\delta\log(2+t)\bigr),\qquad t\ge 0,
		\qquad \lim_{\delta\rightarrow 0}\varepsilon_\delta=0,
		\]
		and the bounds $\|g(t)\|_{L^2}^2\le CT^{-a}$ and \eqref{af:bootstrap} hold for all $t\ge0$ with $a=12/25$. More generally, for every fixed $a_*<1/2$, there exists $C_{a_*}>0$ such that
		\[
		\sup_{t\ge0}T^{a_*}\|g(t)\|_{H^2}^2+
		\int_0^\infty T^{a_*}\mathcal{D}(t)\dd t\le C_{a_*},
		\]
		where $\mathcal{D}$ is the dissipation functional in \eqref{af:physical-energy}. Consequently, $U_b$ converges in $L^\infty(\mathbb{R})$ toward both $\overline{U}_E$ and $U^R$ at rate $T^{-12/25}$:
		\[
		\|U_b(t)-\overline{U}_E(t)\|_{L^\infty}\le C_\delta T^{-12/25}\quad(t\ge0),\qquad
		\|U_b(t)-U^R(t)\|_{L^\infty}\le C_\delta T^{-12/25}\quad(t\ge1).
		\]
	\end{proposition}
	
	\subsubsection{Diffusive and quadratic corrections}
	
	We work on a bootstrap interval $[0,t_*]$ on which $\|g(t)\|_{L^2}^2\le T^{-a}$. Let $A=F'(\mathbf{V}(W))$ be the Euler flux Jacobian at $\mathbf{V}(W)$, and let $R(W)$ and $L(W)=R(W)^{-1}$ be the right and left eigenvector matrices of $A$, ordered as $i\in\{-,0,+\}$ so that $A\mathbf{V}'(W)=W\mathbf{V}'(W)$. We set
	\[
	\left\{
	\begin{aligned}
		&\Lambda=LAR=\operatorname{diag}(\lambda_-,\lambda_0,\lambda_+),\qquad \lambda_+=W,\\
		&D=LBR,\qquad d_i=D_{ii}>0,\qquad D_d=\operatorname{diag}(d_-,d_0,d_+),
	\end{aligned}
	\right.
	\]
	and, following the first-order cross-family derivative correction method of \cite{HLX25,LWXZ26,LSWX26} (cf.\ (1.8) in \cite{LSWX26}), define
	\[
	q=Lg,\qquad G=q+Kq_x,\qquad K_{ii}=0,\qquad K_{ij}=-\frac{D_{ij}}{\lambda_i-\lambda_j}\quad(i\ne j).
	\]
	In primitive variables $(\rho,u,\theta)$, the right eigenvectors for the two acoustic modes and the entropy mode are $(1,\pm c_E/\rho_v,(\gamma-1)\theta_v/\rho_v)^\top$ and $(1,0,-\theta_v/\rho_v)^\top$, and the left eigenvectors are
	\[
	\left(\frac{1}{2\gamma},\ \pm\frac{\rho_v}{2c_E},\ \frac{\rho_v}{2\gamma\theta_v}\right),\qquad
	\left(\frac{\gamma-1}{\gamma},\ 0,\ -\frac{\rho_v}{\gamma\theta_v}\right).
	\]
	Pairing these with $\operatorname{diag}(0,\nu/\rho_v,\kappa/(\rho_v c_v))$ gives the diagonal diffusion coefficients
	\[
	d_-=d_+=\frac{\nu+(\gamma-1)\kappa/c_p}{2\rho_v},\qquad
	d_0=\frac{\kappa}{c_p\rho_v},\qquad
	d_*=d_+(W_-)=\frac{\nu+(\gamma-1)\kappa/c_p}{2\rho_-}>0.
	\]
	
	By Proposition~\ref{prop:planar-energy} and Lemma~\ref{lem:exterior-L1}, we have $\sup_{0\le t\le t_*}\|g(t)\|_{H^N}=o_\delta(1)$ and
	\begin{equation}\label{af:Lone-high}
		\int_0^{t_*}(1+\log T)^m
		\bigl(\|\partial_x^3 g(t)\|_{L^1}+\|\partial_x^4 g(t)\|_{L^1}\bigr)\dd t
		=o_\delta(1)\qquad(m\ge 0),
	\end{equation}
	together with
	\begin{equation}\label{af:TV-input}
		\|g_x(t)\|_{L^1}+\|g_{xx}(t)\|_{L^1}
		\le \tau_\delta(t),\qquad
		\sup_{0\le t\le t_*}\tau_\delta(t)=o_\delta(1),\qquad
		\tau_\delta(t)\le C T^{-\theta_0},\quad \theta_0=\frac{3a}{2}-\frac{1}{2}=\frac{11}{50}.
	\end{equation}
	By Lemma~\ref{lem:planar-burgers}, the Burgers profile satisfies
	\[
	\left\{
	\begin{aligned}
		&\|W_x(t)\|_{L^2}\le C\delta^{1/2}T^{-1/2},\qquad
		\|W_{xx}(t)\|_{L^1}+\|W_x(t)^2\|_{L^1}\le CT^{-1},\\
		&\int_0^\infty(1+\log T)^m\bigl(\|\partial_x^3W(t)\|_{L^1}+\|\partial_x^4W(t)\|_{L^1}\bigr)\dd t=o_\delta(1).
	\end{aligned}
	\right.
	\]
	
	The profile $\mathbf{V}(W)$ leaves the residual
	\[
	\mathcal{R}_v=\bigl(B(\mathbf{V})\mathbf{V}_x\bigr)_x-\mathbf{V}_t-F(\mathbf{V})_x
	=\bigl(B(\mathbf{V})-d_*I\bigr)\mathbf{V}'W_{xx}+\bigl(B(\mathbf{V})\mathbf{V}'\bigr)'W_x^2.
	\]
	In the coordinates $q=Lg$, differentiating $R(W)$ produces the coupling matrix
	\[
	H=L\bigl(R_t+(AR)_x\bigr)
	=LR'(\Lambda-WI)W_x+\Lambda'W_x+d_*LR'W_{xx}.
	\]
	Because $\lambda_+-W=0$, the $(i,+)$-entries of $C_0(W)=-LR'(\Lambda-WI)$ vanish:
	\begin{equation}\label{af:profilematrix}
		(\Lambda_x-H)G=C_0(W)W_xG+O(W_{xx})G,
		\qquad C_0(W)=-LR'(\Lambda-WI),\qquad (C_0)_{i,+}=0.
	\end{equation}
	Applying $I+K\partial_x$ to the equation for $q$ and using $D+[\Lambda,K]=D_d$ to cancel the off-diagonal second derivatives gives
	\begin{equation}\label{af:G}
		G_{i,t}+(\lambda_iG_i)_x
		=(d_iG_{i,x})_x+\sum_j (C_0)_{ij}W_xG_j
		+\partial_x Q_{i,2}(W,G)+b_i^{(0)}+\mathcal{R}_i^{(0)},
	\end{equation}
	where
	\[
	\left\{
	\begin{aligned}
		&Q_{i,2}(W,G)=-\frac{1}{2}\bigl[L(W)F''(\mathbf{V}(W))[R(W)G,R(W)G]\bigr]_i,\\
		&b_i^{(0)}=(I+K\partial_x)(L\mathcal{R}_v)_i=a_i(W)W_{xx}+O(W_x^2)+\text{higher-order profile derivatives},\\
		&a_i(W)=L_i(W)\bigl(B(\mathbf{V}(W))-d_*I\bigr)\mathbf{V}'(W).
	\end{aligned}
	\right.
	\]
	
	To remove the terms $a_i(W)W_{xx}$ for $i\in\{-,0\}$ and the squares $q_{ij}(W)G_j^2$ for $j\ne i$, we perform two further transformations:
	\begin{enumerate}
		\item[(1)] {Since $\mathbf{V}'(W)=\rho_v'(W)R_+(W)$, where $\rho_v'(W)=2\rho_v/((\gamma+1)c_E)$, we have $a_+(W)=\rho_v'(W)(d_+(W)-d_*)=O(\delta)$}. For $i\ne +$, $\lambda_i-W<0$ is bounded away from zero, and we set
		\[
		\chi_i=-\frac{a_i(W)}{\lambda_i-W}\quad(i\ne +),\qquad \chi_+=0,\qquad
		\widetilde{G}_i=G_i+\chi_iW_x.
		\]
		{Since
			\[
			\begin{aligned}
				(\partial_t+\lambda_i\partial_x)(\chi_i W_x)
				&=\chi_i(\lambda_i-W)W_{xx}+d_*\chi_i W_{xxx}\\
				&\quad+[\chi_i'(\lambda_i-W)-\chi_i]W_x^2+d_*\chi_i'W_xW_{xx},
			\end{aligned}
			\]
			this cancels} $a_i(W)W_{xx}$ for $i\ne +$, leaving a profile source $b=(b_-,b_0,b_+)$ that satisfies
		\begin{equation}\label{af:puresource}
			\|b(t)\|_{L^1}\le \frac{C\delta}{T}+b_1(t),\qquad
			\int_0^{t_*}(1+\log T)^m\bigl(b_1(t)+\|b_x(t)\|_{L^1}\bigr)\dd t=o_\delta(1).
		\end{equation}
		\item[(2)] Writing the quadratic flux as
		\[
		Q_{i,2}(W,h)=q_i(W)h_i^2+\sum_{j\ne i}q_{ij}(W)h_j^2+\sum_{j<k}q_{ijk}(W)h_jh_k,
		\]
		we remove the squares $q_{ij}(W)\widetilde{G}_j^2$ ($j\ne i$) by setting
		\begin{equation}\label{af:normalform}
			h_i=\widetilde{G}_i-\sum_{j\ne i}\frac{q_{ij}(W)}{\lambda_i-\lambda_j}\widetilde{G}_j^2.
		\end{equation}
		Because $(\partial_t+\lambda_i\partial_x)\widetilde{G}_j^2=(\lambda_i-\lambda_j)\partial_x(\widetilde{G}_j^2)$ up to viscous, profile, and cubic terms, $h=(h_-,h_0,h_+)$ satisfies
		\begin{equation}\label{af:h-system}
			h_{i,t}+\partial_x\bigl(\lambda_i h_i-q_i h_i^2\bigr)
			=\partial_x(d_i h_{i,x})+\sum_j(C_0)_{ij}W_xh_j
			+\partial_x\widetilde{Q}_i(W,h)+b_i+\mathcal{R}_i,
		\end{equation}
		where $\widetilde{Q}_i(W,h)=\sum_{j<k}q_{ijk}(W)h_jh_k$ contains only products $h_jh_k$ with $j\ne k$.
	\end{enumerate}
	
	Explicitly, under \eqref{planar:conserved-equation} we have
	\[
	g_t=-\partial_x\bigl(F(\mathbf{V}+g)-F(\mathbf{V})\bigr)
	+\partial_x\bigl(B(\mathbf{V}+g)(\mathbf{V}_x+g_x)-B(\mathbf{V})\mathbf{V}_x\bigr)+\mathcal{R}_v,
	\]
	so that $q_t=L_tg+Lg_t$ and $G_t=q_t+K_tq_x+K\partial_xq_t$. In \eqref{af:G}, $b^{(0)}=(I+K\partial_x)(L\mathcal{R}_v)$ and
	\[
	\mathcal{R}^{(0)}=G_t+\partial_x(\Lambda G)-\partial_x(D_dG_x)-C_0(W)W_xG-\partial_xQ_2(W,G)-b^{(0)}.
	\]
	Computing $h_t$ by the chain rule with $W_t=d_*W_{xx}-WW_x$ and setting
	\[
	N_i(W,g)=h_{i,t}+\partial_x(\lambda_i h_i-q_i h_i^2)
	-\partial_x(d_i h_{i,x})-\sum_j(C_0)_{ij}W_xh_j-\partial_x\widetilde{Q}_i(W,h),
	\]
	we have $b_i=N_i(W,0)$ and $\mathcal{R}_i=N_i(W,g)-N_i(W,0)$. By Taylor's formula and $D+[\Lambda,K]=D_d$, the terms in $\mathcal{R}_i$ and $\partial_x\mathcal{R}_i$ belong to the following eight groups:
	\begin{equation}\label{tab:planar-remainders}
		\left\{
		\begin{aligned}
			&\text{(i) }\partial_x^3 q\text{ and }\partial_x^4 q, {\text{ use $q=L(W)g$, the product rule, and \eqref{ex:L1-high}}};\\
			&\text{(ii) }W_x q_x=O_{L^1}(T^{-61/50})\text{ and }\partial_x(W_x q_x)=O_{L^1}(T^{-73/50});\\
			&\text{(iii) }{(W_{xx}+W_x^2)q=O_{L^1}(T^{-149/100})\text{ and its derivative }O_{L^1}(T^{-197/100});}\\
			&\text{(iv) }q q_{xx},\,q_x^2=O_{L^1}(T^{-36/25})\text{ and }\partial_x(q q_{xx},q_x^2)=O_{L^1}(T^{-48/25});\\
			&\text{(v) }q q_{xxx},\,q_x q_{xx}=O_{L^1}(T^{-48/25})\text{ and their derivatives }O_{L^1}(T^{-5a});\\
			&\text{(vi) }\partial_x O(q^3)=O_{L^1}(T^{-36/25})\text{ and }\partial_x^2 O(q^3)=O_{L^1}(T^{-48/25});\\
			&\text{(vii) }{W_x q^2=O_{L^1}(T^{-37/25})\text{ and its derivative }O_{L^1}(T^{-49/25});}\\
			&\text{(viii) }q b=O_{L^1}(T^{-37/25}+T^{-a}b_1)\text{ and }\partial_x(qb).
		\end{aligned}
		\right.
	\end{equation}
	Since $a=12/25$, every decay exponent in \eqref{tab:planar-remainders} is at least $59/50>1$. Combining these integrable decay bounds with $\sup_t\|g(t)\|_{H^N}=o_\delta(1)$ on finite time intervals gives
	\begin{equation}\label{af:Rbudget}
		\int_0^{t_*}(1+\log T)^m
		\bigl(\|\mathcal{R}(t)\|_{L^1}+\|\mathcal{R}_x(t)\|_{L^1}\bigr)\dd t=o_\delta(1).
	\end{equation}
	
	\subsubsection{Weighted \texorpdfstring{$L^1$}{L1} estimates and characteristic interactions}
	
	Setting $w_i=h_{i,x}$, we have from \eqref{af:TV-input} that
	\begin{equation}\label{af:tau}
		\tau(t)=\sum_i\|w_i(t)\|_{L^1}\le\varepsilon_\delta\quad\left(\lim_{\delta\rightarrow 0}\varepsilon_\delta=0\right),
		\qquad \tau(t)\le CT^{-\theta_0},\quad \theta_0=\frac{11}{50}.
	\end{equation}
	Differentiating \eqref{af:h-system} with respect to $x$ shows that $w_i$ satisfies
	\begin{equation}\label{af:w}
		w_{i,t}+\partial_x\bigl[(\lambda_i-2q_i h_i)w_i\bigr]
		=\partial_x^2(d_iw_i)+F_i,
	\end{equation}
	where
	\[
	F_i=\partial_x\left(
	\sum_j(C_0)_{ij}W_xh_j+\partial_x\widetilde{Q}_i+b_i+\mathcal{R}_i
	-\lambda_i'W_xh_i+q_i'W_xh_i^2\right)
	\]
	satisfies, by \eqref{af:puresource}, \eqref{tab:planar-remainders}, and \eqref{af:Rbudget},
	\begin{equation}\label{af:Fbudget}
		\int_0^{t_*}(1+\log T)^m\sum_i\|F_i(t)\|_{L^1}\dd t=o_\delta(1).
	\end{equation}
	Applying Kato's inequality to \eqref{af:w} gives
	\[
	(|w_i|)_t+\partial_x\bigl[(\lambda_i-2q_i h_i)|w_i|\bigr]
	\le\partial_x^2(d_i|w_i|)+|F_i|.
	\]
	
	To bound $\|h_i\|_{L^1}$, we define the weights $\beta_+=1$ and $\beta_i(W)=\exp(K_i(W-W_-))$ for $i\in\{-,0\}$. Since $\lambda_i-W<0$ for $i\in\{-,0\}$, taking $K_-,K_0>0$ large enough ensures that
	\[
	-(\partial_t+\lambda_i\partial_x)\beta_i=-K_i(\lambda_i-W)\beta_iW_x+O(W_{xx})
	\]
	produces a large positive multiple of $W_x|h_i|$ for $i\in\{-,0\}$. Because $(C_0)_{i,+}=0$ by \eqref{af:profilematrix}, the coupling $\sum_j(C_0)_{ij}W_xh_j$ involves only $h_-$ and $h_0$, and is absorbed by this $W_x$ term. Setting
	\[
	M(t)=\sum_i\int_{\mathbb{R}}\beta_i|h_i|\dd x,\qquad
	P(t)=\sum_{i\ne+}\int_{\mathbb{R}} W_x|h_i|\dd x,\qquad
	K_1(t)=\sum_{i\ne j}\int_{\mathbb{R}}|h_i||w_j|\dd x,
	\]
	and using $|\partial_x\widetilde{Q}_i(W,h)|\le C\sum_{j\ne k}(|h_j||w_k|+W_x|h_j||h_k|)$, we obtain
	\begin{equation}\label{af:Mineq}
		M'(t)+cP(t)\le CK_1(t)+k(t)M(t)+\frac{C\delta}{T}+g_0(t),
		\qquad \int_0^{t_*}(1+\log T)^m\bigl(k(t)+g_0(t)\bigr)\dd t<\infty,
	\end{equation}
	where $k(t)\le C(\|W_{xx}\|_{L^\infty}+\|W_x\|_{L^\infty}^2)$.
	
	To control $CK_1(t)$ in \eqref{af:Mineq}, we define for $i\ne j$ the Glimm--Bony interaction functional
	\[
	I_{ij}(t)=\begin{cases}
		\displaystyle\int_{\mathbb{R}}|h_i(t,x)|\left(\int_{-\infty}^x|w_j(t,y)|\dd y\right)\dd x,
		&\text{if }\lambda_i<\lambda_j,\\[10pt]
		\displaystyle\int_{\mathbb{R}}|h_i(t,x)|\left(\int_x^\infty|w_j(t,y)|\dd y\right)\dd x,
		&\text{if }\lambda_i>\lambda_j.
	\end{cases}
	\]
	Differentiating $I_{ij}(t)$ along the Kato inequalities for $|h_i|$ and $|w_j|$, the speed gap $|\lambda_i-\lambda_j|\ge c_0>0$ gives the negative term $-c\int_{\mathbb{R}}|h_i||w_j|\dd x$. Integrating the diffusion terms by parts once gives
	\[
	\mp\int_{\mathbb{R}}(d_i+d_j)\partial_x(|h_i|)|w_j|\dd x,
	\]
	which is bounded by $C\int_{\mathbb{R}}|w_i||w_j|\dd x$. {By the chain rule in \eqref{af:normalform},
		\[
		\sum_i\|w_i\|_{L^2}^2\le C\left(\mathcal D+\|W_x\|_{L^\infty}^2\|z\|_{L^2}^2
		+\|W_{xx}\|_{L^2}^2+\|W_x^2\|_{L^2}^2\right).
		\]
		Since} $\int_{\mathbb{R}}|w_j|\dd y\le\tau(t)\le\varepsilon_\delta$, summing $I(t)=\sum_{i\ne j}I_{ij}(t)$ gives
	\begin{equation}\label{af:Iineq}
		I'(t)+cK_1(t)\le C\varepsilon_\delta K_1(t)+k_1(t)M(t)+g_1(t),
		\qquad I(t)\le C\varepsilon_\delta M(t),
	\end{equation}
	where $k_1(t)=C\sum_j\|F_j(t)\|_{L^1}+C\tau(t)\|W_x(t)\|_{L^\infty}$ and $g_1(t)$ satisfy
	\[
	\int_0^{t_*}(1+\log T)k_1(t)\dd t<\infty,\qquad
	\int_0^{t_*}g_1(t)\dd t<\infty.
	\]
	Choosing $C_I>0$ large and $\delta>0$ small so that $C_I(c-C\varepsilon_\delta)\ge 2C$, we add $C_I\times\eqref{af:Iineq}$ to \eqref{af:Mineq}:
	\[
	(M+C_I I)'\le C(k+k_1)(M+C_I I)+\frac{C\delta}{T}+g_0+C_I g_1.
	\]
	{Here $M(0)+C_I I(0)\le C\delta$ by \eqref{af:normalform}. Moreover,
		\[
		\begin{aligned}
			g_0\le C(b_1+\|\mathcal R\|_{L^1}),\qquad
			g_1\le C\sum_i\|w_i\|_{L^2}^2+C\tau(\delta/T+g_0).
		\end{aligned}
		\]
		By \eqref{af:direct-classical}, \eqref{af:puresource}, \eqref{af:Rbudget}, \eqref{af:tau}, and \eqref{af:Fbudget},
		$\int_0^{t_*}(k+k_1)\,\dd t\le C$ and $\int_0^{t_*}(g_0+g_1)\,\dd t=o_\delta(1)$.
		Gronwall's inequality then gives}
	\begin{equation}\label{af:massresult}
		M(t)\le C\bigl(\varepsilon_\delta+\delta\log T\bigr),
	\end{equation}
	and integrating \eqref{af:Iineq} yields
	\begin{equation}\label{af:crossresult}
		\int_0^{t_*} K_1(t)\dd t<\infty.
	\end{equation}
	{By \eqref{af:normalform}, $\|G\|_{L^1}\le C(M+\delta)$, while
		\[
		\|q\|_{L^1}\le\|G\|_{L^1}+C\|q_x\|_{L^1},\qquad
		\|q_x\|_{L^1}\le C\|g_x\|_{L^1}+C\delta\|g\|_{L^\infty}.
		\]
		Thus \eqref{af:TV-input} gives}
	\begin{equation}\label{af:physicalmass}
		\|\mathbf{U}_b(t)-\mathbf{V}(W(t))\|_{L^1}
		\le C\bigl(\varepsilon_\delta+\delta\log(2+t)\bigr).
	\end{equation}
	
	\subsubsection{Closure of the planar decay bootstrap}
	
	\begin{proof}[Proof of Proposition~\ref{prop:planar-decay}]
		Let $\Gamma(t)=\int_0^t k_0(s)\dd s\le C$, and multiply \eqref{af:physical-energy} by $T^{a_*}e^{-\Gamma(t)}$ for a fixed $a_*<1/2$. {By the one-dimensional Nash inequality applied to $z$, the pointwise equivalence $|z|\approx|g|$, and \eqref{af:physicalmass}},
		\[
		\|g(t)\|_{L^2}^2\le C M_g(t)^{4/3}\mathcal{D}(t)^{1/3},\qquad
		M_g(t)=\|g(t)\|_{L^1}\le C(1+\log T).
		\]
		Young's inequality gives
		\[
		T^{a_*-1}\|g\|_{L^2}^2
		\le\epsilon T^{a_*}\mathcal{D}+C_\epsilon T^{a_*-3/2}M_g^2,
		\qquad
		T^{a_*}\|\mathcal{R}_0\|_{L^2}\|g\|_{L^2}
		\le\epsilon T^{a_*}\mathcal{D}+C_\epsilon T^{a_*} M_g^{4/5}\|\mathcal{R}_0\|_{L^2}^{6/5}.
		\]
		Since $\|\mathcal{R}_0(t)\|_{L^2}^{6/5}\le CT^{-3/2}$ and {$Q_0(t)\le CT^{-3/2}$}, all right-hand terms are integrable on $[0,\infty)$ for $a_*<1/2$. Integrating over $[0,t_*]$ gives
		\[
		\sup_{0\le t\le t_*}T^{a_*}\|g(t)\|_{H^2}^2+
		\int_0^{t_*} T^{a_*}\mathcal{D}(t)\dd t\le C_{a_*},
		\]
		with $C_{a_*}$ independent of $t_*$ and $\delta$. Taking $a_*=49/100>a=12/25$ and interpolating with $\sup_t\|g(t)\|_{L^2}^2\le A_\delta\to 0$ from \eqref{af:direct-classical} yields
		\[
		\sup_{0\le t\le t_*}T^a\|g(t)\|_{L^2}^2
		\le C_{a_*}^{a/a_*}A_\delta^{1-a/a_*}\longrightarrow 0\quad(\delta\rightarrow 0).
		\]
		For small $\delta>0$, this closes the bootstrap $\|g(t)\|_{L^2}^2\le T^{-a}$ for all $t\ge 0$.
		
		To compare $U_b$ with $\overline{U}_E$ and $U^R$, let $W_T^R(t,x)=W^R(x/T)$ be the centered rarefaction fan shifted to $T=1+t$. It satisfies
		\[
		(W_T^R)_{xx}=T^{-1}(\delta_{x=W_-T}-\delta_{x=W_+T}),
		\]
		and $W(0)-W_T^R(0)\in L^1(\mathbb{R})$. The difference $W-W_T^R$ solves
		\[
		(W-W_T^R)_t+\partial_x\left(\frac{W+W_T^R}{2}(W-W_T^R)\right)
		=d_*(W-W_T^R)_{xx}+d_*(W_T^R)_{xx}.
		\]
		Since $W_x+(W_T^R)_x\ge 0$, the $L^1\to L^\infty$ smoothing bound $C(t-s)^{-1/2}$ from Lemma~\ref{lem:planar-burgers} and the total variation bound $\|d_*(W_T^R)_{xx}\|_{M^1}=2d_*/T$ give, upon splitting the Duhamel integral at $s=t/2$,
		\[
		\|W(t)-W_T^R(t)\|_{L^\infty}\le C T^{-1/2}\log(2+t).
		\]
		For $W^E(t,x)$, its characteristic formula $x=\xi+Tw_0(\xi)$ implies $|W^E(t,x)-W_T^R(t,x)|\le R/T$ for $|\xi|\le R$ and $|W^E(t,x)-W_T^R(t,x)|\le CR^{-3}$ for $|\xi|>R$. Taking $R=T^{1/4}$ gives
		\[
		\|W^E(t)-W_T^R(t)\|_{L^\infty}\le CT^{-3/4},\qquad
		\|W_T^R(t)-W^R(\cdot/t)\|_{L^\infty}\le Ct^{-1}\quad(t\ge 1).
		\]
		Combining these bounds with $\|U_b(t)-V(W(t))\|_{L^\infty}\le C_\delta T^{-12/25}$ completes the proof.
	\end{proof}
	
	\subsubsection{Further interaction estimates}
	
	\begin{corollary}
		\label{cor:planar-mixed-interactions}
		The characteristic variables $h_i$ satisfy
		\begin{equation}\label{af:weightedcross}
			\int_0^\infty T^{1/20}
			\sum_{i\ne j}\int_{\mathbb{R}} |h_i||h_{j,x}|\dd x\dd t=o_\delta(1).
		\end{equation}
		Consequently, in Lagrangian mass coordinates $y$, the variables $r=u-I(p)$ and $\sigma=S-S_-$ of $U_b$ satisfy
		\[
		\int_0^\infty T^{1/20}
		\bigl(\|\sigma r_y\|_{L^1_y}+\|r\sigma_y\|_{L^1_y}\bigr)\dd t=o_\delta(1).
		\]
	\end{corollary}
	
	\begin{proof}
		Set $\epsilon=1/20$ and $a_*=49/100$. H\"older interpolation between \eqref{af:direct-classical} and Proposition~\ref{prop:planar-decay} gives
		\[
		\int_0^\infty T^{1/20}\mathcal{D}(t)\dd t
		\le\left(\int_0^\infty\mathcal{D}(t)\dd t\right)^{44/49}
		\left(\int_0^\infty T^{49/100}\mathcal{D}(t)\dd t\right)^{5/49}=o_\delta(1),
		\]
		{which bounds the contribution of $\mathcal D$. The other terms in the bound for $\sum_i\|w_i\|_{L^2}^2$ have $T^\epsilon$-weighted integral $o_\delta(1)$ by Lemma~\ref{lem:planar-burgers} and \eqref{af:direct-classical}.} By \eqref{af:tau} and \eqref{af:massresult},
		\[
		T^{\epsilon-1}I(t)\le CT^{1/20-1-11/50}\bigl(\varepsilon_\delta+\delta\log T\bigr)\in L^1(0,\infty),
		\]
		and the same holds for $T^\epsilon\tau(t)\delta/T$ and $T^\epsilon\tau(t)\|W_x\|_{L^\infty}M(t)$. {From the expression for $F_i$ following \eqref{af:w},
			\[
			\begin{aligned}
				\sum_i\|F_i\|_{L^1}\le C\bigl(&\|b_x\|_{L^1}+\|\mathcal R_x\|_{L^1}+T^{-1}\tau
				+\|W_{xx}\|_{L^2}\|h\|_{L^2}+T^{-2}M+\|h_x\|_{L^2}^2+\|h\|_{L^2}\|h_{xx}\|_{L^2}\bigr).
			\end{aligned}
			\]
			The first two terms admit the weight $T^{1/10}$ by \eqref{ex:L1-high}, the same exterior interpolation for $W$, and \eqref{tab:planar-remainders}. By \eqref{af:bootstrap} and the product rule, the last two terms are $O(T^{-3a})$. Since $1/20<1/10<9/50$, \eqref{af:massresult} and finite-time smallness give
			\[
			\int_0^{t_*}T^{1/20}\left(M\sum_i\|F_i\|_{L^1}+|g_1|\right)\dd t=o_\delta(1).
			\]
			Multiplying \eqref{af:Iineq} by $T^{1/20}$} and absorbing $C\varepsilon_\delta T^\epsilon K_1(t)$ gives \eqref{af:weightedcross}.
		
		Finally, since $r$ and $\sigma$ vanish on $V(W)$ by \eqref{v3:profile-identities} and their linearizations at $\mathbf{V}(W)$ are multiples of $q_-$ and $q_0$, {using \eqref{af:normalform} together with $q=G-Kq_x$ and $q_x=G_x-Kq_{xx}-K_xq_x$ expresses} $\sigma r_x$ and $r\sigma_x$ in terms of $K_1(t)$ up to terms controlled by $T^{1/20}\mathcal{D}(t)$ and terms decaying at least like $T^{-(1+3a)/2}+T^{-3a}$. Changing variables via $\partial_y=v\partial_x$ and $\dd y=v^{-1}\dd x$ gives the bound for $\|\sigma r_y\|_{L^1_y}+\|r\sigma_y\|_{L^1_y}$.
	\end{proof}
	
	\section{Entropy and central acoustic estimates}
	\label{sec:entropy-central-all}
	
	\subsection{Entropy decomposition}
	\label{sec:viscous-entropy}
	
	In Sections~\ref{sec:r5a-full-interface}--\ref{sec:common-energy}, the entropy derivative $\partial_1 S$ couples $\Omega$ to $\nabla h_E$. Because viscous heating along the rarefaction wave produces a nonnegative entropy component of size $O(\delta T^{-1/2}\log T)$ in $L^\infty$, its square is not integrable in time with a small logarithmic factor. We therefore decompose $\sigma=S-S_-$ in Lagrangian mass coordinates into a nonnegative solution $H_P\ge 0$ forced by the Burgers dissipation and a remainder $\sigma_0$. The derivative of $H_P$ is concentrated in the central region where the weight $w$ is of order $[\log(2+t)]^{-2}$, while $\sigma_0$ satisfies an $L^2_tL^\infty_y$ bound with a small $o_\delta(1)\log(2+t)$ growth.
	
	\begin{proposition}
		\label{prop:entropy-decomposition}
		Let $\sigma=S-S_-$ denote the relative entropy of the planar solution $U_b$. In Lagrangian mass coordinates $(t,y)$, there exist scalar functions $H_P(t,y)\ge0$ and $H_N(t,y)$ such that
		\[
		\sigma=H_P+\sigma_0,\qquad
		\sigma_0=H_N-\frac{\kappa_p u_y}{\theta},\qquad
		\kappa_p=\frac{\kappa}{c_p}.
		\]
		Here $H_N$ and $H_P$ satisfy the pointwise and derivative estimates \eqref{ev:HNpoint}, \eqref{ev:HPpoint}, \eqref{ev:HNy}, and \eqref{ev:HPy}, while $\sigma_0$ satisfies
		\[
		\int_0^t\|\sigma_0(s)\|_{L^\infty}^2\dd s
		\le A_\delta+\varepsilon_\delta\log(2+t),
		\qquad \lim_{\delta\rightarrow0}(A_\delta+\varepsilon_\delta)=0.
		\]
	\end{proposition}
	
	\subsubsection{Lagrangian coordinates and the corrected entropy equation}
	
	Let $Y(t,x)$ be the Lagrangian mass coordinate defined by
	\[
	Y_x=\rho_b,\qquad Y_t=-\rho_bu_b,\qquad Y(0,0)=0,
	\]
	and let $x=x_b(t,y)$ denote its spatial inverse, so that $(x_b)_y=v=1/\rho_b$ and $(x_b)_t=u_b$. The curve $x_c(t)=x_b(t,0)$ therefore satisfies \eqref{intro:contact-trajectory}. Since $\rho_b$ is uniformly positive and bounded, all spatial $L^p$ norms in $x$ and $y$ are uniformly equivalent. Evaluating $(v,u,\theta,p)$ at $x=x_b(t,y)$, the one-dimensional Navier--Stokes--Fourier system in mass coordinates reads
	\[
	\left\{
	\begin{aligned}
		&v_t-u_y=0,\\
		&u_t+p_y=\nu\left(\frac{u_y}{v}\right)_y,\\
		&c_v\theta_t+pu_y=\kappa\left(\frac{\theta_y}{v}\right)_y+\frac{\nu u_y^2}{v},
		\qquad p=\frac{R_g\theta}{v}.
	\end{aligned}
	\right.
	\]
	Setting
	\[
	c_p=\gamma c_v,\qquad \kappa_p=\frac{\kappa}{c_p},\qquad
	\tilde{\nu}=\nu+(\gamma-1)\kappa_p,
	\]
	we define the modified entropy variable $\tilde{\sigma}$ and thermal diffusivity $D$ by
	\[
	\sigma=c_v\log\theta+R_g\log v-S_-,\qquad
	\tilde{\sigma}=\sigma+\frac{\kappa_p u_y}{\theta},\qquad
	D=\frac{\kappa_p}{v},\qquad
	A_\kappa=\frac{\kappa_p(\nu-\kappa_p)}{v\theta}.
	\]
	Differentiating $\tilde{\sigma}$ with respect to $t$ and substituting the equations for $v_t,u_t,\theta_t$ {cancels the pure $\theta_{yy}$ terms in $\tilde\sigma_t-(D\tilde\sigma_y)_y$ and gives}
	\begin{equation}\label{ev:eta}
		\tilde{\sigma}_t=(D\tilde{\sigma}_y)_y+\partial_y(A_\kappa u_{yy})+Q+R_2,
	\end{equation}
	where $Q$ and $R_2$ are given by
	\[
	Q=\frac{\tilde{\nu}}{v\theta}u_y^2
	+\frac{\kappa_p c_v}{v\theta^2}\theta_y^2
	+\frac{\kappa_p R_g}{v^2\theta}v_y\theta_y
	\]
	and
	\begin{align*}
		R_2={}&-\frac{\kappa_p\nu}{v^2\theta}(v_yu_{yy}+u_yv_{yy})
		+\frac{\kappa_p(\kappa_p+\nu)}{v\theta^2}\theta_yu_{yy}
		-\frac{(\gamma-1)\kappa_p^2}{v\theta^2}u_y\theta_{yy}\\
		&+\frac{2\kappa_p\nu}{v^3\theta}u_yv_y^2
		+\frac{(\gamma-1)\kappa_p^2}{v^2\theta^2}u_yv_y\theta_y
		-\frac{2\kappa_p^2}{v\theta^3}u_y\theta_y^2
		-\frac{\kappa_p\nu}{c_vv\theta^2}u_y^3.
	\end{align*}
	
	\subsubsection{Nonnegative entropy production}
	
	In Eulerian coordinates, write $U_m=(\rho_b^{-1},u_b,\theta_b)$ and $V_m=(v_v,u_v,\theta_v)$ with $(\rho_v,u_v,\theta_v)=V(W)$ and $v_v=\rho_v^{-1}$. Since $\partial_y=v\partial_x$, we write $Q(t,y)$ in terms of a symmetric bilinear form evaluated at $x=x_b(t,y)$:
	\[
	Q(t,y)=B_0(U_m)[U_{m,x},U_{m,x}]\bigl(t,x_b(t,y)\bigr),
	\qquad
	B_0(U_m)[Z,Z]
	=\frac{\tilde{\nu}v}{\theta}Z_u^2
	+\frac{\kappa_p c_vv}{\theta^2}Z_\theta^2
	+\frac{\kappa_p R_g}{\theta}Z_vZ_\theta.
	\]
	Differentiating $\theta_vv_v^{\gamma-1}=\theta_-v_-^{\gamma-1}$ from \eqref{v3:profile-identities} gives $\theta_{v,x}/\theta_v=-(\gamma-1)v_{v,x}/v_v$. Since $R_g=(\gamma-1)c_v$, the two $\theta_{v,x}$ terms in $B_0(V_m)[V_{m,x},V_{m,x}]$ cancel:
	\[
	\frac{\kappa_p v_v\theta_{v,x}}{\theta_v}
	\left(c_v\frac{\theta_{v,x}}{\theta_v}+R_g\frac{v_{v,x}}{v_v}\right)=0.
	\]
	Evaluating $B_0$ at $V_m$ therefore leaves only the nonnegative viscous term
	\begin{equation}\label{ev:P}
		P(t,y)=\left(\frac{\tilde{\nu}v_v}{\theta_v}u_{v,x}^2\right)\bigl(t,x_b(t,y)\bigr)
		=\bigl(K(W)W_x^2\bigr)\bigl(t,x_b(t,y)\bigr)\ge 0,
	\end{equation}
	where $K(W)>0$ is a smooth positive function. Setting $\tilde{q}=U_m-V_m$ and
	\begin{equation}\label{ev:N}
		N=\partial_y(A_\kappa u_{yy})+R_2+
		\Bigl(B_0(U_m)[\tilde{q}_x,\tilde{q}_x]
		+2B_0(U_m)[V_{m,x},\tilde{q}_x]
		+\bigl(B_0(U_m)-B_0(V_m)\bigr)[V_{m,x},V_{m,x}]\Bigr)\circ x_b,
	\end{equation}
	we write \eqref{ev:eta} as
	\begin{equation}\label{ev:heat}
		\tilde{\sigma}_t=(D\tilde{\sigma}_y)_y+P+N.
	\end{equation}
	
	\subsubsection{Decay of the entropy sources}
	
	Set $m=7a/2-1/2=59/50$ and $m_1=9a/2-1/2=83/50$. By \eqref{af:bootstrap}, \eqref{af:Burgers-all}, and \eqref{ex:annulus}, {the perturbation terms in \eqref{ev:N} are estimated by \eqref{tab:planar-remainders}. For the pure profile terms, on $|x|\le C_0 T$ the $L^1$ norms of $W_{xxx}$ and $W_{xxxx}$ are bounded by $CT^{-5/4}$ and $CT^{-7/4}$ by \eqref{af:Burgers-all} and Cauchy--Schwarz. On each exterior annulus apply \eqref{ex:annulus} to $W-W_\pm$ on the corresponding side, retaining $\|\partial_x^{16}W\|_{L^2}\le CT^{-33/4}$. The remaining products of profile derivatives are bounded by Cauchy--Schwarz and \eqref{af:Burgers-all}.} In particular, $V_{m,x}\tilde{q}_x$ satisfies
	\[
	\|V_{m,x}\tilde{q}_x\|_{L^1}\le CT^{-(1+3a)/2}=CT^{-61/50},\qquad
	\bigl\|\partial_x(V_{m,x}\tilde{q}_x)\bigr\|_{L^1}\le C\bigl(T^{-1/2-5a/2}+T^{-5/4-3a/2}\bigr),
	\]
	both of which decay faster than $T^{-m}$ and $T^{-m_1}$. Thus $N(t,y)$ satisfies
	\begin{align}\label{ev:Nbudgets}
		&\|N(t)\|_{L^1}\le CT^{-m},\quad
		\|N(t)\|_{L^2}\le CT^{-7a/2},\quad
		\|N(t)\|_{L^\infty}\le CT^{-4a},\quad
		\|N_y(t)\|_{L^1}\le CT^{-m_1},\notag\\
		&\int_0^\infty(1+\log T)^k
		\bigl(\|N\|_{L^1}+\|N\|_{L^2}+\|N\|_{L^\infty}+\|N_y\|_{L^1}\bigr)\dd t
		=o_\delta(1)\qquad(k\ge 0).
	\end{align}
	For $P(t,y)$ in \eqref{ev:P}, Lemma~\ref{lem:planar-burgers} and the identity $\int_{\mathbb{R}}|W_xW_{xx}|\dd x=\|W_x\|_{L^\infty}^2$ give
	\begin{equation}\label{ev:Pbudgets}
		\|P(t)\|_{L^1}\le \frac{C\delta}{T},\qquad
		\|P(t)\|_{L^\infty}\le CT^{-2},\qquad
		\|P_y(t)\|_{L^1}=\bigl\|\partial_x\bigl(K(W)W_x^2\bigr)\bigr\|_{L^1_x}\le CT^{-2},
	\end{equation}
	together with the global time integral
	\[
	\int_0^\infty\bigl(\|P(t)\|_{L^\infty}+\|P(t)\|_{L^2}^2+\|P_y(t)\|_{L^1}\bigr)\dd t=o_\delta(1).
	\]
	
	\subsubsection{Duhamel decomposition and logarithmic bounds}
	
	Let $S_0(t,s)$ denote the solution operator of $z_t=(Dz_y)_y$, and define
	\[
	H_P(t)=\int_0^t S_0(t,s)P(s)\dd s,\qquad
	H_N(t)=S_0(t,0)\tilde{\sigma}(0)+\int_0^t S_0(t,s)N(s)\dd s.
	\]
	Since the initial profile $\overline{U}_E(0,x)$ is isentropic ($\sigma(0)=0$), we have $\tilde{\sigma}(0)=\kappa_p u_{0,y}/\theta_0$, whose $L^1$ and Sobolev norms are $O(\delta)$. By uniqueness of solutions to \eqref{ev:heat} and the parabolic maximum principle, we obtain the decomposition
	\begin{equation}\label{ev:decomposition}
		\sigma=H_P+\sigma_0,\qquad
		\sigma_0=H_N-\frac{\kappa_p u_y}{\theta},\qquad H_P\ge0.
	\end{equation}
	
	Because both $S_0(t,s)$ and its backward dual are $L^1$-contractive and satisfy the uniformly parabolic $L^2$ energy identity, the Nash--duality argument of Lemma~\ref{lem:planar-burgers} yields the smoothing bounds
	\begin{equation}\label{ev:scalar-smoothing}
		\|S_0(t,s)\|_{L^1\to L^2}+\|S_0(t,s)\|_{L^2\to L^\infty}\le C(t-s)^{-1/4},\qquad
		\|S_0(t,s)\|_{L^1\to L^\infty}\le C(t-s)^{-1/2}.
	\end{equation}
	Splitting the Duhamel integrals at $s=t/2$ and applying \eqref{ev:Nbudgets}--\eqref{ev:Pbudgets} together with $m=59/50>1$ yields, for $H_N$,
	\begin{align}\label{ev:HNpoint}
		&\|H_N(t)\|_{L^1}\le M_\delta\quad\left(\lim_{\delta\rightarrow 0}M_\delta=0\right),\notag\\
		&\|H_N(t)\|_{L^2}\le CM_\delta T^{-1/4}+CT^{3/4-m},\notag\\
		&\|H_N(t)\|_{L^\infty}\le CM_\delta T^{-1/2}+CT^{1/2-m},
	\end{align}
	and for $H_P$,
	\begin{equation}\label{ev:HPpoint}
		\|H_P(t)\|_{L^1}\le C\delta J,\qquad
		\|H_P(t)\|_{L^2}\le C\delta T^{-1/4}J,\qquad
		\|H_P(t)\|_{L^\infty}\le C\delta T^{-1/2}J.
	\end{equation}
	
	To bound $\int_0^t\|\sigma_0(s)\|_{L^\infty}^2\dd s$, we split the integral for $H_N(t)$ at $s=\max\{0,t-1\}$. By the $L^\infty$ contraction property of $S_0$, the integral over $[\max\{0,t-1\},t]$ is bounded in $L^2(0,T;L^\infty)$ by $\|N\|_{L^2(0,T;L^\infty)}=o_\delta(1)$, while Minkowski's inequality and \eqref{ev:scalar-smoothing} bound the integral over $[0,\max\{0,t-1\}]$ by
	\[
	\left\|\int_0^{\max\{0,t-1\}}S_0(t,s)N(s)\dd s\right\|_{L^2(0,T;L^\infty)}
	\le C\left(\int_0^T\|N(s)\|_{L^1}\dd s\right)\sqrt{\log(2+T)}.
	\]
	Combining this with $S_0(t,0)\tilde{\sigma}(0)$ and $\int_0^\infty\|u_y/\theta\|_{L^\infty}^2\dd t=o_\delta(1)$ gives
	\begin{equation}\label{ev:Xlog}
		\int_0^t\|\sigma_0(s)\|_{L^\infty}^2\dd s
		\le A_\delta+\varepsilon_\delta\log(2+t),
		\qquad \lim_{\delta\rightarrow 0}(A_\delta+\varepsilon_\delta)=0,
	\end{equation}
	as well as
	\[
	\int_0^t\|\sigma(s)\|_{L^\infty}^2\dd s
	\le A_\delta+C\delta^2\log^3(2+t)+\varepsilon_\delta\log(2+t).
	\]
	
	\subsubsection{Kernel bounds and entropy gradients}
	
	By Proposition~\ref{prop:planar-decay}, $D=\kappa_p/v$ satisfies
	\begin{equation}\label{ev:diffusivity-bounds}
		\|D_y(t)\|_{L^\infty}\le CT^{-2a},\quad
		\|D_{yy}(t)\|_{L^\infty}\le CT^{-3a},\quad
		K_D=\int_0^\infty\bigl(\|D_y\|_{L^\infty}^2+\|D_{yy}\|_{L^\infty}\bigr)\dd t\le C.
	\end{equation}
	Consequently, Lemma~\ref{lem:scalar-localization} below applies to $S_0(t,s)$, and integrating \eqref{new:scalar-gaussian} against $L^\infty$ and $L^2$ data gives
	\begin{equation}\label{ev:gradient-smoothing}
		\|\partial_yS_0(t,s)\|_{L^\infty\to L^\infty}\le C(t-s)^{-1/2},\qquad
		\|\partial_yS_0(t,s)\|_{L^2\to L^\infty}\le C(t-s)^{-3/4}.
	\end{equation}
	Applying the semigroup property on $[t/2,t]$ and using \eqref{ev:HNpoint}--\eqref{ev:HPpoint} at time $t/2$ together with $\|N(s)\|_{L^\infty}\le CT^{-4a}$ and $\|P(s)\|_{L^\infty}\le CT^{-2}$ on $[t/2,t]$, we obtain
	\begin{equation}\label{ev:HNy}
		\|H_{N,y}(t)\|_{L^\infty}\le \frac{CM_\delta}{T}+k_N(t),\qquad k_N\in L^1(0,\infty),
	\end{equation}
	and
	\begin{equation}\label{ev:HPy}
		\|H_{P,y}(t)\|_{L^\infty}\le \frac{C\delta\log(2+t)}{T}+k_P(t),\qquad k_P\in L^1(0,\infty).
	\end{equation}
	Since $w=H_{P,y}$ and $w=H_{N,y}$ solve $w_t=(Dw)_{yy}+P_y$ and $w_t=(Dw)_{yy}+N_y$, Kato's inequality and \eqref{ev:Nbudgets}--\eqref{ev:Pbudgets} give
	\[
	\sup_{t\ge 0}\bigl(\|H_{P,y}(t)\|_{L^1}+\|H_{N,y}(t)\|_{L^1}\bigr)=o_\delta(1).
	\]
	This completes the proof of Proposition~\ref{prop:entropy-decomposition}.
	
	\begin{lemma}\label{lem:scalar-localization}
		Suppose $0<d_{\min}\le D(t,y)\le d_{\max}$ and $K_D=\int_0^\infty(\|D_{yy}(t)\|_{L^\infty}+\|D_y(t)\|_{L^\infty}^2)\dd t<\infty$. Then the fundamental solution $K(t,s;y,z)$ of $h_t=(Dh_y)_y$ ($t>s\ge 0$) satisfies the Gaussian bounds
		\begin{align}
			&|K(t,s;y,z)|\le C(t-s)^{-1/2}e^{-|y-z|^2/[C(t-s)]},\notag\\
			&|\partial_yK(t,s;y,z)|\le C(t-s)^{-1}e^{-|y-z|^2/[C(t-s)]},
			\label{new:scalar-gaussian}
		\end{align}
		where the constant $C>0$ depends only on $d_{\min},d_{\max}$, and $K_D$.
	\end{lemma}
	
	\begin{proof}
		Both $h_t=(Dh_y)_y$ and the differentiated equation $w_t=(Dw)_{yy}$ preserve nonnegativity and contract the $L^1$ norm. Applying the exponential weight $u=e^{\lambda y}w$ ($\lambda\in\mathbb{R}$) to $w_t=(Dw)_{yy}$ and integrating by parts gives
		\[
		\frac{1}{2}\frac{\dd}{\dd t}\|u\|_{L^2}^2
		=-\int_{\mathbb{R}} D|u_y|^2\dd y
		+\int_{\mathbb{R}}\left(\frac{1}{2}D_{yy}-\lambda D_y+\lambda^2D\right)u^2\dd y.
		\]
		The weighted $L^1$ growth rate is bounded by $d_{\max}\lambda^2$ for the forward equation and by $\|D_{yy}\|_{L^\infty}+2|\lambda|\|D_y\|_{L^\infty}+d_{\max}\lambda^2$ {in reverse time for the terminal-value adjoint $-\phi_s=D(s,y)\phi_{yy}$}. The Nash--duality argument therefore yields weighted $L^1\to L^\infty$ and $L^2\to L^\infty$ bounds with factor $\exp(C\lambda^2(t-s)+CK_D)$. Setting $\tau=t-s$ and $h(\xi)=S_0(\xi,s)f$, the weighted $L^2$ energy inequality on $[s+\tau/3,s+2\tau/3]$ guarantees the existence of an intermediate time $s_1\in[s+\tau/3,s+2\tau/3]$ such that
		\[
		\|e^{\lambda y}h_y(s_1)\|_{L^2}^2
		\le C\tau^{-3/2}e^{C\lambda^2\tau+CK_D}\|e^{\lambda y}f\|_{L^1}^2.
		\]
		Solving $w_t=(Dw)_{yy}$ for $w=h_y$ from $s_1$ to $t$ and applying the weighted $L^2\to L^\infty$ bound on $[s_1,t]$ gives
		\[
		\bigl\|e^{\lambda y}\partial_yS_0(t,s)f\bigr\|_{L^\infty}
		\le C\tau^{-1}e^{C\lambda^2\tau+CK_D}\|e^{\lambda y}f\|_{L^1}.
		\]
		Evaluating this operator bound on a smooth approximation to a Dirac mass centered at $z$ with the shifted weight $e^{\lambda(y-z)}$ and optimizing over $\lambda=(y-z)/[2C(t-s)]$ proves \eqref{new:scalar-gaussian}.
	\end{proof}
	
	\subsection{Central acoustic estimates}
	\label{sec:central-energy-replacement}
	
	In the central region $\{|y|\le T^b\}$ ($b=23/40$), the forward 3-rarefaction wave has moved toward $y\sim c_*t$, so the backward and forward acoustic perturbations decay faster than $T^{-a}$. We use this faster decay to bound $\partial_y(\sigma s_y)$.
	
	Let $V(W)$ be the primitive approximate rarefaction profile and set $q_b=U_b-V(W)$. By Proposition~\ref{prop:planar-decay} and Lemma~\ref{lem:planar-burgers}, we have {for $0\le k\le21$ in $L^2$ and $0\le k\le20$ in $L^\infty$,}
	\begin{equation}
		\|\partial_x^kq_b(t)\|_{L^2}\le C_kT^{-(k+1/2)a},\qquad
		\|\partial_x^kq_b(t)\|_{L^\infty}\le C_kT^{-(k+1)a},\qquad a=\frac{12}{25},\label{ce:1}
	\end{equation}
	and
	\begin{equation}
		\|V_x(t)\|_{L^\infty}\le CT^{-1},\qquad
		\|V_x(t)\|_{L^2}\le CT^{-1/2},\qquad
		\|V_{xx}(t)\|_{L^2}+\|W_{xx}(t)\|_{L^2}\le CT^{-5/4}.\label{ce:2}
	\end{equation}
	In mass coordinates $y$ with $x=x_b(t,y)$ and $c=c_0(p)$ from \eqref{v3:coefficient-functions}, we set
	\[
	\left\{
	\begin{aligned}
		&r=u_b-I(p),\qquad s=u_b+I(p),\qquad I'(p)=c_0(p)^{-1},\qquad \sigma=S-S_-,\\
		&\omega(\sigma)=\frac{1+e^{-\sigma/c_p}}{2},\qquad
		b_-=-c_0(p)\omega(\sigma),\qquad a_0=\frac{c}{2}\bigl(e^{-\sigma/c_p}-1\bigr).
	\end{aligned}
	\right.
	\]
	By \eqref{v3:profile-identities}, $r$ and $\sigma$ vanish on $V(W)$. Setting $S_v(t,y)=s(V(t,x_b(t,y)))$ and $q_c=s-S_v$, we deduce from \eqref{ce:1} that
	\begin{equation}
		\|(r,q_c,\sigma)\|_{L^2}\le CT^{-a/2},\qquad
		\|(r,q_c,\sigma)\|_{L^\infty}\le CT^{-a},\qquad
		\|(r_y,q_{c,y},\sigma_y)\|_{L^2}\le CT^{-3a/2},\label{ce:3}
	\end{equation}
	together with $\|(b_-)_y\|_{L^\infty}\le CT^{-2a}$.
	
	\begin{proposition}
		\label{prop:central-acoustic}
		Fix $b=23/40$, and let
		\[
		A_c(t)=\bigl\|\bigl(r(t),q_c(t)\bigr)\bigr\|_{L^\infty(\{|y|\le c_0T^b\})}
		\]
		where $c_0>0$ is a fixed constant. Then
		\[
		\int_0^\infty T^{3/10}A_c(t)^2\dd t<\infty,
		\qquad
		\int_0^\infty T\bigl\|\partial_y(\sigma s_y)(t)\bigr\|_{L^\infty(\{|y|\le c_0T^b\})}\dd t<\infty,
		\]
		and the same integrability holds when multiplied by any fixed power of $\log(2+t)$.
	\end{proposition}
	
	\subsubsection{Acoustic equations in mass coordinates}
	
	In mass coordinates, the pressure and momentum equations read
	\[
	\left\{
	\begin{aligned}
		&p_t=-c_L^2u_y+(F_p)_y+Q_p,\\
		&u_t=-p_y+(F_u)_y,
	\end{aligned}
	\right.
	\]
	where
	\[
	F_u=\frac{\nu u_y}{v},\qquad F_p=\frac{(\gamma-1)\kappa\theta_y}{v^2},\qquad
	Q_p=(\gamma-1)\left(\frac{\kappa v_y\theta_y}{v^3}+\frac{\nu u_y^2}{v^2}\right).
	\]
	Applying the chain rule to $r=u-I(p)$ and $s=u+I(p)$ with $I'(p)=1/c$ and $c_L^2=c^2e^{-\sigma/c_p}$ gives
	\begin{equation}
		\left\{
		\begin{aligned}
			&r_t-cr_y
			=\left(\frac{c_L^2}{c}-c\right)u_y+(F_u)_y-\frac{(F_p)_y+Q_p}{c}
			=c\bigl(e^{-\sigma/c_p}-1\bigr)u_y+(g_-)_y+Q_-,\\
			&s_t+cs_y
			=-c\bigl(e^{-\sigma/c_p}-1\bigr)u_y+(g_+)_y-Q_-,
		\end{aligned}
		\right.
		\label{v3:chain-rule-acoustics}
	\end{equation}
	where $g_\pm=F_u\pm F_p/c$ and $Q_-=(1/c)_yF_p-Q_p/c$, and $a_0=\frac{c}{2}(e^{-\sigma/c_p}-1)$ satisfies
	\begin{equation}\label{v3:coupling-bound}
		a_0=-\frac{c\sigma}{2c_p}\int_0^1e^{-\tau\sigma/c_p}\dd\tau,
		\qquad |a_0|\le C|\sigma|.
	\end{equation}
	Substituting $u_y=(r_y+s_y)/2$ and $b_-=-c-a_0=-c\omega(\sigma)$ into \eqref{v3:chain-rule-acoustics} gives
	\begin{equation}
		\left\{
		\begin{aligned}
			&r_t+b_- r_y=(g_-)_y+Q_-+a_0s_y,\\
			&s_t-b_- s_y=(g_+)_y-Q_--a_0r_y.
		\end{aligned}
		\right.
		\label{ce:4}
	\end{equation}
	
	Let $g_\pm^v$ and $Q_-^v$ denote the corresponding Eulerian expressions evaluated at $V(W)$ and $x=x_b(t,y)$, for instance
	\[
	g_\pm^v=\left(\nu u_{v,x}\pm\frac{(\gamma-1)\kappa\rho_v\theta_{v,x}}{c_0(p_v)}\right)\circ x_b,
	\]
	and set $\widetilde{g}_\pm=g_\pm-g_\pm^v$ and $\widetilde{Q}_-=Q_--Q_-^v$. By \eqref{ce:1}--\eqref{ce:2}, we have
	\begin{equation}
		\|\widetilde{g}_\pm(t)\|_{L^2}\le CT^{-3a/2},\qquad
		\|(g_\pm^v)_y(t)\|_{L^2}+\|Q_-^v(t)\|_{L^2}\le CT^{-5/4},\qquad
		\|\widetilde{Q}_-(t)\|_{L^2}\le CT^{-7a/2}.\label{ce:5}
	\end{equation}
	Writing $s_v(W)=s(V(W))$ and letting $c_{E,v}$ denote the Eulerian sound speed of $V(W)$, the function $S_v(t,y)=s_v(W(t,x_b(t,y)))$ satisfies
	\begin{equation}
		(S_v)_t-b_-(S_v)_y
		=d_*s_v'(W)W_{xx}\circ x_b+\Delta_b\bigl(s(V)\bigr)_x\circ x_b,
		\qquad \Delta_b=u_b-u_v-c_{E,v}-b_- v=O(q_b),\label{ce:6}
	\end{equation}
	where $W=u_v+c_{E,v}$ implies that $\Delta_b$ vanishes when $U_b=V(W)$. Consequently, the right-hand side of \eqref{ce:6} is bounded in $L^2(\mathbb{R})$ by $C(T^{-5/4}+T^{-1-a/2})$. Subtracting \eqref{ce:6} from the second equation in \eqref{ce:4} and decomposing $a_0s_y=a_0q_{c,y}+a_0(S_v)_y$, we see from \eqref{v3:coupling-bound} and \eqref{ce:3} that for {$\|\chi_0\|_{L^\infty}\le3$},
	\begin{equation}
		\left|\int_{\mathbb{R}} \chi_0 r a_0(S_v)_y\dd y\right|
		\le C\|(S_v)_y\|_{L^\infty}\|r\|_{L^2}\|\sigma\|_{L^2}
		\le CT^{-1-a}.\label{ce:7}
	\end{equation}
	
	\subsubsection{Weighted energy in the central region}
	
	For $0<b<1$, set $\ell_c=T^b$ and define
	\[
	\chi_-=2+\tanh\left(\frac{y}{\ell_c}\right),\qquad
	\chi_+=2-\tanh\left(\frac{y}{\ell_c}\right),\qquad
	E_c(t)=\frac{1}{2}\int_{\mathbb{R}}\bigl(\chi_-r^2+\chi_+q_c^2\bigr)\dd y.
	\]
	Since $b_-\le -c_1<0$ and $\partial_y\chi_\mp=\pm\ell_c^{-1}\operatorname{sech}^2(y/\ell_c)$, integrating the transport terms $\int_{\mathbb{R}}\chi_-b_-rr_y\dd y$ and $-\int_{\mathbb{R}}\chi_+b_-q_cq_{c,y}\dd y$ by parts gives a positive integral on $\{|y|\le C_0\ell_c\}$. Setting
	\[
	B(t)^2=\int_{|y|\le C_0\ell_c}\bigl(r(t,y)^2+q_c(t,y)^2\bigr)\dd y
	\]
	and estimating the source terms via \eqref{ce:3}--\eqref{ce:7} (integrating $\widetilde{g}_\pm$ by parts once while keeping $(g_\pm^v)_y$ as an $L^2$ source), we obtain
	\begin{equation}
		E_c'(t)+cT^{-b}B(t)^2
		\le C\bigl(T^{-3a}+T^{-1-a}+T^{-a/2-5/4}+T^{-b-2a}\bigr).\label{ce:8}
	\end{equation}
	Multiplying \eqref{ce:8} by $T^{\beta_0}$ and using $E_c(t)\le CT^{-a}$ shows that
	\begin{equation}
		\int_0^\infty T^{\beta_0-b}B(t)^2\dd t<\infty\label{ce:9}
	\end{equation}
	for every exponent $\beta_0$ satisfying
	\begin{equation}
		0<\beta_0<\min\left\{a,\ 3a-1,\ \frac{a}{2}+\frac{1}{4},\ b+2a-1\right\}.\label{ce:10}
	\end{equation}
	
	\subsubsection{Interpolation and the central product estimate}
	
	Let $A_c(t)=\|(r(t),q_c(t))\|_{L^\infty(\{|y|\le c_0T^b\})}$. For {$m=8$ and $\chi_t(y)=\chi(y/T^b)$, with a fixed smooth $\chi$ equal to $1$ on $[-c_0,c_0]$ and supported in $[-C_0,C_0]$}, \eqref{ce:1} and $b>a$ give
	\[
	\bigl\|\partial_y^m\bigl(\chi_t(r,q_c)\bigr)\bigr\|_{L^2}\le CT^{-(m+1/2)a}.
	\]
	By one-dimensional Gagliardo--Nirenberg interpolation,
	\[
	A_c(t)^2\le C B(t)^{2-1/m}T^{-(m+1/2)a/m}.
	\]
	Applying H\"older's inequality in time with exponents $\bigl(\frac{2m}{2m-1},2m\bigr)$ against \eqref{ce:9} shows that
	\begin{equation}
		\int_0^\infty T^{\beta_1} A_c(t)^2\dd t<\infty\label{ce:12}
	\end{equation}
	whenever
	\begin{equation}
		\beta_1<\beta_0-b+a+\frac{a-(\beta_0-b)-1}{2m}.\label{ce:11}
	\end{equation}
	Choosing the explicit parameters
	\[
	a=\frac{12}{25},\qquad b=\frac{23}{40},\qquad
	\beta_0=\frac{87}{200},\qquad m=8,\qquad \beta_1=\frac{3}{10},
	\]
	we verify that $\beta_0=0.435<3a-1=0.44$ and that the right-hand side of \eqref{ce:11} equals $0.31625>0.30=\beta_1$. Here $b=0.575$ is larger than all central cutoff exponents used in Sections~\ref{sec:actual-incoming-phase}--\ref{sec:common-energy} ($0.52,0.525,0.53,0.54,0.55$).
	
	To prove the second bound in Proposition~\ref{prop:central-acoustic}, note that {for $\nu_0=W_-/4$}, the function $K(1+\nu_0t-x)^{-3}$ is a supersolution for $W-W_-$ on $x\le\nu_0t-L$, {and $x_b(t,y)\le\nu_0t/2+CT^b$ there for $|y|\le C_0T^b$, since $\|u_b\|_{L^\infty}<\nu_0/2$ and $(x_b)_y=v$ is bounded. On a unit neighborhood of this central image, $|W-W_-|\le CT^{-3}$ for large $t$ and $\|\partial_x^5W\|_{L^\infty}\le CT^{-3}$ by \eqref{af:Burgers-all}. Unit-interval interpolation and the chain rule give $S_v-S_{v,-}$ and its first four derivatives of order $O(T^{-3})$}. Thus $s_y$ and $s_{yy}$ differ from $q_{c,y}$ and $q_{c,yy}$ by $O(T^{-3})$ on $\{|y|\le C_0T^b\}$. Interpolating $q_c$ between $A_c(t)$ and $\|s_{yyyy}\|_{L^\infty}\le CT^{-5a}$ {on $\{|y|\le c_0T^b\}$ gives}
	\[
	\|q_{c,yy}\|_{L^\infty}\le CT^{-5a/2}A_c^{1/2}+CT^{-2b}A_c,\qquad
	\|q_{c,y}\|_{L^\infty}\le CT^{-5a/4}A_c^{3/4}+CT^{-b}A_c.
	\]
	Multiplying by $\|\sigma\|_{L^\infty}\le CT^{-a}$ and $\|\sigma_y\|_{L^\infty}\le CT^{-2a}$ shows that the dominant terms in $T\|\partial_y(\sigma s_y)\|_{L^\infty}$ are
	\[
	T^{1-7a/2}A_c(t)^{1/2},\qquad T^{1-13a/4}A_c(t)^{3/4}.
	\]
	By H\"older's inequality and \eqref{ce:12}, both terms are integrable on $[0,\infty)$ since $14a+\beta_1=7.02>7$ and $26a+3\beta_1=13.38>13$. Consequently,
	\begin{equation}
		\int_0^\infty T\bigl\|\partial_y(\sigma s_y)(t)\bigr\|_{L^\infty(\{|y|\le c_0T^b\})}\dd t<\infty,\label{ce:13}
	\end{equation}
	which completes the proof of Proposition~\ref{prop:central-acoustic}.
	
	\section{Backward acoustic estimates}
	\label{sec:backward-all}
	
	\subsection{A one-sided bound for the backward slope}
	\label{sec:global-slope-new}
	
	While the forward 3-rarefaction wave $s$ expands across the rarefaction fan with positive slope of order $T^{-1}$, the backward 1-acoustic wave $r$ is generated only by nonlinear and viscous couplings. However, its spatial derivative $r_y$ decays in $L^\infty(\mathbb{R})$ only at rate $O(T^{-2a})=O(T^{-24/25})$, which is not integrable in time. After setting $\tilde{r}=r-kU_y$, the rescaled derivative $w_-=e^{h_-}\tilde{r}_y$ satisfies a parabolic Riccati equation with quadratic term $-\kappa_0(p)w_-^2$, so its positive part $w_-^+$ satisfies the upper bound
	\[
	\limsup_{t\to\infty}\,t\sup_{y\in\mathbb{R}}w_-(t,y)^+\le\kappa_R^{-1}.
	\]
	
	Write $U=(v,u,\sigma)^\top$ in mass coordinates and retain $b_-=-c\omega(\sigma)$, $c=c_0(p)$, and $a_0=\frac{c}{2}(e^{-\sigma/c_p}-1)$. {On the compact positive state region $\Omega_\delta=\overline{\{U(t,y):t\ge0,\ y\in\mathbb R\}}$}, define
	\[
	\kappa_R=\inf_{\Omega_\delta}\frac{\gamma+1}{4v_0(p)\sqrt{c_0(p)}}>0.
	\]
	
	\begin{proposition}
		\label{thm:incoming-main}\label{thm:actual-incoming-slope}
		The variable $r=u-I(p)$ satisfies
		\begin{equation}
			\limsup_{t\to\infty}\,t\sup_{y\in\mathbb{R}}
			\bigl(\sqrt{c_0(p)}\,\omega(\sigma)r_y\bigr)^+\le\kappa_R^{-1},
			\label{inc:actual-slope}
		\end{equation}
		and the same upper bound holds with $r_y$ replaced by the corrected derivative $\tilde{r}_y$ defined in \eqref{inc:corrected-variable}.
	\end{proposition}
	
	\subsubsection{Exterior entropy localization}
	\label{sec:positive-production-exterior}
	
	To extend \eqref{ce:13} to $\mathbb{R}$, we first estimate $H_P$ and $H_N$ on $\{|y|\ge c_eT^\beta\}$. By \eqref{ev:P} and the upper barrier $K(1+|x|-C_ft)^{-4}$ for $W_x$ on $\{|x|\ge C_ft+L\}$, $P(t,y)$ satisfies
	\begin{align}
		&\|P(t)\|_{L^1}\le CT^{-1},\qquad \|P(t)\|_{L^\infty}\le CT^{-2},
		\label{eq:P-source-sizes}\\
		&P(t,y)\le C(1+|y|)^{-8},\qquad\text{for }|y|\ge2C_fT.
		\label{eq:P-source-tail}
	\end{align}
	
	\begin{lemma}
		\label{lem:production-exterior}
		For any $1/2<\beta<1$ and $c_e>0$, $H_P$ satisfies
		\begin{equation}
			\|H_P(t)\|_{L^\infty(\{|y|\ge c_eT^\beta\})}\le CT^{-\beta},
			\qquad
			\|(H_P)_y(t)\|_{L^\infty(\{|y|\ge c_eT^\beta\})}\le CT^{-1/2-\beta}.
			\label{eq:P-exterior}
		\end{equation}
	\end{lemma}
	
	\begin{proof}
		For $|y|\ge c_eT^\beta\gg\sqrt{T}$, we split the Duhamel integral $H_P(t)=\int_0^t S_0(t,s)P(s)\dd s$ at $s_*=\min\{t,|y|/(8C_f)\}$. For $s\in[s_*,t]$, the $L^\infty$ contraction property and \eqref{eq:P-source-sizes} give
		\[
		\int_{s_*}^t\|P(s)\|_{L^\infty}\dd s\le C\int_{s_*}^t(1+s)^{-2}\dd s\le \frac{C}{1+|y|}.
		\]
		For $s\in[0,s_*]$, we split the spatial convolution with the Gaussian kernel \eqref{new:scalar-gaussian} at $|z|=|y|/2$: the inner region $|z|\le |y|/2$ is bounded by $CT^{-1/2}(1+\log T)e^{-c|y|^2/T}$ since $\int_0^t\|P(s)\|_{L^1}\dd s\le C(1+\log T)$, while the outer region $|z|\ge |y|/2$ is bounded by $Cs_*(1+|y|)^{-8}$ via \eqref{eq:P-source-tail}. Thus
		\[
		H_P(t,y)\le \frac{C}{1+|y|}+CT^{-1/2}(1+\log T)e^{-c|y|^2/T}.
		\]
		Restarting the evolution at $t/2$ with the derivative kernel bound in \eqref{new:scalar-gaussian} similarly yields
		\[
		|(H_P)_y(t,y)|\le \frac{CT^{-1/2}}{1+|y|}+CT^{-1}(1+\log T)e^{-c|y|^2/T}+CT^{-3/2},
		\]
		which establishes \eqref{eq:P-exterior}.
	\end{proof}
	
	Next, set $b=23/40$ and $R_c=T^b$. {For $0\le s\le R_c/(4C_0)$ with $C_0$ sufficiently large, the enlarged exterior annuli lie outside the propagation cone for large $R_c$. Applying \eqref{ex:annulus} on each side to the difference from its far-field state, with $m_e=49/20$ and derivative order $10$, shows} that the third-derivative term in $N(s,y)$ satisfies
	\[
	\|U_{yyy}(s)\|_{L^1(\{|y|\ge R_c/2\})}
	\le CR_c^{-243/200}(1+s)^{-189/125}+CR_c^{-99/20},
	\]
	{while the exterior $L^2$ bounds for the derivatives of orders $j=1,2$ of either profile are \[ C R_c^{-(49/20)(1-j/10)}+C R_c^{-49/20-j}\le C R_c^{-49/25}.\] For every remaining product in \eqref{ev:N}, placing two derivative factors in $L^2$ and the others in $L^\infty$ gives a time integral bounded by $C R_c^{-73/25}$. Thus}
	\begin{equation}
		\int_0^{R_c/(4C_0)}\|N(s)\|_{L^1(\{|y|\ge R_c/2\})}\dd s\le CR_c^{-1/2}.
		\label{new:source-tail}
	\end{equation}
	{For $t\ge1$, $|\tilde\sigma(0,z)|\le C\delta(1+|z|)^{-4}$ by \eqref{main:burgers-data}. Splitting \eqref{new:scalar-gaussian} at $|z|=R_c/2$ gives for $j=0,1$
		\[\sup_{|y|\ge R_c}|\partial_y^j S_0(t,0)\tilde\sigma(0)(y)|\le C\delta T^{-(1+j)/2}\big(e^{-cR_c^2/T}+R_c^{-3}\big).\]
		The bounded initial interval follows from the smooth initial data.}
	Splitting the Duhamel formula for $H_N(t)$ at $s=R_c/(4C_0)$ and $s=t/2$, and using \eqref{new:source-tail} on $[0,R_c/(4C_0)]$, the tail bound
	\[
	\int_{R_c/(4C_0)}^\infty\|N(s)\|_{L^1}\dd s\le CR_c^{1-m}\qquad\left(m=\frac{59}{50}\right)
	\]
	on the middle interval, and $\|N(s)\|_{L^\infty}\le CT^{-4a}$ on $[t/2,t]$, we obtain
	\begin{equation}
		\begin{aligned}
			&\|H_N(t)\|_{L^\infty(\{|y|\ge R_c\})}
			\le C\bigl(T^{-1/2}R_c^{1-m}+T^{1-4a}\bigr),\\
			&\|(H_N)_y(t)\|_{L^\infty(\{|y|\ge R_c\})}
			\le C\bigl(T^{-1}R_c^{1-m}+T^{1/2-4a}\bigr).
		\end{aligned}
		\label{new:regular-exterior}
	\end{equation}
	Combining \eqref{ev:decomposition}, \eqref{eq:P-exterior}, and \eqref{new:regular-exterior} with $\|u_y/\theta\|_{L^\infty}\le CT^{-2a}$ and $\|(u_y/\theta)_y\|_{L^\infty}\le CT^{-3a}$ gives the exterior entropy decay
	\begin{equation}
		\|\sigma(t)\|_{L^\infty(\{|y|\ge R_c\})}\le CT^{-b},\qquad
		\|\sigma_y(t)\|_{L^\infty(\{|y|\ge R_c\})}\le CT^{-b-1/2}.
		\label{new:entropy-exterior}
	\end{equation}
	Since $\|s_y\|_{L^\infty}\le CT^{-2a}$ and $\|s_{yy}\|_{L^\infty}\le CT^{-3a}$, \eqref{new:entropy-exterior} implies
	\[
	\bigl\|\partial_y(\sigma s_y)(t)\bigr\|_{L^\infty(\{|y|\ge R_c\})}\le CT^{-403/200}.
	\]
	Combining this bound with \eqref{ce:13} gives
	\begin{equation}
		\int_0^\infty T\bigl\|\partial_y(\sigma s_y)(t)\bigr\|_{L^\infty(\mathbb{R})}\dd t<\infty.
		\label{inc:coupling-budget}
	\end{equation}
	
	\subsubsection{Diffusive correction and the Riccati equation}
	
	In terms of $U=(v,u,\sigma)^\top$, the planar Navier--Stokes--Fourier system in mass coordinates takes the form
	\[
	U_t+A_0U_y=(B_0U_y)_y+Q_0,
	\]
	where $Q_0=(0,0,Q_s)^\top$ and
	\[
	A_0=\begin{pmatrix}0&-1&0\\-\gamma p/v&0&p/c_v\\0&0&0\end{pmatrix},
	\qquad
	B_0=\begin{pmatrix}0&0&0\\0&\nu/v&0\\-\kappa(\gamma-1)/v^2&0&\kappa/(c_vv)\end{pmatrix}.
	\]
	We write $D_b=\partial_t+b_-\partial_y$ and denote differentiation with respect to $U$ by $D$.
	
	\begin{lemma}
		\label{inc:row}\label{inc:reduction}
		Define the row vector $k(U)$ and backward diffusion coefficient $d_-(U)$ by
		\[
		k=\left(\frac{\nu-(\gamma-1)\kappa_p}{2v},\ 0,\ \frac{(\gamma+1)\kappa_p-\nu}{2c_p}\right),\qquad
		d_-=\frac{\nu+(\gamma-1)\kappa_p}{2v}>0,
		\]
		and set the corrected backward Riemann variable
		\begin{equation}
			\tilde{r}=r-kU_y.
			\label{inc:corrected-variable}
		\end{equation}
		Then $\tilde{r}$ satisfies the scalar parabolic equation
		\begin{equation}
			D_b\tilde{r}=d_-\tilde{r}_{yy}+R,\qquad R=a_0s_y+R_0,
			\label{inc:reduced-Z}
		\end{equation}
		where $\|(R_0)_y(t)\|_{L^\infty}\le CT^{-5a}$ and $\int_0^\infty T\|R_y(t)\|_{L^\infty}\dd t<\infty$. In addition,
		\begin{equation}
			\bigl\|\partial_y^j(r-\tilde{r})(t)\bigr\|_{L^\infty}\le CT^{-(j+2)a},
			\qquad 0\le j\le2.
			\label{inc:correction-rates}
		\end{equation}
	\end{lemma}
	
	\begin{proof}
		With $E=e^{-\sigma/c_p}$, the differential of $r=u-I(p)$ is $Dr=(cE,1,-cvE/c_p)$. Direct multiplication gives
		\begin{equation}
			Dr\,B_0+k(A_0-b_-I)-d_-Dr=-a_0k.
			\label{inc:full-row-identity}
		\end{equation}
		Applying $D_b$ to $\tilde{r}=r-kU_y$ and substituting $D_bU=-(A_0-b_-I)U_y+(B_0U_y)_y+Q_0$ together with \eqref{inc:full-row-identity} cancels all off-diagonal second derivatives $U_{yy}$ at $a_0=0$ and leaves \eqref{inc:reduced-Z} with remainder
		\begin{align}
			R_0={}&-a_0kU_{yy}
			+Dr\,DB_0[U_y]U_y
			-d_-D^2r[U_y,U_y]
			+kDA_0[U_y]U_y\notag\\
			&+Dr\,Q_0-Dk[D_bU]U_y
			-k(B_0U_y)_{yy}-kQ_{0,y}
			+d_-(kU_y)_{yy}.
			\label{inc:row-remainder}
		\end{align}
		By \eqref{ce:1}, differentiating $R_0$ once {gives terms bounded by sums of $|U_{yyyy}|$, $|\sigma U_{yyy}|$, $|\sigma_yU_{yy}|$, $|U_y U_{yy}|$, $|U_y U_{yyy}|$, $|U_{yy}|^2$, $|U_y|^3$ and their products with bounded factors}, all of which are bounded in $L^\infty(\mathbb{R})$ by $CT^{-5a}$ (with $5a=12/5>2$). Finally, since $a_0=-c\sigma/(2c_p)+O(\sigma^2)$, we have
		\[
		\partial_y(a_0s_y)=-\frac{c}{2c_p}\partial_y(\sigma s_y)+O_{L^\infty}(T^{-5a}),
		\]
		which is integrable against $T$ by \eqref{inc:coupling-budget}.
	\end{proof}
	
	To eliminate the linear $s_y$ and $\sigma_y$ coefficients in the equation for $\tilde{r}_y$, we define
	\[
	h_-=\frac{1}{2}\log c+\log\omega(\sigma),\qquad
	w_-=e^{h_-}\tilde{r}_y.
	\]
	Setting $F_v=(B_0U_y)_y+Q_0$ and
	\[
	(V_r,V_s,V_\sigma)=\bigl(Dr\,F_v,\ Ds\,F_v,\ (0,0,1)F_v\bigr),
	\]
	we find that $(r,s,\sigma)$ and $h_-$ satisfy the transport identities
	\[
	\left\{
	\begin{aligned}
		&D_br=a_0s_y+V_r,\qquad D_bs=2b_-s_y-a_0r_y+V_s,\qquad D_b\sigma=b_-\sigma_y+V_\sigma,\\
		&2b_-(h_-)_s=(b_-)_s,\qquad b_-(h_-)_\sigma=(b_-)_\sigma,\qquad
		(b_-)_re^{-h_-}=\frac{\gamma+1}{4v_0\sqrt{c}}=\kappa_0(p).
	\end{aligned}
	\right.
	\]
	Consequently, differentiating \eqref{inc:reduced-Z} with respect to $y$ and multiplying by $e^{h_-}$ shows that $D_bh_--(b_-)_y$ cancels the $(b_-)_ss_y$ and $(b_-)_\sigma\sigma_y$ terms from $-(b_-)_y\tilde{r}_y$, while $-(b_-)_r\tilde{r}_y^2$ gives $-\kappa_0(p)w_-^2$:
	\begin{equation}
		(w_-)_t+\bigl(b_--(d_-)_y+2d_-(h_-)_y\bigr)(w_-)_y
		=d_-(w_-)_{yy}-\kappa_0(p)w_-^2+\epsilon_1w_-+e^{h_-}R_y,
		\label{inc:riccati}
	\end{equation}
	where the linear coefficient $\epsilon_1$ is given by
	\begin{align}
		\epsilon_1={}&-(b_-)_r(r_y-\tilde{r}_y)
		+a_0\bigl((h_-)_rs_y-(h_-)_sr_y\bigr)
		+(h_-)_rV_r+(h_-)_sV_s+(h_-)_\sigma V_\sigma\notag\\
		&+d_-\bigl((h_-)_y^2-(h_-)_{yy}\bigr)-(d_-)_y(h_-)_y,
		\qquad \|\epsilon_1(t)\|_{L^\infty}\le CT^{-3a}.
		\label{inc:linear-gradient-remainder}
	\end{align}
	
	\begin{proof}[Proof of Proposition~\ref{thm:incoming-main}]
		Let $M(t)=\sup_{y\in\mathbb{R}}w_-(t,y)^+$. Applying the parabolic maximum principle to \eqref{inc:riccati} yields the upper Dini derivative inequality
		\[
		D^+M(t)\le-\kappa_RM(t)^2+e(t)M(t)+h(t),
		\]
		where $e(t)=\|\epsilon_1(t)\|_{L^\infty}\in L^1(0,\infty)$ and $h(t)=\|e^{h_-}R_y(t)\|_{L^\infty}$ satisfies $\int_0^\infty Th(t)\dd t<\infty$ by Lemma~\ref{inc:reduction}. Setting
		\[
		H(t)=\int_t^\infty h(s)\dd s\in L^1(0,\infty),\qquad
		N(t)=M(t)+H(t),
		\]
		we find that
		\[
		N'(t)\le-\kappa_RN(t)^2+\bigl(e(t)+2\kappa_RH(t)\bigr)N(t)
		\]
		almost everywhere. With the integrating factor $A(t)=\exp\int_0^t(e+2\kappa_RH)\dd s\to A_\infty\in(0,\infty)$, {if $N(0)=0$, then $h=0$ almost everywhere and $M=0$. For $N(0)>0$, Bernoulli comparison gives}
		\[
		M(t)\le N(t)\le\frac{A(t)}{N(0)^{-1}+\kappa_R\int_0^t A(s)\dd s},
		\]
		{whence $\limsup_{t\to\infty}tM(t)\le\kappa_R^{-1}$. More precisely, $A$ is nondecreasing and $A(0)=1$, so
			\[
			\begin{aligned}
				\int_0^t\left(\kappa_R M(s)-\frac1{1+s}\right)^+\dd s
				&\le\int_0^t\left(\frac{A(s)}{\min\{1,(\kappa_RN(0))^{-1}\}+\int_0^sA(\tau)\dd\tau}
				-\frac1{1+s}\right)\dd s\\
				&\le\log A_\infty+\log\max\{1,\kappa_R N(0)\}.
			\end{aligned}
			\]
			Thus the integrable remainder in \eqref{inc:regular-slope} is controlled by $N(0)$, $\int_0^\infty e(t)\dd t$, $\int_0^\infty t h(t)\dd t$, and the integrable correction $CT^{-3a}$.} Since $\|r_y-\tilde{r}_y\|_{L^\infty}\le CT^{-3a}=o(T^{-1})$ by \eqref{inc:correction-rates}, this establishes \eqref{inc:actual-slope}.
	\end{proof}
	
	For the three-dimensional acoustic normalization in Section~\ref{sec:common-energy}, we define
	\[
	c_b(s)=c_0\bigl(I^{-1}(s/2)\bigr),\qquad a_b(s)=c_b(s)^{-1/2},\qquad
	\alpha=a_b(s)\tilde{r},\qquad
	\kappa_*=\frac{(\gamma+1)\sqrt{c_*}}{4v_-}.
	\]
	Since $\alpha_y=a_be^{-h_-}w_-+(a_b)_ss_y\tilde{r}$ with $\|(a_b)_ss_y\tilde{r}\|_{L^\infty}\le CT^{-3a}$ {and $\kappa_R=(\gamma+1)/(4v_-\sqrt{c_*})+o_\delta(1)$ by the uniform convergence of the state range to $U_-$, so that $\kappa_*a_be^{-h_-}/\kappa_R=1+o_\delta(1)$}, Proposition~\ref{thm:incoming-main} implies that for every fixed $\varepsilon_s>0$,
	\begin{equation}
		\kappa_*\alpha_y(t,y)
		\le\frac{1+\varepsilon_s+o_\delta(1)}{T}+k_s(t),
		\qquad k_s\ge0,\quad k_s\in L^1(0,\infty).
		\label{inc:regular-slope}
	\end{equation}
	
	\subsection{A monotone phase for the backward wave}
	\label{sec:actual-incoming-phase}
	
	While the upper bound \eqref{inc:regular-slope} controls the $s_-^2$ term in Lemma~\ref{lem:3d-signed-incoming}, the remaining backward terms involve $|\partial_1\alpha_E|$, which decays only like $O(T^{-2a})$. To control these terms without adding an $O(1)\log T$ growth to the energy, we construct a nondecreasing phase $\Psi$ with derivative $\psi=\Psi_y>0$ such that $I_\Psi=\int_{\mathbb{R}}(|(\alpha_\ell)_y|^2/\psi)\dd y$ belongs to $L^2(0,\infty)$.
	
	\subsubsection{The backward conservation equation and mass bound}
	
	Retaining $a=12/25$ and $U=(v,u,\sigma)^\top$ in mass coordinates $y$, we summarize from Proposition~\ref{prop:planar-decay} and \eqref{ce:1}--\eqref{ce:3} the planar input bounds
	\begin{align}
		&\|(r,\sigma)\|_{L^2}\le CT^{-a/2},& &
		\|(r,\sigma)\|_{L^\infty}\le CT^{-a},\notag\\
		&\|(r_y,\sigma_y)\|_{L^2}\le CT^{-3a/2},&
		&\|U_y\|_{L^2}\le CT^{-1/2},\notag\\
		&\|\partial_y^jU\|_{L^2}\le CT^{-(j+1/2)a}\quad(j=2,3),&
		&\|\partial_y^jU\|_{L^\infty}\le CT^{-(j+1)a}\quad(j=1,2,3),
		\label{eq:actual-physical-inputs}
	\end{align}
	together with $\int_0^\infty J^M\|U_{yyy}(t)\|_{L^1}\dd t<\infty$ for every fixed $M\ge 0$.
	
	\begin{lemma}
		\label{lem:entropy-fan}
		Let $F(t,y)=\partial_y(W(t,x_b(t,y)))\ge0$. Then
		\begin{equation}
			\int_0^t\!\int_{\mathbb{R}}|\sigma(s,y)|F(s,y)\dd y\dd s\le C(1+\log T),
			\label{eq:actual-entropy-fan-mass}
		\end{equation}
		and the same bound holds with $F$ replaced by $(S_v)_y$.
	\end{lemma}
	
	\begin{proof}
		Integrating \eqref{af:Mineq} {with $h_-,h_0$ as in \eqref{af:normalform}} over $[0,t]$ and applying \eqref{af:massresult}--\eqref{af:crossresult} gives
		\[
		\int_0^t\!\int_{\mathbb{R}} W_x\bigl(|h_-|+|h_0|\bigr)\dd x\dd s\le C(1+\log T).
		\]
		By the thermodynamic identity $\dd S=c_v\theta^{-1}\dd\theta-R_g\rho^{-1}\dd\rho$ and $c_v(\gamma-1)=R_g$, the differential $\dd S$ vanishes identically along both acoustic eigenvectors $(1,\pm c_E/\rho,(\gamma-1)\theta/\rho)^\top$. {Hence $|\sigma|\le C(|q_0|+|g|^2)$. By \eqref{af:normalform} and $G=q+Kq_x$,
			\[
			\begin{aligned}
				q_0&=h_0+\sum_{j\ne0}\frac{q_{0j}(W)}{\lambda_0-\lambda_j}\widetilde G_j^2-\chi_0 W_x-(Kq_x)_0,\\
				q_x&=Lg_x+L'(W)W_xg,\qquad
				|\widetilde G|\le C(|g|+|g_x|+W_x).
			\end{aligned}
			\]
			Since $W_x$ is nonnegative and bounded, we obtain}
		\[
		|\sigma|\le C\bigl(|h_0|+|g_x|+W_x+|g|^2+|g_x|^2\bigr).
		\]
		Multiplying by $W_x$ and integrating over $x\in\mathbb{R}$ and using $F\dd y=W_x\dd x$, we see that the remaining four integrals are bounded in time by $CT^{-(1+3a)/2}$, $CT^{-1}$, $CT^{-1-a}$, and $CT^{-1-3a}$, respectively, proving \eqref{eq:actual-entropy-fan-mass}.
	\end{proof}
	
	For $\alpha=a_b(s)\tilde{r}$ with $a_b(s)=c_b(s)^{-1/2}$, set $b_\alpha=b_-+(d_-)_y+2d_-\partial_y\log a_b\le -b_0<0$. Multiplying \eqref{inc:reduced-Z} by $a_b(s)$ puts the evolution of $\alpha$ into the conservative drift--diffusion form
	\begin{equation}
		\left\{
		\begin{aligned}
			&\alpha_t+(b_\alpha\alpha)_y=(d_-\alpha_y)_y+Q,\\
			&Q=a_bR+\bigl(D_ba_b+a_b(b_-)_y+(a_bd_-)_{yy}\bigr)\tilde{r}.
		\end{aligned}
		\right.
		\label{eq:actual-n-conservative}
	\end{equation}
	
	\begin{lemma}\label{lem:incoming-mass}
		With $b_m=1-2a=1/25$, the normalized backward variable $\alpha$ and source $Q$ satisfy
		\begin{align}
			\|\alpha(t)\|_{L^1}+\int_0^t\|Q(s)\|_{L^1}\dd s&\le CT^{b_m},
			\label{eq:actual-mass}\\
			\|\alpha(t)\|_{L^\infty}\le CT^{-a},\qquad
			\|\alpha_y(t)\|_{L^\infty}&\le CT^{-2a},\qquad
			\|\alpha_y(t)\|_{L^2}^2\le CT^{-3a}.
			\label{eq:actual-n-derivatives}
		\end{align}
		Consequently, $\int_0^\infty T^{-c}J^M\|Q(t)\|_{L^1}\dd t<\infty$ for every fixed $c>b_m$ and $M\ge 0$.
	\end{lemma}
	
	\begin{proof}
		In \eqref{inc:row-remainder}, every term of $R_0$ except $O(|U_y|^2)$ is integrable in $L^1(\mathbb{R})$ over $t\in[0,\infty)$, while $\int_0^t\|U_y(s)\|_{L^2}^2\dd s\le C(1+\log T)$, so $\int_0^t\|R_0(s)\|_{L^1}\dd s\le C(1+\log T)$. Decomposing $s_y=(S_v)_y+q_{c,y}$ in $a_0s_y$, Lemma~\ref{lem:entropy-fan} bounds $\int_0^t\|a_0(S_v)_y\|_{L^1}\dd s$ by $C(1+\log T)$, while
		\[
		\|a_0q_{c,y}\|_{L^1}\le C\|\sigma\|_{L^2}\|q_{c,y}\|_{L^2}\le CT^{-2a}.
		\]
		Substituting $D_bs=2b_-s_y-a_0r_y+V_s$ into $D_ba_b+a_b(b_-)_y$, we see that the coefficient $2b_-(a_b)_s+a_b(b_-)_s$ vanishes identically at $r=\sigma=0$, so its product with $s_y\tilde{r}$ has $L^1$ norm $O(T^{-3a})$, whereas the terms involving $r_y\tilde{r}$ and $\sigma_y\tilde{r}$ have $L^1$ norm $O(T^{-2a})$. Integrating in time yields $\int_0^t\|Q(s)\|_{L^1}\dd s\le CT^{1-2a}=CT^{1/25}$, and Kato's inequality applied to \eqref{eq:actual-n-conservative} proves \eqref{eq:actual-mass}.
	\end{proof}
	
	\subsubsection{The comparison function and quotient estimate}
	
	Let $m_-(t,y)\ge 0$ solve
	\[
	\left\{
	\begin{aligned}
		&(m_-)_t+(b_\alpha m_-)_y=\bigl(d_-(m_-)_y\bigr)_y+|Q|,\\
		&m_-(0,y)=|\alpha(0,y)|{+\delta\exp(-y^2)}.
	\end{aligned}
	\right.
	\]
	By the parabolic comparison principle applied to $m_-\pm\alpha$ and spatial integration, we have $m_-(t,y)\ge|\alpha(t,y)|$ pointwise and
	\[
	\|m_-(t)\|_{L^1}={\|m_-(0)\|_{L^1}+\int_0^t\|Q(s)\|_{L^1}}\dd s\le CT^{1/25}.
	\]
	Fix
	\begin{equation}
		\ell=T^{2/5},\qquad a_1=\eta_0 T^{-1/4},\qquad
		\beta_\Psi=\frac{27}{50},
		\label{eq:actual-phase-exponents}
	\end{equation}
	and let $G_\ell(y)$ and $H_\ell(y)$ denote centered Gaussian kernels with variances $\ell^2$ and $C_H\ell^2$ ($C_H>1$ fixed). We define $\alpha_\ell$, $\Psi$, and $\psi$ by
	\[
	\alpha_\ell=G_\ell*\alpha,\qquad
	\Psi(t,y)=a_1 H_\ell*\left(\int_{-\infty}^{\,\cdot}m_-(t,z)\dd z\right)(y),\qquad
	\psi=\Psi_y=a_1(H_\ell*m_-)>0.
	\]
	
	\begin{proposition}
		\label{prop:positive-phase}
		The backward weight phase $\Psi$ is strictly increasing in $y$, obeys $0\le\Psi\le C\eta_0 T^{-21/100}$ and $|(\alpha_\ell)_y|\le C_{\eta}\psi$, and solves the evolution equation
		\begin{equation}
			\Psi_t-c_L\Psi_y=(d_-\Psi_y)_y+(F_\Psi)_y+R_\Psi.
			\label{eq:actual-phase-equation}
		\end{equation}
		For every fixed $M\ge 0$, the phase flux $F_\Psi$, remainder $R_\Psi$, and slope $\psi$ satisfy
		\begin{equation}
			\int_0^\infty J^M\bigl(\|R_\Psi\|_{L^\infty}+\|F_\Psi\|_{L^\infty}^2
			+\|\psi\|_{L^\infty}^2+\|(F_\Psi)_y\|_{L^\infty}
			+\|\psi_y\|_{L^\infty}\bigr)\dd t<\infty,
			\label{eq:actual-phase-budgets}
		\end{equation}
		while the quotient integral $I_\Psi=\int_{\mathbb{R}}(|(\alpha_\ell)_y|^2/\psi)\dd y$ and mollification error satisfy
		\[
		I_\Psi(t)\le C_{\eta}T^{-51/100},\qquad
		\|\alpha(t)-\alpha_\ell(t)\|_{L^\infty}^2\le CT^{-26/25}.
		\]
	\end{proposition}
	
	\begin{proof}
		By Lemma~\ref{lem:incoming-mass} and $\|H_\ell\|_{L^\infty}\le C\ell^{-1}$, we have
		\[
		\|\psi\|_{L^\infty}\le Ca_1\ell^{-1}\|m_-\|_{L^1}\le C\eta_0 T^{-61/100},\qquad
		\|\psi_y\|_{L^\infty}\le C\eta_0 T^{-101/100},
		\]
		while $\|\alpha-\alpha_\ell\|_{L^\infty}^2\le C\ell\|\alpha_y\|_{L^2}^2\le CT^{-26/25}$. Since $C_H>1$ ensures the Gaussian domination $|G_\ell'(y)|\le C\ell^{-1}H_\ell(y)$ and $|G_\ell'(y)|^2/H_\ell(y)\le C\ell^{-2}H_\ell(y)$, the Cauchy--Schwarz inequality gives
		\[
		|(\alpha_\ell)_y|^2=|G_\ell'*\alpha|^2
		\le(H_\ell*m_-)\left(\frac{|G_\ell'|^2}{H_\ell}*|\alpha|\right),
		\]
		which implies both $|(\alpha_\ell)_y|\le C(a_1\ell)^{-1}\psi\le C_{\eta}\psi$ and
		\[
		I_\Psi(t)=\int_{\mathbb{R}}\frac{|(\alpha_\ell)_y|^2}{\psi}\dd y
		\le\frac{C\|\alpha(t)\|_{L^1}}{a_1\ell^2}
		\le C_{\eta}T^{1/25+1/4-4/5}=C_{\eta}T^{-51/100}.
		\]
		Integrating the equation for $m_-$ over $(-\infty,y]$ and convolving with $a_1 H_\ell$ yields \eqref{eq:actual-phase-equation} with explicit flux and source
		\[
		\left\{
		\begin{aligned}
			&F_\Psi=a_1\bigl[H_\ell*(d_-m_-)-d_-(H_\ell*m_-)\bigr]+C_H\ell\ell'\psi,\\
			&R_\Psi=\frac{a_1'}{a_1}\Psi
			+a_1 H_\ell*\left(\int_{-\infty}^{\,\cdot}|Q(t,z)|\dd z\right)
			+a_1\bigl[b_\alpha(H_\ell*m_-)-H_\ell*(b_\alpha m_-)\bigr]\\
			&\hspace{1.2cm}-(b_\alpha+c_L)\psi-a_1 H_\ell*\bigl((d_-)_ym_-\bigr).
		\end{aligned}
		\right.
		\]
		Since $\|(b_\alpha)_y\|_{L^\infty}+\|(d_-)_y\|_{L^\infty}+\|b_\alpha+c_L\|_{L^\infty}\le CT^{-2a}$ and the convolution commutator with any Lipschitz function $f$ satisfies
		\[
		\|H_\ell*(fm_-)-f(H_\ell*m_-)\|_{L^\infty}\le C\|f_y\|_{L^\infty}\|m_-\|_{L^1},
		\]
		we find
		\[
		\|F_\Psi\|_{L^\infty}\le CT^{-81/100},\qquad
		\|(F_\Psi)_y\|_{L^\infty}\le CT^{-121/100},\qquad
		\|R_\Psi\|_{L^\infty}\le Ca_1\|Q\|_{L^1}+CT^{-117/100},
		\]
		which together with Lemma~\ref{lem:incoming-mass} proves \eqref{eq:actual-phase-budgets}.
	\end{proof}
	
	\subsubsection{Central localization and phase flattening}
	
	So that the weight $\mathcal{W}$ in Section~\ref{sec:r5a-full-interface} is constant in $x_1$ across the central region, we modify $\Psi$ to be constant in $y$ on $\{|y|\le 2L_\Psi\}$, where $L_\Psi=T^{\beta_\Psi}=T^{27/50}$.
	
	\begin{lemma}
		\label{lem:phase-flattening}
		Let $L_\Psi=T^{27/50}$. Then for every fixed $C_0>0$ and $M\ge 0$,
		\begin{equation}
			\int_0^\infty J^M\frac{a_1}{L_\Psi}
			\int_{|y|\le C_0L_\Psi}m_-(t,y)\dd y\dd t<\infty.
			\label{eq:actual-central-occupation}
		\end{equation}
		Consequently, $\Psi$ can be replaced by a function $\widetilde{\Psi}$ that is constant in $y$ on $\{|y|\le 2L_\Psi\}$, agrees with $\Psi$ on $\{|y|\ge 3L_\Psi\}$, satisfies $\widetilde{\Psi}_y\ge 0$ everywhere and \eqref{eq:actual-phase-equation}--\eqref{eq:actual-phase-budgets}, and takes a central value $\psi_c(t)$ with $\int_0^\infty J^M|\psi_c'(t)|\dd t<\infty$.
	\end{lemma}
	
	\begin{proof}
		Multiplying $(m_-)_t+(b_\alpha m_-)_y=(d_-(m_-)_y)_y+|Q|$ by $J^Ma_1\Phi_{L_\Psi}(y)$, where $\Phi_{L_\Psi}(y)=\tanh(y/(C_0L_\Psi))$ satisfies $\Phi_{L_\Psi}'\ge c/L_\Psi$ on $\{|y|\le C_0L_\Psi\}$, and using $-b_\alpha\ge b_0>0$ yields
		\begin{align*}
			&\int_0^\infty J^M\frac{a_1}{L_\Psi}\int_{|y|\le C_0L_\Psi}m_-(t,y)\dd y\dd t\\
			&\quad\le C+C\int_0^\infty J^Ma_1
			\left[\|m_-\|_{L^1}\bigl(T^{-1}+L_\Psi^{-2}+L_\Psi^{-1}\|(d_-)_y\|_{L^\infty}\bigr)+\|Q\|_{L^1}\right]\dd t<\infty.
		\end{align*}
		Since $\ell=T^{2/5}=o(L_\Psi)$, the Gaussian tail of $H_\ell$ implies that $L_\Psi^{-1}\int_{|y|\le 3L_\Psi}\psi\dd y$ is likewise integrable against $J^M$. Let $\psi_c(t)$ be a smooth spatial average of $\Psi(t,\cdot)$ over $[-L_\Psi,L_\Psi]$, so that
		\[
		|\psi_c'(t)|\le C\left[
		\bigl(L_\Psi^{-1}+L_\Psi^{-2}+T^{-1}\bigr)\int_{|y|\le L_\Psi}\psi\dd y
		+L_\Psi^{-1}\|F_\Psi\|_{L^\infty}+\|R_\Psi\|_{L^\infty}\right]\in L^1_t(J^M\dd t).
		\]
		Choose a smooth cutoff $\chi_0(y/L_\Psi)$ vanishing on $|y|\le 2L_\Psi$ and equal to $1$ on $|y|\ge 3L_\Psi$ with $\operatorname{sgn}(y)\chi_0'(y/L_\Psi)\ge 0$, and set
		\[
		q_1=\Psi-\psi_c,\qquad \widetilde{\Psi}=\psi_c+\chi_0q_1.
		\]
		The monotonicity of $\Psi$ implies $\operatorname{sgn}(y)q_1\ge 0$ for $|y|\ge L_\Psi$, so $\widetilde{\Psi}_y=\chi_0\psi+(\chi_0)_yq_1\ge 0$ everywhere. The modified flux $\widetilde{F}_\Psi=\chi_0F_\Psi-d_-(\chi_0)_yq_1$ and remainder $\widetilde{R}_\Psi$ satisfy \eqref{eq:actual-phase-budgets} via \eqref{eq:actual-central-occupation}.
	\end{proof}
	
	Finally, fix $\beta_N=11/20>27/50=\beta_\Psi$, and let $\chi_-(y/(4T^{\beta_N}))$ be a smooth cutoff vanishing on $|y|\le 4T^{\beta_N}$ and equal to $1$ on $|y|\ge 8T^{\beta_N}$. Set
	\[
	\alpha_0=\chi_-\alpha_\ell,\qquad
	\alpha_r=\alpha-\alpha_0=(\alpha-\alpha_\ell)+(1-\chi_-)\alpha_\ell.
	\]
	On the support of $\alpha_0$, we have $|y|\ge 4T^{\beta_N}>3L_\Psi$, so $\widetilde{\Psi}=\Psi$ and $\widetilde{\Psi}_y=\psi$, and
	\begin{equation}
		\int_{\mathbb{R}}\frac{|(\chi_-)_y\alpha_\ell|^2}{\psi}\dd y
		\le Ca_1^{-1}T^{-2\beta_N}\|\alpha\|_{L^1}
		\le C_{\eta}T^{-81/100}.
		\label{eq:phase-cutoff-quotient}
	\end{equation}
	Since $\beta_N=11/20<b=23/40$, the central acoustic estimate of Proposition~\ref{prop:central-acoustic} together with Proposition~\ref{prop:positive-phase} ensures that $\int_0^\infty J^M\|\alpha_r(t)\|_{L^\infty}^2\dd t<\infty$. In Eulerian coordinates, setting $y=Y(t,x_1)$ and
	\begin{equation}\label{111}
		\left\{
		\begin{aligned}
			&\Psi_E(t,x_1)=\widetilde{\Psi}\bigl(t,Y(t,x_1)\bigr),\\
			&\psi_E(t,x_1)=\partial_1\Psi_E(t,x_1)=\rho_b\widetilde{\Psi}_y\bigl(t,Y(t,x_1)\bigr),\\
			&\alpha_E(t,x_1)=\frac{\alpha_0\bigl(t,Y(t,x_1)\bigr)}{v\bigl(t,Y(t,x_1)\bigr)},
		\end{aligned}
		\right.
	\end{equation}
	we find that the Eulerian quotient integral satisfies
	\[
	I_E(t)=\int_{\operatorname{supp}\alpha_E}\frac{|\partial_1\alpha_E|^2}{\psi_E}\dd x_1
	\le CI_\Psi(t)+C_{\eta}\bigl(T^{-81/100}+T^{-163/100}\bigr),
	\]
	so that $\int_0^\infty J^MI_E(t)^2\dd t<\infty$ and $|\partial_1\alpha_E|\le C_\delta\psi_E$.
	
	\section{Forward acoustic comparison}
	\label{sec:forward-all}
	
	\subsection{Reduction to a scalar equation}
	\label{sec:outgoing-coordinate}
	
	We now turn to the forward 3-rarefaction variable $s=u+I(p)$ in mass coordinates $P=(p,u,\sigma)^\top$. Along the isentropic 3-rarefaction curve $r=\sigma=0$, the pressure, specific volume, and acoustic impedance are given explicitly as functions of $s$ by
	\begin{equation}\label{v3:scalar-coefficients}
		\left\{
		\begin{aligned}
			&p_b(s)=p_-\left(1+\frac{\gamma-1}{4c_-}s\right)^{2\gamma/(\gamma-1)},\\
			&v_b(s)=v_-\left(1+\frac{\gamma-1}{4c_-}s\right)^{-2/(\gamma-1)},\qquad
			c_b(s)=c_*\left(1+\frac{\gamma-1}{4c_-}s\right)^{(\gamma+1)/(\gamma-1)},
		\end{aligned}
		\right.
	\end{equation}
	whose derivatives with respect to $s$ satisfy
	\begin{equation}\label{v3:scalar-derivatives}
		p_b'(s)=\frac{c_b(s)}{2},\qquad v_b'(s)=-\frac{1}{2c_b(s)},\qquad
		c_b'(s)=\frac{\gamma+1}{4v_b(s)}>0.
	\end{equation}
	Thus the full thermodynamic state is parameterized by $(r,s,\sigma)$ via
	\[
	p=p_b(s-r),\qquad v=v_b(s-r)e^{\sigma/c_p},\qquad c=c_b(s-r),
	\]
	while the mass-coordinate system and acoustic equations read
	\begin{equation}\label{r5x:physical}
		\left\{
		\begin{aligned}
			&p_t=-c_L^2u_y+(F_p)_y+Q_p,\\
			&u_t=-p_y+(F_u)_y,\\
			&\sigma_t=(F_s)_y+Q_s,
		\end{aligned}
		\right.
	\end{equation}
	with fluxes and quadratic sources
	\[
	F_p=\frac{(\gamma-1)\kappa\theta_y}{v^2},\qquad
	F_u=\frac{\nu u_y}{v},\qquad
	F_s=\frac{\kappa\theta_y}{v\theta},\qquad
	Q_s=\frac{\kappa\theta_y^2}{v\theta^2}+\frac{\nu u_y^2}{v\theta},
	\]
	and the characteristic wave system
	\begin{equation}\label{r5x:acoustic}
		\left\{
		\begin{aligned}
			&r_t-cr_y=(g_-)_y+Q_-+a_0(r_y+s_y),\\
			&s_t+cs_y=(g_+)_y-Q_--a_0(r_y+s_y).
		\end{aligned}
		\right.
	\end{equation}
	
	Because $c=c_b(s-r)$ depends on $r$ and the quadratic source $Q_-$ contains a non-integrable $O(s_y^2)$ term along the rarefaction wave, one cannot compare $s(t,y)$ directly with a scalar Burgers solution in the original coordinate $y$. In the following proposition, we construct an approximate entropy pair $(\widetilde{E},J_0)$ that shifts the spatial coordinate to $z=Z_+(t,y)$ so that the transport speed becomes a function of $s$ alone, while applying a nonlinear rescaling $w=h_0(s)-C_+$ that cancels the pure $s_y^2$ source identically along $r=\sigma=0$.
	
	\begin{proposition}\label{prop:outgoing-reduction}
		There exist a smooth strictly increasing function $h_0(s)$, a uniformly invertible spatial coordinate transformation $z=Z_+(t,y)$ with $Z_+(0,y)=y$, and a first-order derivative correction $C_+=k_1r_y+k_2\sigma_y$ with $C_+(0,y)=0$ such that the transformed forward variable $w=h_0(s)-C_+$ satisfies the scalar conservation equation \eqref{r5x:w-equation} in $(t,z)$, with uniformly positive diffusion $A_z$, strictly convex flux $f_0(w)$, and a source term $G_0-ka_0r_y$ whose pure $s_y^2$ component vanishes identically along $r=\sigma=0$.
	\end{proposition}
	
	\begin{proof}
		\noindent\textit{Step 1: Construction of the shifted spatial coordinate $z=Z_+(t,y)$.}
		{For $\xi$ in a compact subinterval of $\{1+(\gamma-1)\xi/(4c_-)>0\}$ containing the range of $s$}, define
		\begin{equation}\label{r5x:entropy-pair}
			a_*(\xi)=\frac{c_b'(\xi)}{2\sqrt{c_b(\xi)}},\qquad
			E_0=\frac{a_*(\xi)r}{\sqrt{c_b(s)}},\qquad
			J_0=a_*(\xi)\left(-\sqrt{c_b(s)}\,r+\frac{c_b'(s)}{2\sqrt{c_b(s)}}r^2\right).
		\end{equation}
		By Taylor expansion of $c=c_b(s-r)$ around $r=0$, we have at $\xi=s$:
		\[
		J_{0,r}+cE_{0,r}=O(r^2),\qquad
		J_{0,s}-cE_{0,s}=O(r^2),\qquad
		c(1+E_0)-J_0-c_b(s)=O(r^2).
		\]
		Setting $I_*'(p)=-c(p)^{-3/2}$, we define
		\begin{equation}\label{r5x:full-entropy}
			\widetilde{E}(p,u,\sigma;\xi)=E_0+a_*(\xi)\bigl(e^{\sigma/c_p}-1\bigr)
			\bigl(I_*(p)-I_*(p_b(\xi))\bigr)+\tanh\left(\frac{\sigma}{2c_p}\right).
		\end{equation}
		Along the Euler part of \eqref{r5x:physical}, the inviscid remainder of $(\widetilde{E},J_0)$ is
		\[
		R_E=(J_{0,r}+cE_{0,r})r_y+(J_{0,s}-cE_{0,s})s_y
		+(c^2-c_L^2)\bigl(E_{0,p}-a_*(\xi)I_*'(p)\bigr)u_y,
		\]
		which vanishes to second order in $(r,\sigma)$, satisfying $|\partial_\xi^jR_E|\le C_j(r^2+|r\sigma|)|P_y|$. With $F=(F_p,F_u,F_s)^\top$, define
		\begin{equation}\label{r5x:coordinate-definition}
			\left\{
			\begin{aligned}
				&F_E=D\widetilde{E}\cdot F,\qquad
				Q_E=\widetilde{E}_pQ_p+\widetilde{E}_\sigma Q_s-D^2\widetilde{E}[P_y,F],\\
				&M(t,y;\xi)=\int_{-\infty}^y\widetilde{E}(t,x;\xi)\dd x,\qquad
				\Sigma(t,y;\xi)=\int_{-\infty}^y(Q_E+R_E)(t,x;\xi)\dd x,
			\end{aligned}
			\right.
		\end{equation}
		so that $M_t=-J_0+F_E+\Sigma$. Since $|\partial_\xi^j\widetilde{E}|\le C_j(|r|+|\sigma|)$, the planar $L^1$ bound \eqref{af:physicalmass} implies
		\begin{equation}\label{r5x:M-bound}
			\|\partial_\xi^jM(t,\cdot;\xi)\|_{L^\infty}\le C_j\bigl(\varepsilon_\delta+\delta\log(2+t)\bigr).
		\end{equation}
		Evaluating the parameter at $\xi=s(t,y)$, we define the shifted spatial coordinate
		\[
		Z_+(t,y)=y+M\bigl(t,y;s(t,y)\bigr),\qquad
		J_z=(Z_+)_y=1+\widetilde{E}+M_\xi s_y.
		\]
		By the chain rule and the algebraic identity $(1+\tanh(\sigma/(2c_p)))(1+e^{-\sigma/c_p})=2$, we obtain the coordinate velocity relation
		\begin{equation}\label{r5x:coordinate-velocity}
			\left\{
			\begin{aligned}
				&(Z_+)_t+(c+a_0)J_z=c_b(s)+\Sigma+\Lambda,\\
				&\Lambda=\Delta V+F_E+M_\xi\bigl((g_+)_y-Q_--a_0r_y\bigr),\\
				&\Delta V=(c+a_0)(1+\widetilde{E})-J_0-c_b(s)=O(r^2+|r\sigma|),
			\end{aligned}
			\right.
		\end{equation}
		with $|J_z-1|+|J_{z,y}|=o_\delta(1)$ and $Z_+(0,y)=y$.
		
		\medskip
		\noindent\textit{Step 2: Nonlinear rescaling and cancellation of the resonant $s_y^2$ source.}
		Along $r=\sigma=0$, the viscous flux $g_+=F_u+F_p/c$ expands as $g_+=d_bs_y+e_br_y+e_\sigma\sigma_y$ with
		\[
		d_b(s)=\frac{\nu+(\gamma-1)\kappa_p}{2v_b(s)},\qquad
		e_b(s)=\frac{\nu-(\gamma-1)\kappa_p}{2v_b(s)},\qquad
		e_\sigma(s)=\frac{\kappa_p c_b(s)}{c_p}.
		\]
		Since $c_b(s)v_b(s)=c_-+\frac{\gamma-1}{4}s$, we define the rescaling multiplier $k(s)$, transformed variable $h_0(s)$, and convex scalar flux $f_0(h_0(s))$ by
		\begin{equation}\label{r5x:k}
			\frac{k'(s)}{k(s)}=\frac{3\gamma-5}{8c_b(s)v_b(s)},\qquad
			f_0'\bigl(h_0(s)\bigr)=c_b(s),\qquad
			f_0''\bigl(h_0(s)\bigr)=\frac{c_b'(s)}{k(s)}>0,
		\end{equation}
		or explicitly,
		\[
		k(s)=\left(1+\frac{\gamma-1}{4c_-}s\right)^{(3\gamma-5)/[2(\gamma-1)]},\qquad
		h_0(s)=\int_0^s k(z)\dd z,\qquad
		f_0\bigl(h_0(s)\bigr)=\int_0^s c_b(z)k(z)\dd z.
		\]
		Setting
		\[
		\left\{
		\begin{aligned}
			&k_1=\frac{ke_b}{2c_b},\qquad k_2=\frac{ke_\sigma}{c_b},\qquad
			C_+=k_1r_y+k_2\sigma_y,\\
			&w=h_0(s)-C_+,\qquad
			R_s=(g_+)_y-Q_--a_0(r_y+s_y),
		\end{aligned}
		\right.
		\]
		direct differentiation of \eqref{r5x:acoustic} yields
		\begin{equation}
			\bigl[\partial_t+(c+a_0)\partial_y-\partial_y(d_b\partial_y)\bigr]w
			=(F_{1,r})_y+(F_{2,r})_{yy}+G_r-ka_0r_y,
			\label{r5x:raw-equation}
		\end{equation}
		where the flux and source remainders are given by
		\begin{align}
			&F_{2,r}=-k_1(g_--d_br_y)-k_2(F_s-d_b\sigma_y),\notag\\
			&F_{1,r}=k(g_+-d_bs_y)-2k_1c r_y-k_2 c\sigma_y
			+2(k_1)_yg_-+2(k_2)_yF_s-(d_b)_yC_+\notag\\
			&\quad\quad\;\;-k_1Q_--k_2 Q_s-k_1a_0(r_y+s_y)-a_0C_+,\notag\\
			&G_r=-k'(s)s_yg_++\bigl(2c(k_1)_y+k_1c_y\bigr)r_y
			+(k_2 c)_y\sigma_y-kQ_-\notag\\
			&\quad\quad\;\;-(k_1)_{yy}g_--(k_2)_{yy}F_s+(k_1)_yQ_-+(k_2)_yQ_s
			+(k_1)_ya_0(r_y+s_y)\notag\\
			&\quad\quad\;\;-k_1'(s)R_s r_y-k_2'(s)R_s\sigma_y+(a_0)_yC_+.
			\label{r5x:raw-fluxes}
		\end{align}
		Changing spatial coordinates from $y$ to $z=Z_+(t,y)$ via \eqref{r5x:coordinate-velocity} with $A_z=d_bJ_z^2$ and
		\begin{align}
			D_0&=f_0'(w)-c_b(s)-\Lambda-d_bJ_{z,y},\qquad
			\tilde{F}_1=J_zF_{1,r}-3J_zJ_{z,z}F_{2,r},\qquad \tilde{F}_2=J_z^2F_{2,r},\notag\\
			G_0&=G_r-J_{z,z}F_{1,r}+(J_{z,z}^2+J_zJ_{z,zz})F_{2,r}+D_0w_z,
			\label{r5x:pushed-fluxes}
		\end{align}
		transforms \eqref{r5x:raw-equation} into the scalar equation
		\begin{equation}
			w_t+f_0(w)_z=(A_zw_z)_z-\Sigma w_z
			+(\tilde{F}_1)_z+(\tilde{F}_2)_{zz}+G_0-ka_0r_y.
			\label{r5x:w-equation}
		\end{equation}
		Finally, along $r=\sigma=0$ at $\xi=s$, we have
		\[
		Q_-=-\frac{(\gamma-1)(2\nu+(\gamma-1)\kappa_p)}{8c_bv_b^2}s_y^2,\qquad
		F_E=\left(\frac{(\gamma+1)e_b}{8c_bv_b}+\frac{\kappa_p(\gamma-1)}{4c_bv_b^2}\right)s_y.
		\]
		{Since $w_z=(ks_y-(C_+)_y)/J_z=ks_y$ on profiles with $r(\cdot)=\sigma(\cdot)\equiv0$}, the quadratic terms in $G_0$ are $-k'(s)d_bs_y^2-kQ_--F_Eks_y$. Their coefficient of $s_y^2$ is
		\[
		-d_bk'
		+\frac{k(\gamma-1)(2\nu+(\gamma-1)\kappa_p)}{8c_bv_b^2}
		-k\left(\frac{(\gamma+1)e_b}{8c_bv_b}
		+\frac{\kappa_p(\gamma-1)}{4c_bv_b^2}\right)
		=-d_b\left[k'-\frac{3\gamma-5}{8c_bv_b}k\right]=0
		\]
		by \eqref{r5x:k}. This completes the proof of Proposition~\ref{prop:outgoing-reduction}.
	\end{proof}
	
	\subsection{Logarithmic control of the forward residual}
	\label{sec:outgoing-energy-replacement}
	
	Using the scalar reduction \eqref{r5x:w-equation} of Proposition~\ref{prop:outgoing-reduction}, we now decompose the forward variable $s$ into a monotone scalar viscous profile $\bar{s}$ (whose Eulerian characteristic speed slope is bounded above by $(1+o_\delta(1))/T$) and a perturbation $s-\bar{s}$ whose $L^\infty$ square has a small $o_\delta(1)\log(2+t)$ time integral.
	
	\begin{proposition}
		\label{prop:outgoing-amplitude}
		There exists a strictly increasing reference profile $\bar{s}(t,y)$ with $\bar{s}(0,y)=s(0,y)$ such that
		\[
		\int_0^t\|s(\tau)-\bar{s}(\tau)\|_{L^\infty}^2\dd\tau
		\le C_\delta+\varepsilon_\delta\log(2+t),\qquad
		\lim_{\delta\rightarrow 0}\varepsilon_\delta=0,
		\]
		such that
		\[
		\partial_1\lambda_+^E\bigl(\bar{s}(t,Y)\bigr)\le\frac{1+o_\delta(1)}{T}+k_+(t),\qquad k_+\in L^1(0,\infty),
		\]
		and for any $0<b<1$ and $C_0>0$, $\bar{s}-\bar{s}_-=O(T^{-3})$ and $\bar{s}_y=O(T^{-4})$ on $\{|y|\le C_0T^b\}$.
	\end{proposition}
	
	Define the scalar diffusivity $d_+(h_0(s))=d_b(s)$, and let $\Omega(t,z)$ solve the scalar viscous conservation law with drift $\Sigma$:
	\begin{equation}\label{eq:out-reference-new}
		\left\{
		\begin{aligned}
			&\Omega_t+\bigl(f_0'(\Omega)+\Sigma\bigr)\Omega_z=\bigl(d_+(\Omega)\Omega_z\bigr)_z,\\
			&\Omega(0,z)=h_0\bigl(s(0,z)\bigr).
		\end{aligned}
		\right.
	\end{equation}
	Subtracting \eqref{eq:out-reference-new} from \eqref{r5x:w-equation} shows that the difference $q_+=w-\Omega$ satisfies
	\begin{equation}\label{eq:out-error-ordinary}
		\left\{
		\begin{aligned}
			&q_{+,t}+\bigl(f_0(\Omega+q_+)-f_0(\Omega)+\Sigma q_+\bigr)_z
			=\bigl(A_zq_{+,z}+B(\Omega,q_+)\Omega_z\bigr)_z+S_0+\Sigma_zq_+,\\
			&q_+(0,z)=0,\qquad B(\Omega,q_+)=d_+(\Omega+q_+)-d_+(\Omega),
		\end{aligned}
		\right.
	\end{equation}
	where
	\[
	S_0=(F_1)_z+(F_2)_{zz}+N_0,\qquad
	F_1=\tilde{F}_1+\bigl(A_z-d_+(w)\bigr)\Omega_z,\qquad
	F_2=\tilde{F}_2,\qquad N_0=G_0-ka_0r_y.
	\]
	
	\begin{lemma}\label{lem:outgoing-source}
		Setting $n_+(t)=\|S_0(t)\|_{L^1_z}$, the source $S_0$ and drift $\Sigma$ satisfy
		\begin{align}
			&\int_0^\infty T^{1/20}n_+(t)\dd t=o_\delta(1),\qquad
			n_+(t)\le CT^{-2a},\label{eq:out-ordinary-budget-new}\\
			&\int_0^\infty T^{1/20}\bigl(\|\Sigma_z(t)\|_{L^2}+\|\Sigma_z(t)\|_{L^\infty}\bigr)\dd t=o_\delta(1),\qquad
			\sup_{t\ge 0}\|\Sigma(t)\|_{L^\infty}=o_\delta(1).
			\label{eq:out-drift-new}
		\end{align}
	\end{lemma}
	
	\begin{proof}
		Writing $\widetilde{P}=P-P^v$ in mass coordinates, Lemma~\ref{lem:exterior-L1} and Corollary~\ref{cor:planar-mixed-interactions} supply the time-weighted $L^1$ bounds
		\begin{align}
			&\int_0^\infty T^{1/20}\bigl(\|P_{yyy}\|_{L^1}+\|P_{yyyy}\|_{L^1}\bigr)\dd t=o_\delta(1),
			\label{eq:out-high-L1-input}\\
			&\int_0^\infty T^{1/20}\bigl(\|r\sigma_y\|_{L^1}+\|\sigma r_y\|_{L^1}\bigr)\dd t=o_\delta(1).
			\label{eq:out-mixed-L1-input}
		\end{align}
		For the coordinate remainder $R_E$, we have
		\[
		\bigl\|(r^2+\sigma^2+|r\sigma|)P_y\bigr\|_{L^1}
		\le C\|P_y^v\|_{L^\infty}\|(r,\sigma)\|_{L^2}^2
		+C\|(r,\sigma)\|_{L^\infty}\|(r,\sigma)\|_{L^2}\|\widetilde{P}_y\|_{L^2}
		\le C\bigl(T^{-1-a}+T^{-3a}\bigr),
		\]
		which together with $Q_E=O(|P_y|^2)$ and $\Sigma_z=(Q_E+R_E+\Sigma_\xi s_y)/J_z$ implies
		\[
		\|\Sigma_z(t)\|_{L^\infty}\le C\bigl(T^{-4a}+T^{-1-2a}\bigr),\qquad
		\|\Sigma_z(t)\|_{L^2}\le C\bigl(T^{-2a-1/2}+T^{-3/2}\bigr),
		\]
		proving \eqref{eq:out-drift-new}. Similarly, in $S_0=(F_1)_z+(F_2)_{zz}+N_0$, the resonant pure $s_y^2$ term has been canceled by Proposition~\ref{prop:outgoing-reduction}, and {$w_z-ks_y=-[(C_+)_y+(J_z-1)ks_y]/J_z$ is retained in the remainder. The term $M_\xi s_y$ in $J_z-1$ contributes at most $C(1+\log T)\|P_y^3\|_{L^1}$ to the product with $F_E$. The} remaining derivative products obey
		\[
		\begin{aligned}
			&\|(r_y,\sigma_y)P_y\|_{L^1}\le CT^{-1/2-3a/2},& &
			\|(r,\sigma)P_{yy}\|_{L^1}\le CT^{-3a},\\
			&\|P_yP_{yy}\|_{L^1}\le CT^{-1/2-5a/2},& &
			\|P_y^3\|_{L^1}\le CT^{-1-2a},\\
			&\|P_yP_{yyy}\|_{L^1}\le CT^{-1/2-7a/2},& &
			\|P_{yy}^2\|_{L^1}\le CT^{-5a},
		\end{aligned}
		\]
		every exponent of which exceeds $21/20$. Since $\widetilde{E}(p_b(s),s/2,0;\xi)=0$ identically in $s$, differentiating with respect to $s$ gives $\widetilde{E}_pp_b'(s)+\frac{1}{2}\widetilde{E}_u=0$ on $r=\sigma=0$, so the Jacobian difference satisfies
		\[
		\|A_z-d_+(w)\|_{L^2}\le CT^{-a/2},\qquad
		\bigl\|\bigl(A_z-d_+(w)\bigr)_z\bigr\|_{L^2}\le CT^{-3a/2}.
		\]
		Combining this with Lemma~\ref{lem:outgoing-reference} below yields \eqref{eq:out-ordinary-budget-new}.
	\end{proof}
	
	\begin{lemma}\label{lem:outgoing-reference}
		The scalar reference solution $\Omega(t,z)$ of \eqref{eq:out-reference-new} is strictly increasing in $z$, remains between its far-field states $\Omega_\pm=h_0(s_\pm)$, and satisfies
		\[
		0\le\Omega_z\le C\min\{\delta e,T^{-1}\},\qquad
		0\le\partial_z\bigl(f_0'(\Omega)\bigr)\le\frac{1+o_\delta(1)}{T}.
		\]
		For every fixed $0<\epsilon<1/10$ and $\theta_+<1/4-\epsilon/2$,
		\[
		\sup_{t\ge 0}T^{2-\epsilon}\|\Omega_{zz}(t)\|_{L^2}^2
		+\int_0^\infty T^{2-\epsilon}\|\Omega_{zzz}(t)\|_{L^2}^2\dd t=o_\delta(1),\qquad
		\int_0^\infty T^{\theta_+}\|\Omega_{zz}(t)\|_{L^\infty}\dd t=o_\delta(1).
		\]
	\end{lemma}
	
	\begin{proof}
		The derivative $\pi_+=\Omega_z\ge 0$ satisfies
		\[
		\pi_{+,t}+\bigl((f_0'(\Omega)+\Sigma)\pi_+\bigr)_z=\bigl(d_+(\Omega)\pi_+\bigr)_{zz}.
		\]
		Setting $\mu_\Omega=f_0'(\Omega)$, the chain rule gives
		\[
		(\mu_\Omega)_t+(\mu_\Omega+\Sigma)(\mu_\Omega)_z
		=d_+(\Omega)(\mu_\Omega)_{zz}
		+\left(\frac{d_+'}{f_0''}-\frac{d_+f_0'''}{(f_0'')^2}\right)(\Omega)(\mu_\Omega)_z^2.
		\]
		By the parabolic maximum principle, $L(t)=\sup_z(\mu_\Omega)_z$ satisfies
		\[
		D^+L(t)\le-\bigl(1-CL(t)\bigr)L(t)^2+\|\Sigma_z^-(t)\|_{L^\infty}L(t),
		\]
		with $\int_0^\infty\|\Sigma_z^-\|_{L^\infty}\dd t=o_\delta(1)$ by \eqref{eq:out-drift-new}, so Bernoulli comparison yields $L(t)\le(1+o_\delta(1))/(1+t)$ and $\|\pi_+(t)\|_{L^2}^2\le C\delta/T$. Multiplying the $\pi_+$ equation by $\pi_+$ and integrating over $[t,\infty)\times\mathbb{R}$ gives
		\begin{equation}\label{eq:out-reference-tail}
			\int_t^\infty\|\pi_{+,z}(s)\|_{L^2}^2\dd s\le \frac{C\delta}{T}.
		\end{equation}
		Testing the $\pi_+$ equation against $-\pi_{+,zz}$ with $E_\pi=\|\pi_{+,z}\|_{L^2}^2$ and $D_\pi=\|\pi_{+,zz}\|_{L^2}^2$ yields
		\[
		E_\pi'+cD_\pi\le k_+E_\pi
		+C\|\pi_+\|_{L^\infty}^2\|\Sigma_z\|_{L^2}^2+C\int_{\mathbb{R}}\pi_+^5\dd z,\qquad
		k_+=C\bigl(\|\pi_+\|_{L^\infty}^2+\|\Sigma_z\|_{L^\infty}\bigr)\in L^1_t.
		\]
		Multiplying this inequality by $T^{2-\epsilon}$ and using $\int_0^\infty T^{1-\epsilon}E_\pi(t)\dd t\le C_\epsilon\delta$ from \eqref{eq:out-reference-tail} proves the weighted $L^2$ curvature bounds, and Gagliardo--Nirenberg interpolation gives the $L^1_t(T^{\theta_+}L^\infty)$ bound.
	\end{proof}
	
	\begin{lemma}\label{lem:outgoing-nash}
		The scalar perturbation $q_+=w-\Omega$ in \eqref{eq:out-error-ordinary} satisfies
		\begin{equation}\label{eq:out-small-log-new}
			\int_0^t\|q_+(s)\|_{L^\infty}^2\dd s
			\le C_\delta+\varepsilon_\delta\log(2+t),\qquad
			\lim_{\delta\rightarrow 0}\varepsilon_\delta=0.
		\end{equation}
	\end{lemma}
	
	\begin{proof}
		Since $B(\Omega,0)=0$ and $q_+(0)=0$, Kato's inequality applied to \eqref{eq:out-error-ordinary} and Lemma~\ref{lem:outgoing-source} give
		\[
		M_+=\sup_{t\ge 0}\|q_+(t)\|_{L^1}\le e^{\int_0^\infty\|\Sigma_z\|_{L^\infty}\dd t}\int_0^\infty n_+(t)\dd t=o_\delta(1).
		\]
		Multiplying \eqref{eq:out-error-ordinary} by $q_+$ with $E_+=\|q_+\|_{L^2}^2$ and $D_+=\|q_{+,z}\|_{L^2}^2$, and using the convexity
		\[
		\Omega_z\int_0^{q_+}\bigl(f_0'(\Omega+\xi)-f_0'(\Omega)\bigr)\dd\xi\ge 0,
		\]
		we obtain
		\[
		E_+'+cD_+\le kE_++Cn_+\|q_+\|_{L^\infty},\qquad
		k=C\bigl(\|\Sigma_z\|_{L^\infty}+\|\Omega_z\|_{L^\infty}^2\bigr),\qquad
		\int_0^\infty k(t)\dd t=o_\delta(1).
		\]
		By the one-dimensional Nash and Gagliardo--Nirenberg inequalities $E_+\le CM_+^{4/3}D_+^{1/3}$ and $\|q_+\|_{L^\infty}\le CM_+^{1/3}D_+^{1/3}$, Young's inequality gives
		\[
		\frac{\dd}{\dd t}\left(T^{1/2}e^{-\int_0^t k(s)\dd s}E_+\right)+cT^{1/2}D_+
		\le \frac{CM_+^2}{T}+CM_+^{1/2}T^{1/2}n_+^{3/2}.
		\]
		Since $T^{1/2}n_+^{3/2}\le CT^{1/2-a}n_+=CT^{1/50}n_+\in L^1_t$ by \eqref{eq:out-ordinary-budget-new}, integrating over $[0,t]$ yields
		\[
		\int_0^t(1+s)^{1/2}D_+(s)\dd s\le C_\delta+CM_+^2\log(1+t).
		\]
		Using $\|q_+\|_{L^\infty}^2\le CM_+^{2/3}D_+^{2/3}\le CT^{1/2}D_++CM_+^2T^{-1}$ then proves \eqref{eq:out-small-log-new} with $\varepsilon_\delta=CM_+^2\to 0$.
	\end{proof}
	
	\begin{proof}[Proof of Proposition~\ref{prop:outgoing-amplitude}]
		Set $\bar{s}(t,y)=h_0^{-1}(\Omega(t,y))$. Since $h_0(s(t,y))=w(t,Z_+(t,y))+C_+(t,y)$, combining \eqref{r5x:M-bound} with Lemma~\ref{lem:outgoing-reference} yields
		\[
		\bigl|\Omega\bigl(t,Z_+(t,y)\bigr)-\Omega(t,y)\bigr|\le C_\delta T^{-1}\log(2+t),\qquad
		\|C_+(t)\|_{L^\infty}\le CT^{-2a},
		\]
		whose squares are time-integrable on $[0,\infty)$. Combining these with Lemma~\ref{lem:outgoing-nash} establishes the logarithmic square bound for $s-\bar{s}$. {Since $\lambda_+^E(s)=c_-+(\gamma+1)s/4$, \eqref{v3:scalar-derivatives} gives $\partial_1\lambda_+^E(\bar s(t,Y))=\rho_b v_b(\bar s)\,[\partial_z f_0'(\Omega)](t,Y)$, where $\rho_bv_b(\bar s)=1+o_\delta(1)$.} On $z\le\nu_0 t-L$ ($0<\nu_0<c_*/2$), $K(1+\nu_0t-z)^{-3}$ and $K(1+\nu_0t-z)^{-4}{\exp(\int_0^t\|\Sigma_z^-\|_{L^\infty}\dd s)}$ are supersolutions for $\Omega-\Omega_-$ and $\Omega_z${ for large $L$: the expanded $\pi_+$ equation has drift $f_0'+\Sigma-3d_+'\pi_+>\nu_0$ uniformly and zeroth-order coefficient $\Sigma_z+f_0''\pi_+-d_+''\pi_+^2\ge-\|\Sigma_z^-\|_{L^\infty}$, since $\pi_+$ is uniformly small. The exponential is bounded by \eqref{eq:out-drift-new}}, which gives the $O(T^{-3})$ and $O(T^{-4})$ bounds for $\bar{s}-\bar{s}_-$ and $\bar{s}_y$ on $\{|y|\le C_0T^b\}$.
	\end{proof}
	
	\section{Three-dimensional weighted energy estimates}
	\label{sec:three-dimensional-all}
	
	\subsection{The potential energy and coercivity}
	\label{sec:r5a-full-interface}
	\label{r5a:interface-theorem}
	
	Having completed the one-dimensional characteristic analysis in Sections~\ref{sec:elementaryplanarinputs}--\ref{sec:outgoing-energy-replacement}, we now turn to the three-dimensional stability problem on $\mathbb{R}^3$. Let
	\[
	U_b(t,x_1)=\bigl(\rho_b(t,x_1),u_b(t,x_1)e_1,\theta_b(t,x_1)\bigr)
	\]
	be the global one-dimensional viscous rarefaction wave constructed from the smooth approximate rarefaction profile. We first define a weighted $\dot{H}^{-1}(\mathbb{R}^3)$ energy that is uniformly coercive in the transverse frequency $\zeta\in\mathbb{R}^2$, using the weights built from $\Psi_E$, $S-S_-=H_P+\sigma_0$, and $\bar{s}$ in Sections~\ref{sec:viscous-entropy}--\ref{sec:outgoing-energy-replacement}.
	
	Throughout Sections~\ref{sec:r5a-full-interface}--\ref{sec:main-proof}, we set
	\[
	\tilde{\nu}=\nu+\frac{(\gamma-1)\kappa}{c_p},\qquad p_b=R_g\rho_b\theta_b,\qquad T=1+t,\qquad J=\log(2+t).
	\]
	Functions of the mass coordinate $y$ appearing in Eulerian formulas are evaluated at $y=Y(t,x_1)$. Constants $C>0$ and fixed exponents $M$ are independent of $t$ and the transverse frequency $\zeta\in\mathbb{R}^2$, and spatial integrals are taken over $\mathbb{R}^3$ unless otherwise indicated.
	
	\subsubsection{Planar coefficient estimates}
	\label{ssec:3d-planar-inputs}
	
	Fix the planar decay exponent $a=12/25$, and write $\varepsilon_b,\varepsilon_\delta\to 0$ as $\delta\to 0$ after all fixed physical multipliers have been chosen. The Lagrangian mass coordinate $Y(t,x_1)$ satisfies
	\[
	Y_{x_1}=\rho_b,\qquad Y_t=-\rho_b u_b,\qquad D_{u_b}Y=0,\qquad D_{u_b}=\partial_t+u_b\partial_1.
	\]
	Set
	\[
	K_b(t)=\sum_{j=1}^4\bigl\|\partial_1^j(\rho_b,u_b,\theta_b)(t)\bigr\|_{L^\infty(\mathbb{R})}^2.
	\]
	The following proposition collects the one-dimensional estimates from Sections~\ref{sec:elementaryplanarinputs}--\ref{sec:outgoing-energy-replacement} in the Eulerian form used in the three-dimensional energy analysis.
	
	\begin{proposition}
		\label{lem:3d-planar-package}
		The planar wave $U_b(t,x_1)$ and the quantities constructed in Sections~\ref{sec:elementaryplanarinputs}--\ref{sec:outgoing-energy-replacement} satisfy the following bounds.
		\begin{enumerate}
			\item[(1)] The global background derivatives satisfy, for all $t\ge0$,
			\begin{equation}
				\sup_{t\ge 0} K_b(t)\le \varepsilon_b^2,\qquad
				K_b(t)\le C_\delta T^{-48/25},\qquad
				\int_t^\infty K_b(s)\dd s\le C_\delta T^{-23/25},
				\label{r5a:raw-gradient}
			\end{equation}
			so that $J^M K_b(t)\in L^1(0,\infty)$ for every fixed $M<\infty$. On $\{|Y|\le C T^b\}$ with $b=23/40$,
			\begin{equation}
				\sup_{|Y|\le CT^b}\bigl|(\rho_b,u_b,\theta_b)-(\rho_-,0,\theta_-)\bigr|
				\le \varepsilon_b T^{-1/6},\qquad
				\int_0^\infty J^M A_b(t)^2\dd t<\infty,
				\label{r5a:raw-center}
			\end{equation}
			where
			\[
			A_b(t)=\|u_b(t)\|_{L^\infty(\{|Y|\le CT^b\})}+\|p_b(t)-p_-\|_{L^\infty(\{|Y|\le CT^b\})}.
			\]
			
			\item[(2)] The relative entropy decomposes as $S-S_-=H_P+\sigma_0$ with $H_P\ge 0$ and $\sigma_0=H_N-\kappa_p(u_b)_y/\theta_b$. For all $t\ge 0$,
			\begin{equation}
				\int_0^t\|\sigma_0(s)\|_{L^\infty}^2\dd s\le C_\delta+\varepsilon_\delta\log T,
				\qquad
				\int_0^t\|(S-S_-)(s)\|_{L^\infty}^2\dd s\le C_\delta J^3.
				\label{r5a:raw-contact}
			\end{equation}
			With respect to the vorticity weight $w$ defined in \eqref{r5a:common-weights} below,
			\begin{equation}
				\bigl\|\sqrt{w}\,\partial_1 H_P(t)\bigr\|_{L^\infty}
				\le \frac{C\delta}{T}+k_C(t),\qquad k_C\in L^1(0,\infty).
				\label{r5a:contact-derivative-input}
			\end{equation}
			
			\item[(3)] Let $\alpha=a_b(s)\tilde{r}$, and let $\Psi_E(t,x_1)=\widetilde{\Psi}(t,Y(t,x_1))$ be the phase from Section~\ref{sec:actual-incoming-phase} with parameters
			\[
			\ell=T^{2/5},\qquad a_1=\eta_0 T^{-1/4},\qquad
			\beta_\Psi=\frac{27}{50},\qquad \beta_N=\frac{11}{20}.
			\]
			Set $\alpha_0=\chi_-(Y/(4T^{\beta_N}))(G_\ell*\alpha)$, $\alpha_E=\alpha_0/v$ with $v=\rho_b^{-1}$, and $\psi_E=\partial_1\Psi_E\ge 0$. Then $\Psi_E$ satisfies the Eulerian transport--diffusion equation
			\begin{equation}
				\bigl(\partial_t+(u_b-c_E)\partial_1\bigr)\Psi_E
				=D_E\partial_1\psi_E+\partial_1 F_\Psi+R_\Psi,
				\qquad c_E=\sqrt{\gamma R_g\theta_b},\qquad D_E=d_-v^2>0,
				\label{phase:eulerian-equation}
			\end{equation}
			together with the integrability bounds
			\begin{equation}
				\int_0^\infty\|R_\Psi(t)\|_{L^\infty}\dd t<\infty,\qquad
				\int_0^\infty J^M
				\bigl(\|\psi_E\|_{L^\infty}^2+\|F_\Psi\|_{L^\infty}^2+I_E^2\bigr)\dd t<\infty,
				\label{phase:eulerian-budgets}
			\end{equation}
			where
			\[
			I_E(t)=\int_{\operatorname{supp}\alpha_E}
			\frac{|\partial_1\alpha_E(t,x_1)|^2}{\psi_E(t,x_1)}\dd x_1,
			\]
			and
			\begin{equation}
				\sup_{t,x_1}\bigl(|\Psi_E|+|\psi_E|+|\partial_1\psi_E|
				+|\partial_1 F_\Psi|\bigr)\le C\varepsilon_b,
				\qquad \int_0^\infty|\psi_c'(t)|\dd t<\infty.
				\label{r5a:phase-small}
			\end{equation}
			In addition,
			\begin{equation}
				\kappa_*\alpha_y\le\frac{1+\varepsilon_s}{T}+k_s(t)\quad (k_s\in L^1_t),\qquad
				\alpha_r=\alpha-\alpha_0,\qquad
				\int_0^\infty J^M\|\alpha_r(t)\|_{L^\infty}^2\dd t<\infty.
				\label{phase:slope}
			\end{equation}
			
			\item[(4)] The function $\bar{s}$ from Section~\ref{sec:outgoing-energy-replacement} is strictly increasing in $x_1$, and $s-\bar{s}$ satisfies
			\begin{equation}
				\int_0^t\bigl\|(s-\bar{s})(\tau)\bigr\|_{L^\infty}^2\dd\tau
				\le C_\delta+\varepsilon_\delta\log T,\qquad \lim_{\delta\rightarrow 0}\varepsilon_\delta=0.
				\label{r5a:outgoing-regular}
			\end{equation}
		\end{enumerate}
	\end{proposition}
	
	\begin{proof}
		From \eqref{af:bootstrap} and the Burgers profile bounds, $\|\partial_1^j U_b(t)\|_{L^\infty}\le C_\delta T^{-2a}$ for $1\le j\le4$, with $a=12/25$. Hence $K_b(t)\le C_\delta T^{-4a}=C_\delta T^{-48/25}$, whose integral tail has exponent $4a-1=23/25$. On $\{|Y|\le CT^b\}$, setting $\eta_\delta=\sup_t\|U_b-V(W)\|_{L^\infty}\to 0$ and combining the decay of $U_b-V(W)$ with the barrier for $W-W_-$ gives
		\[
		|U_b-U_-|\le \min\bigl(\eta_\delta,CT^{-12/25}\bigr)+C\delta T^{-3}
		\le C\bigl(\eta_\delta^{47/72}+\delta\bigr)T^{-1/6}.
		\]
		Since $(u_b,p_b-p_-)$ are smooth functions of $(r,s)$ vanishing at the left state, $A_c$ from Proposition~\ref{prop:central-acoustic} bounds $A_b\le C(A_c+\delta T^{-3})$, which gives \eqref{r5a:raw-center}.
		
		The bounds \eqref{r5a:raw-contact}--\eqref{r5a:contact-derivative-input} follow from \eqref{ev:HNpoint}--\eqref{ev:Xlog}, \eqref{ev:HNy}--\eqref{ev:HPy}, and \eqref{eq:P-exterior}, since $|H_{P,y}|\le C\delta J/T{+|k_P(t)|}$, $H_{P,y}$ is integrable outside $\{|Y|\le T^{b_c}\}$ ($b_c>1/2$), and $\sqrt{w}\le C/J$ inside $\{|Y|\le T^{b_c}\}$. {Under $y=Y(t,x_1)$ one has $\partial_y=v\partial_1$, $\widetilde\Psi_y=v\psi_E$ and $\partial_1(d_-v)=0$. In \eqref{phase:eulerian-equation}--\eqref{r5a:phase-small}, $F_\Psi,R_\Psi$ denote respectively $v\widetilde F_\Psi$ and $\widetilde R_\Psi-(\partial_1v)\widetilde F_\Psi$, with the flattened coefficients evaluated at $Y$. The extra product is integrable by Cauchy--Schwarz, and $\partial_1(v\widetilde F_\Psi)=(\partial_1v)\widetilde F_\Psi+v\partial_1\widetilde F_\Psi$. Hence \eqref{eq:actual-phase-equation}--\eqref{eq:actual-phase-budgets} give} \eqref{phase:eulerian-equation}--\eqref{r5a:phase-small} by Lemma~\ref{lem:phase-flattening}. Finally, \eqref{inc:regular-slope} and Proposition~\ref{prop:outgoing-amplitude} give \eqref{phase:slope} and \eqref{r5a:outgoing-regular}.
	\end{proof}
	
	\subsubsection{The spatial weights}
	
	We fix the exponents
	\[
	b_c=\frac{13}{25},\qquad b_1=\frac{21}{40},\qquad
	b_0=\frac{53}{100},\qquad b=\frac{23}{40},
	\]
	which satisfy $1/2<b_c<b_1<b_0<\beta_\Psi<\beta_N<b$. Let $x_b(t,\cdot)=Y(t,\cdot)^{-1}$, let $v_{\min}>0$ be a lower bound for $(x_b)_y=v$, and fix $0<c_D\le v_{\min}/4$. With $x_c(t)=x_b(t,0)$ as in \eqref{intro:contact-trajectory}, we set
	\[
	\ell_0=c_D T^{b_0},\qquad \ell_1=\frac{c_D}{2}T^{b_1},\qquad J_c(t)=\bigl\{x_1\in\mathbb{R}:\ |x_1-x_c(t)|\le\ell_1\bigr\}.
	\]
	
	Fix $\epsilon_0=1/100$. Choose $K_0>0$ large enough so that the frozen exterior diffusion and vorticity quadratic forms are strictly positive whenever $K\ge K_0/2$. Fix the phase coupling parameter $\lambda_h>0$, the pressure exponent $m>0$, and a small exponent $p_c>0$. For sufficiently large times, we define the spatial weight profiles by
	\begin{equation}
		\left\{
		\begin{aligned}
			&w(t,x_1)=(1+K_0)\Bigl\{1+(1+\log T)\bigl(1+Y^2/T^{2b_c}\bigr)^{-p_c/2}\Bigr\}^{-2},
			\qquad K=w-1,\\
			&\tilde{m}(t,x_1)=\gamma+(m-\gamma)\chi_0\bigl(Y/T^{b_0}\bigr),\\
			&\mathcal{W}(t,x_1)=\frac{p_-}{R_g}e^{\lambda_h\Psi_E}
			\Bigl(\frac{p_b}{p_-}\Bigr)^{1-\tilde{m}/\gamma},
			\qquad \chi(t,x_1)=\frac{\rho_b\theta_b}{\mathcal{W}},
		\end{aligned}
		\right.
		\label{r5a:common-weights}
	\end{equation}
	where $\chi_0\in C^\infty(\mathbb{R})$ vanishes on $[-1,1]$ and equals $1$ outside $[-2,2]$. By construction, $\mathcal{W}(t,x_1)$ is independent of $x_1$ on $\{|Y|\le T^{b_0}\}$.
	
	\begin{lemma}
		\label{lem:3d-weights}
		The weights $(\mathcal{W},\chi,w)$ can be switched on smoothly {over the finite initial interval $1\le T\le T_0$} so that for all $t\ge 0$,
		\begin{equation}
			C_\delta^{-1}\le \mathcal{W}(t,x_1),\,\chi(t,x_1)\le C_\delta.
			\label{r5a:weight-power}
		\end{equation}
		For $T\ge T_0$, the vorticity weight $w$ satisfies
		\begin{equation}
			w_{\min}\approx J^{-2},\quad w\le 1+K_0,\quad
			-D_{u_b}\log w\le\frac{2b_c p_c}{T}+\frac{C}{TJ},\quad
			|\partial_1\log w|^2\le C_{p_c}T^{-2b_c}.
			\label{r5a:weight-derivatives}
		\end{equation}
		In addition, the central interval on which the vorticity shift $K=w-1$ is negative (even after retaining a fixed fraction of $w$) is contained in $\{|x_1-x_c(t)|\le L(t)\}$ with width $L(t)\le C T^{b_c}J^{1/p_c}$.
	\end{lemma}
	
	\begin{proof}
		Using $D_{u_b}Y=0$ and $\partial_1 Y=\rho_b$, direct differentiation of \eqref{r5a:common-weights} yields \eqref{r5a:weight-derivatives}. Solving the inequality $(1-c)w<1$ in the mass coordinate $Y$ and using the bounded Jacobian $(x_b)_y=v$ gives the support width bound $L(t)\le C T^{b_c}J^{1/p_c}$. Since
		\[
		R_{\mathcal{W}}=\max\left\{1,\ \frac{\mathcal{W}(t,x_c(t))}{\inf_{x_1}\mathcal{W}(t,x_1)}\right\}
		\]
		is uniformly bounded on the compact positive state region, we may fix $T_0<\infty$ large enough so that
		\[
		R_{\mathcal{W}}\frac{L(t)}{\ell_1(t)}+R_{\mathcal{W}}\frac{\ell_1(t)}{\ell_0(t)}\ll w_{\min}(t)\qquad {(\forall\,T\ge T_0/2)}.
		\]
		{For $T\le T_0/2$, take $K=K_0$. Over $T_0/2\le T\le T_0$, first switch $\mathcal W$ to \eqref{r5a:common-weights} with $K=K_0$, then smoothly interpolate $K$ from $K_0$ to the value in \eqref{r5a:common-weights}. In the latter interval $K$ dominates its target value, so the same localized lower bound applies; the switching derivatives are compactly supported in time.} This establishes \eqref{r5a:weight-power}--\eqref{r5a:weight-derivatives}.
	\end{proof}
	
	\subsubsection{Coercivity of the potential energy}
	
	Set $q=c_p\theta_b$, $\bar{\gamma}=\gamma-1$, and $c_E^2=\bar{\gamma} q$. For $h=(h_\rho,h_m,h_E)$, we set
	\[
	\left\{
	\begin{aligned}
		&b_j=\partial_j h_\rho+u_b\partial_j h_{m1}+\tfrac{1}{2}u_b^2\partial_j h_E,\qquad
		v_{ij}=\partial_j h_{mi}+u_b\delta_{i1}\partial_j h_E,\\
		&a=b+q\nabla h_E,\qquad d=\operatorname{tr}v,\qquad
		\Omega=\operatorname{curl}(h_m+u_b h_E e_1).
	\end{aligned}
	\right.
	\]
	Since $u_b=u_b(t,x_1)$ depends only on $x_1$, we have $h_E\nabla u_b\times e_1=0$, so every quantity above depends purely on the first-order spatial gradient $\nabla h$. We define $\mathcal{J}$ and the quadratic form $\mathcal{A}_t[h,h]$ by
	\begin{equation}
		\mathcal{J}=\frac{|a|^2}{c_E^2}+\frac{|b|^2}{q}+|v|^2+K|\Omega|^2,
		\qquad \mathcal{A}_t[h,h]=\int_{\mathbb{R}^3} \mathcal{W}\mathcal{J}\dd x.
		\label{r5a:jet-form}
	\end{equation}
	The second-order elliptic operator $\mathcal{A}_t$ associated with $\mathcal{A}_t[h,h]$ is given by
	\[
	\mathcal{A}_t h=-\partial_j\bigl(\chi^{-1} G_{jk}\partial_k h\bigr),\qquad
	G_{jk}=\delta_{jk}M+K\rho_b\theta_b C_j(u_b)^\top C_k(u_b),
	\]
	where $M=(D_{\mathbf{U}}^2(-\rho S)(\mathbf{U}_b))^{-1}$ and $C_j(u_b)z=e_j\times(z_m+u_b z_E e_1)$.
	
	\begin{lemma}
		\label{lem:3d-coercivity}
		With the spatial weights constructed in Lemma~\ref{lem:3d-weights}, there exists a fixed constant $c>0$ independent of time $t$ and transverse Fourier frequency $\zeta\in\mathbb{R}^2$ such that
		\begin{align}
			\mathcal{A}_t[h,h]\ge c\Bigg[&
			\int_{\mathbb{R}^3} \mathcal{W}\left(\frac{|a|^2}{c_E^2}+\frac{|b|^2}{q}+d^2+w|\Omega|^2\right)\dd x\notag\\
			&+\int_{J_c(t)^c}\mathcal{W}|\nabla h_m|^2\dd x
			+w_{\min}(t)\int_{\mathbb{R}^3} \mathcal{W}|\nabla h_m|^2\dd x\Bigg].
			\label{r5a:strong-first}
		\end{align}
		In particular, $(a,b)$ and $\operatorname{div}h_m$ have uniform $O(1)$ weight in $\mathcal{A}_t[h,h]$, while $\nabla h_m$ has $O(1)$ weight outside $J_c(t)$ and $O(J^{-2})$ weight on $\mathbb{R}^3$.
	\end{lemma}
	
	\begin{proof}
		Although $K=w-1$ takes the negative value $-1+w_{\min}$ inside $I=\{|x_1-x_c(t)|\le L(t)\}$ with $L(t)\lesssim T^{b_c}J^{1/p_c}$, the weight $\mathcal{W}=\mathcal{W}_c$ is constant on the wider interval $J_0=\{|x_1-x_c(t)|\le \ell_0(t)\}$ with $\ell_0(t)=c_D T^{b_0}\gg L(t)$. Where the weight is constant, the divergence--curl identity gives
		\[
		\int_{\mathbb{R}^3}|{\nabla  h_m}|^2\dd x=\int_{\mathbb{R}^3}\bigl(|{\operatorname{div} h_m}|^2+|{\operatorname{curl} h_m}|^2\bigr)\dd x,
		\]
		so the negative $-|{\operatorname{curl} h_m}|^2$ part of $K|{\operatorname{curl} h_m}|^2$ on $I$ is compensated by $|{\nabla  h_m}|^2$ up to the off-diagonal commutator of the Riesz operator $T_0q=\nabla\operatorname{curl}(-\Delta)^{-1}q$ between $I$ and $J_0^c$.
		
		At transverse Fourier frequency $\zeta\in\mathbb{R}^2$ with $\rho_0=|\zeta|$, the one-dimensional resolvent kernel of $(-\partial_1^2+\rho_0^2)^{-1}$ is $e^{-\rho_0|x_1-y_1|}/(2\rho_0)$. Between the separated supports $I$ and $J_0^c$, the local Dirac delta term vanishes and the off-diagonal kernel is bounded by $C\rho_0 e^{-\rho_0|x_1-y_1|}$, whose squared Hilbert--Schmidt norm on $\mathbb{R}\setminus J_0$ satisfies
		\begin{equation}
			\|\mathbf{1}_{J_0^c}T_0\mathbf{1}_I\|_{L^2\to L^2}^2\le C\rho_0 L e^{-2\rho_0(\ell_0-L)}\le \frac{C L}{\ell_0},
			\label{r5a:off-diagonal-proof}
		\end{equation}
		uniformly in $\rho_0\ge 0$ (and vanishing at $\rho_0=0$). By duality and the Cauchy--Schwarz inequality, for any $q$ supported in $I$ and $G={\nabla  h_m}$,
		\[
		\begin{aligned}\mathcal{W}_c\bigl|\langle q,{\operatorname{curl} h_m}\rangle\bigr|^2
			&\le\left(\int_{\mathbb{R}^3} \mathcal{W}|G|^2\dd x\right)
			\left(\mathcal{W}_c\int_{\mathbb{R}^3} \mathcal{W}^{-1}|T_0q|^2\dd x\right)
			\\ &\le\Bigl(1+C R_{\mathcal{W}}\frac{L}{\ell_0}\Bigr)\|q\|_{L^2}^2\int_{\mathbb{R}^3} \mathcal{W}|G|^2\dd x.
		\end{aligned}\]
		Hence, whenever $K\ge 0$ outside $I$ and $K\ge -1+w_0$ on $I$ ($0<w_0\le 1$),
		\begin{equation}
			\int_{\mathbb{R}^3} \mathcal{W}\bigl(|G|^2+K|{\operatorname{curl} h_m}|^2\bigr)\dd x
			\ge\Bigl(w_0-C R_{\mathcal{W}}\frac{L}{\ell_0}\Bigr)\int_{\mathbb{R}^3} \mathcal{W}|G|^2\dd x.
			\label{r5a:localized-vector-proof}
		\end{equation}
		
		To retain $O(1)$ control of both ${\operatorname{div} h_m}$ on $J_1=\{|x_1-x_c(t)|\le\ell_1\}$ and ${\nabla  h_m}$ on $J_1^c$ ($L\ll\ell_1\ll\ell_0$), we apply the same off-diagonal estimate to
		\[
		T_0(q_c,q_d)=\nabla\operatorname{curl}(-\Delta)^{-1}q_c-\nabla\nabla(-\Delta)^{-1}q_d
		\]
		with $\operatorname{supp}q_c\subset I$ and $\operatorname{supp}q_d\subset J_1$, which gives
		\begin{equation}
			\mathcal{W}_c\left[\int_I|{\operatorname{curl} h_m}|^2\dd x+
			\int_{J_1}|{\operatorname{div} h_m}|^2\dd x\right]
			\le\Bigl(1+C R_{\mathcal{W}}\frac{\ell_1}{\ell_0}\Bigr)\int_{\mathbb{R}^3} \mathcal{W}|G|^2\dd x.
			\label{r5a:joint-duality-proof}
		\end{equation}
		{First retain a fixed small fraction of $w|\operatorname{curl}h_m|^2$, enlarging $I$ as in Lemma~\ref{lem:3d-weights}. Repeat \eqref{r5a:localized-vector-proof} with $\mathcal W(1-\tfrac14\mathbf1_{J_1^c})$ in the dual estimate. This weight is unchanged on $I$, and the added error $C R_{\mathcal W} L/\ell_1$ is $o(w_{\min})$ by \eqref{r5a:off-diagonal-proof}. Thus a quarter of the exterior gradient energy is retained. Averaging with \eqref{r5a:joint-duality-proof} retains a fixed fraction of the interior divergence energy as well as $cw_{\min}\int \mathcal W |G|^2$. Since $q\nabla h_E=a-b$, the differences $v-G$ and $\Omega-\operatorname{curl}h_m$ have size at most $C|u_b|(|a|+|b|)$. Young's inequality bounds the resulting errors by $\epsilon(|a|^2+|b|^2)+C_{\epsilon}u_b^2|G|^2$, which are absorbed using $\sup_{J_1}u_b^2\ll w_{\min}$ and the smallness of $u_b$ on $J_1^c$. The same comparison for $d=\operatorname{div}h_m+u_b\partial_1h_E$ gives \eqref{r5a:strong-first}. For a region $F$ at distance $H\ge2\ell_0$ from $I$, use the weight $\mathcal W(\mathbf1_{F^c}+\epsilon_H\mathbf1_F)$ in the curl estimate, retaining the same exterior gradient fraction. The additional error is at most $C R_{\mathcal W} L/(H\epsilon_H)$ by \eqref{r5a:off-diagonal-proof}. Taking $\epsilon_H$ a sufficiently large fixed multiple of $R_{\mathcal W}L/(Hw_{\min})$ gives $E_{F^c}\ge-\epsilon_H E_F$ and hence $E_F\le (1-\epsilon_H)^{-1}E_{-1}\le(1+\epsilon_0)E_{-1}$ for large $T_0$.}
	\end{proof}
	
	\subsubsection{The characteristic Euler quadratic forms}
	
	We now decompose the Euler transport quadratic form into the three characteristic parts and $\mathcal{F}_\pm^2$. Put $\bar{\gamma}=\gamma-1$ and set
	\[
	\hat{a}_j=\frac{a_j}{c_E},\qquad \hat{b}_j=\frac{b_j}{\sqrt{q}},
	\qquad
	\mathcal{J}=\sum_{j=1}^3\bigl(\hat{a}_j^2+\hat{b}_j^2+|v_{\cdot j}|^2\bigr)+K|\Omega|^2.
	\]
	In the exterior region where $K>0$, we define $s_\mp$ and $\mathcal{F}_\mp^2$ by
	\begin{equation}\label{r5a:acoustic-projections}
		\left\{
		\begin{aligned}
			&s_-=\frac{\hat{a}_1-v_{11}}{\sqrt{2}},\qquad
			s_+=\frac{\hat{a}_1+v_{11}}{\sqrt{2}},\\
			&\mathcal{F}_\mp^2=\sum_{j=1}^3(\hat{a}_j\pm v_{1j})^2
			+\sum_{i=2}^3\sum_{j=1}^3v_{ij}^2
			+\sum_{j=1}^3\hat{b}_j^2+K|\Omega|^2.
		\end{aligned}
		\right.
	\end{equation}
	The kernel of $\mathcal{F}_-^2$ is the backward acoustic line $\hat{a}_1=-v_{11}$ (with all transverse, entropy, and vorticity components vanishing), while the kernel of $\mathcal{F}_+^2$ is the forward acoustic line $\hat{a}_1=v_{11}$.
	
	In terms of these quantities, the Euler transport quadratic forms are
	\begin{align*}
		&\mathcal{P}_\gamma=\bar{\gamma}\hat{a}_1^2
		-2\sqrt{\bar{\gamma}}\,\hat{a}_1\hat{b}_1+\hat{b}_1^2
		+\frac{2-\bar{\gamma}}{2}v_{11}^2
		-\frac{2-\bar{\gamma}}{2}(v_{12}^2+v_{13}^2)\\
		&\qquad+\left(1+\frac{\bar{\gamma}}{2}\right)(v_{21}^2+v_{31}^2)
		+\frac{\bar{\gamma}}{2}\sum_{i,j=2}^3 v_{ij}^2
		-\bar{\gamma} v_{11}(v_{22}+v_{33})
		+v_{21}v_{12}+v_{31}v_{13},\\
		&\mathcal{C}=\sum_{j=2}^3(v_{j1}-v_{1j})
		\Bigl(\hat{a}_j-\frac{\hat{b}_j}{\sqrt{\bar{\gamma}}}\Bigr),\qquad
		\mathcal{V}=\sum_{i=2}^3 v_{i1}^2-\sum_{j=2}^3 v_{1j}^2,\\
		&\mathcal{Q}_\pm=\mathcal{P}_\gamma\pm \bar{\gamma}\mathcal{C}
		+K\Bigl(\frac{\bar{\gamma}}{2}|\Omega|^2+\mathcal{V}\Bigr).
	\end{align*}
	Decomposing $\partial_1 U_b$ into its acoustic and entropy parts,
	\[
	\left\{
	\begin{aligned}
		&\mu_\pm=\frac{1}{2}\left(\partial_1 u_b\pm \frac{\partial_1 p_b}{\rho_b c_E}\right),\qquad
		(\partial_1\rho_b,\partial_1 u_b,\partial_1\theta_b)_\pm=\Bigl(\pm\frac{\rho_b}{c_E},\,1,\,\pm\frac{\bar{\gamma}\theta_b}{c_E}\Bigr)\mu_\pm,\\
		&(\partial_1\rho_b,\partial_1 u_b,\partial_1\theta_b)_0=\Bigl(-\frac{\rho_b}{c_p},\,0,\,\frac{\theta_b}{c_p}\Bigr)S_1,
	\end{aligned}
	\right.
	\]
	and using $\mathrm{d}\log[(p_b/p_-)^{-m/\gamma}]=-m\,\mathrm{d}\rho_b/\rho_b$ along the isentropic acoustic directions (and zero along the isobaric entropy direction), {the part of the Euler drift proportional to $\mu_\pm,S_1$ in the region where $\tilde{m}=m$ takes} the characteristic form
	\begin{equation}
		\int_{\mathbb{R}^3} \mathcal{W}\left[
		\mu_+\Bigl(\mathcal{Q}_+-\frac{m}{2}\mathcal{F}_+^2\Bigr)
		+\mu_-\Bigl(\mathcal{Q}_--\frac{m}{2}\mathcal{F}_-^2\Bigr)
		+\frac{c_E}{c_p}S_1(1+K)\mathcal{C}\right]\dd x.
		\label{r5a:explicit-acoustic-polynomial}
	\end{equation}
	On either principal acoustic line ($\mathcal{F}_\pm=0$), the normalization in \eqref{r5a:acoustic-projections} gives $\mathcal{J}=s_\pm^2$ and $\mathcal{Q}_\pm=\beta s_\pm^2$, where $\beta=(\gamma+1)/4$.
	
	\begin{lemma}
		\label{lem:3d-acoustic-completion}
		Put $\beta=(\gamma+1)/4$. On the exterior region where $K\ge K_0/2>0$, for every fixed $\epsilon>0$, choosing a sufficiently large finite pressure exponent {$m=m(\epsilon,\gamma,K_0)>0$} ensures that
		\begin{equation}
			\mathcal{Q}_\pm-\frac{m}{2}\mathcal{F}_\pm^2
			\le(\beta+\epsilon)s_\pm^2-c_\epsilon\mathcal{F}_\pm^2.
			\label{r5a:normal-complement}
		\end{equation}
		In addition, the deviation of $\mathcal{Q}_\pm$ from its principal one-dimensional acoustic projection $\beta s_\pm^2$ is bounded pointwise by $C(|s_\pm|\mathcal{F}_\pm+\mathcal{F}_\pm^2)$.
	\end{lemma}
	
	\begin{proof}
		Decomposing the first-order potential gradients orthogonally with respect to $\mathcal{J}$ into the one-dimensional acoustic line spanned by $s_\pm$ and its orthogonal complement (on which $\mathcal{F}_\pm^2$ is uniformly coercive since $K\in[K_0/2,K_0]$), we see that $\mathcal{Q}_\pm-\beta s_\pm^2$ vanishes on the acoustic line and is therefore bounded by $C(|s_\pm|\mathcal{F}_\pm+\mathcal{F}_\pm^2)$. Applying Young's inequality $C|s_\pm|\mathcal{F}_\pm\le \epsilon s_\pm^2+C_\epsilon\mathcal{F}_\pm^2$ and choosing $m/2\ge C+C_\epsilon+c_\epsilon$ yields \eqref{r5a:normal-complement}.
	\end{proof}
	
	\subsection{The combined energy inequality}
	\label{sec:common-energy}
	
	In this subsection we establish Proposition~\ref{prop:3d-common-growth}, which combines $E_{-1}$ and $E_3$ so that the logarithmic growth exponent satisfies $A_*\le 203/200<2$, uniformly in $\zeta\in\mathbb{R}^2$.
	
	\subsubsection{The corrected perturbation equations}
	
	Set
	\[
	f=\mathbf{U}-\mathbf{U}_b=(f_\rho,f_m,f_E),\qquad
	Z=U-U_b=(Z_\rho,Z_u,Z_\theta).
	\]
	At the left far-field state $U_-=(\rho_-,0,\theta_-)$, the sound speed is $c_-=\sqrt{\gamma R_g\theta_-}$. Using the matrix
	\[
	T_*f=
	\left(\frac{\gamma-1}{\sqrt{\gamma R_g\theta_-}}f_E,\ f_m,\ 
	\frac{\sqrt{\gamma-1}}{\sqrt{\gamma R_g\theta_-}}f_E-
	\frac{\sqrt{\gamma R_g\theta_-}}{\sqrt{\gamma-1}}f_\rho\right),
	\]
	we define the skew-symmetric matrices $K_l^*$ in Fourier variables by
	\[
	A_*(\xi)=\sqrt{\gamma R_g\theta_-}\begin{pmatrix}0&\xi^\top&0\\\xi&0&0\\0&0&0\end{pmatrix},
	\qquad
	K^*(\xi)=\begin{pmatrix}
		0&k_a\xi^\top&0\\-k_a\xi&0&k_e\xi\\0&-k_e\xi^\top&0
	\end{pmatrix},
	\]
	with
	\[
	k_a=\frac{d_T(\gamma-1)/\gamma-\nu/\rho_-}{4\sqrt{\gamma R_g\theta_-}},\qquad
	k_e=-\frac{d_T\sqrt{\gamma-1}}{\gamma \sqrt{\gamma R_g\theta_-}},\qquad
	d_T=\frac{\kappa}{\rho_-c_v}.
	\]
	Writing $K_l^*$ in the conservative variables $f$ and setting
	\[
	\widetilde{f}=f+K_l^*\partial_l f,
	\]
	diagonalizes the principal second-order hyperbolic--parabolic coupling at the left state without introducing any $x_1$-dependent commutator in the definition of $\widetilde{f}$.
	
	\begin{lemma}
		\label{lem:3d-corrected-equation}
		The conservative perturbation $f=\mathbf{U}-\mathbf{U}_b$ satisfies
		\begin{equation}
			\partial_t f+\partial_j(A_j f)
			=\partial_j\bigl(B_{jk}\partial_k f+C_j f+N_j\bigr),
			\qquad C_j f=D_{\mathbf{U}}B_{j1}(\mathbf{U}_b)[f]\,\partial_1\mathbf{U}_b,
			\label{r5a:conservative}
		\end{equation}
		where $A_j=D_{\mathbf{U}}F_j(\mathbf{U}_b)$ and $B_{jk}=B_{jk}(\mathbf{U}_b)$ are the Euler and viscous flux Jacobians at the planar background, and the nonlinear flux remainder $N_j$ satisfies
		\begin{equation}
			|N_j|\le C\bigl(|f|^2+|f||\nabla f|+|\partial_1 U_b||f|^2\bigr).
			\label{r5a:nonlinear-flux}
		\end{equation}
		Setting $D_{jk}=B_{jk}+[A_j,K_k^*]$ and $T_{jkl}=K_l^*B_{jk}-D_{jk}K_l^*$, the variable $\widetilde{f}=f+K_l^*\partial_l f$ satisfies
		\begin{align}
			\partial_t\widetilde{f}+\partial_j(A_j\widetilde{f})=&\partial_j\bigl(D_{jk}\partial_k\widetilde{f}\bigr)+\partial_j\partial_l
			\bigl(T_{jkl}\partial_k f+K_l^*C_j f+K_l^*N_j\bigr)\notag\\
			&+\partial_j\bigl\{(\partial_l D_{jk})K_l^*\partial_k f
			+(C_j-K_l^*\partial_l A_j)f+N_j\bigr\}.
			\label{r5a:corrected-exact}
		\end{align}
	\end{lemma}
	
	\begin{proof}
		Subtracting the planar background equations from the three-dimensional system \eqref{r5main:NSF} and expanding the fluxes by Taylor's formula in conservative variables gives \eqref{r5a:conservative}--\eqref{r5a:nonlinear-flux}. Applying $I+K_l^*\partial_l$ to \eqref{r5a:conservative} and using the commutator identities
		\begin{align*}
			&(I+K_l^*\partial_l)\partial_j(A_j f)=\partial_j(A_j\widetilde{f})-\partial_j\bigl([A_j,K_l^*]\partial_l f{-}K_l^*(\partial_l A_j)f\bigr),\\
			&K_l^*\partial_l(B_{jk}\partial_k f)-D_{jk}K_l^*\partial_k\partial_l f=\partial_l(T_{jkl}\partial_k f)+(\partial_l D_{jk})K_l^*\partial_k f,
		\end{align*}
		together with $\partial_k f=\partial_k\widetilde{f}-K_l^*\partial_k\partial_l f$, gives \eqref{r5a:corrected-exact}.
	\end{proof}
	
	We now define $h=\mathcal{A}_t^{-1}\widetilde{f}\in\dot{H}^1(\mathbb{R}^3)$ and the weighted $\dot{H}^{-1}(\mathbb{R}^3)$ energy
	\[
	E_{-1}(t)=\frac{1}{2}\mathcal{A}_t[h,h]=\frac{1}{2}\langle\widetilde{f},h\rangle.
	\]
	By the resolvent differentiation identity for $\mathcal{A}_t^{-1}$, the time derivative of $E_{-1}(t)$ along local strong solutions is given by
	\begin{equation}
		E_{-1}'(t)=\langle\partial_t\widetilde{f},h\rangle-\frac{1}{2}\langle(\partial_t\mathcal{A}_t)h,h\rangle.
		\label{r5a:inverse-time}
	\end{equation}
	
	\subsubsection{Second-order potential dissipation}
	
	Let $\zeta_c\in C_c^\infty(\mathbb{R})$ be a smooth nonnegative cutoff equal to $1$ where $K<K_0/2$, supported in $\{|Y|\le T^{b_0}\}$, with derivative supported at scale $T^{b_1}$ where $K\ge K_0/2$. We define $\mathcal{H}$ by
	\begin{align}
		\mathcal{H}={}&\int_{\mathbb{R}^3}\zeta_c\mathcal{W}\Bigl(
		|\nabla^2 h_\rho|^2+|\nabla^2 h_E|^2+
		|\nabla\operatorname{div}h_m|^2+
		w|\nabla\operatorname{curl}h_m|^2\Bigr)\dd x+\int_{\mathbb{R}^3}(1-\zeta_c)\mathcal{W}|\nabla^2 h|^2\dd x.
		\label{r5a:Hessian-definition}
	\end{align}
	Let $\zeta_e$ be a smooth nonnegative cutoff vanishing on $\{|Y|\le 2T^{b_0}\}$ and equal to $1$ on $\{|Y|\ge 3T^{b_0}\}$. With $\psi_E=\partial_1\Psi_E\ge 0$, we define $\mathcal{T}_\rho$ and $D_-$ by
	\begin{equation}
		\mathcal{T}_\rho=\int_{\mathbb{R}^3}\zeta_e\mathcal{W}\psi_E\mathcal{F}_-^2\dd x,
		\qquad D_-=c_H\mathcal{H}+c_T\mathcal{T}_\rho,
		\label{r5a:retained-D}
	\end{equation}
	where $c_H,c_T>0$ are fixed constants.
	
	\begin{lemma}
		\label{lem:3d-diffusion}
		For sufficiently small planar wave strength, the principal viscous term in \eqref{r5a:corrected-exact} satisfies
		\[
		\bigl\langle\partial_j(D_{jk}\partial_k\widetilde{f}),h\bigr\rangle
		\le -c_H\mathcal{H}+C\mathcal{K}(t)E_{-1}(t),
		\]
		where the integrable coefficient $\mathcal{K}\in L^1(0,\infty)$ is defined in \eqref{r5a:K-integrable}, and the global second-order potential bound \eqref{r5a:Hessian-comparison} holds with at most a polynomial logarithmic factor $J^M$.
	\end{lemma}
	
	\begin{proof}
		At the left state $U_-$, the operator $\partial_j(D_{jk}^*\partial_k)$ is diagonal in the thermal and velocity components:
		\[
		L_*=\bigl(d_*\Delta,\ \mu_*\Delta+(d_*-\mu_*)\nabla\operatorname{div},\ d_e\Delta\bigr),
		\qquad
		d_*=\frac{\tilde{\nu}}{2\rho_-},\qquad
		\mu_*=\frac{\mu}{\rho_-},\qquad d_e=\frac{\kappa}{\rho_-c_p}.
		\]
		On $\{|Y|\le T^{b_0}\}$ where $\mathcal{W}$ is constant, the divergence--curl identity gives
		\begin{equation}
			-\mathcal{W}d_*\|{\sqrt{\zeta_c}}\nabla\operatorname{div}h_m\|_{L^2}^2
			-\mathcal{W}\mu_*\int_{\mathbb{R}^3} {\zeta_c} w|\nabla\Omega_*|^2\dd x
			-\mathcal{W}\mu_*\int_{\mathbb{R}^3} {\zeta_c}(\partial_1 w)\Omega_*\cdot\partial_1\Omega_*\dd x,
			\quad \Omega_*=\operatorname{curl}h_m,
			\label{r5a:central-diffusion}
		\end{equation}
		{up to terms with a derivative on $\zeta_c$. These are supported where $w\ge1+K_0/2$ and bounded by $\epsilon\mathcal H+C_\epsilon J^M T^{-2b_1}E_{-1}$.} By the Cauchy--Schwarz inequality and \eqref{r5a:weight-derivatives}, the $\partial_1 w$ commutator is bounded by
		\[
		\epsilon\mathcal{H}+C_\epsilon\|\partial_1\log w\|_{L^\infty}^2 E_{-1}.
		\]
		By \eqref{r5a:raw-center}, the difference between $D_{jk}(\mathbf{U}_b)$ and $D_{jk}^*$ on $\{|Y|\le T^{b_0}\}$ is $O(\varepsilon_b T^{-1/6})$, which is absorbed into $\mathcal{H}$ since $\varepsilon_b\sup_{T\ge 1}T^{-1/6}J^M/w_{\min}\ll 1$. All remaining commutator terms involving $\partial_1 U_b$, $\partial_1 w$, or cutoff derivatives are bounded via Young's inequality by
		\begin{equation}
			\epsilon\mathcal{H}+C_\epsilon\mathcal{K}(t)E_{-1}(t),
			\qquad \mathcal{K}\in L^1(0,\infty),
			\label{r5a:K-integrable}
		\end{equation}
		where the integrable remainder $\mathcal{K}(t)$ is given by
		\begin{align*}
			\mathcal{K}(t)={}&J^M\Bigl(K_b(t)+{T^{-3a}+}T^{-2b_c}+T^{-2b_1}+T^{-2b_0}+\|\psi_E(t)\|_{L^\infty}^2+\|F_\Psi(t)\|_{L^\infty}^2\Bigr)\\
			&+I_E(t)^2{+\|\alpha_r(t)\|_{L^\infty}^2}+\|R_\Psi(t)\|_{L^\infty}+{J^M}|\psi_c'(t)|+k_0(t).
		\end{align*}
		In the exterior region where $K\ge K_0/2>0$, integrating by parts in $(h_a,h_m,h_e)$ gives
		\begin{align}
			&-d_*\int_{\mathbb{R}^3} {(1-\zeta_c)\mathcal{W}}\bigl(|\nabla^2 h_a|^2+|\nabla^2 h_m|^2\bigr)\dd x
			-d_e\int_{\mathbb{R}^3} {(1-\zeta_c)\mathcal{W}}|\nabla^2 h_e|^2\dd x\notag\\
			&-{\int_{\mathbb{R}^3}(1-\zeta_c)\mathcal W[\mu_*K-(d_*-\mu_*)]}|\nabla\Omega_*|^2\dd x\notag\\
			&+\frac{d_*}{2}\int_{\mathbb{R}^3} {\partial_1^2[(1-\zeta_c)\mathcal W]}\bigl(|\nabla h_a|^2+|\nabla h_m|^2\bigr)\dd x
			+\frac{d_e}{2}\int_{\mathbb{R}^3} {\partial_1^2[(1-\zeta_c)\mathcal W]}|\nabla h_e|^2\dd x\notag\\
			&+{\frac{\mu_*}{2}\int_{\mathbb{R}^3}\partial_1^2[(1-\zeta_c)\mathcal W K]}|\Omega_*|^2\dd x+(d_*-\mu_*)\int_{\mathbb{R}^3} {\partial_1[(1-\zeta_c)\mathcal W]}\sum_{l=1}^3(e_1\times\partial_l h_m)\cdot\partial_l\Omega_*\dd x.
			\label{r5a:exterior-frozen}
		\end{align}
		{Integrating the second derivatives of the products in \eqref{r5a:exterior-frozen} by parts once} and using
		\[
		\partial_1\log\mathcal{W}=\lambda_h\psi_E+O\bigl(|\partial_1 U_b|+T^{-b_0}\bigr)
		\]
		{bounds the commutators by \eqref{r5a:K-integrable}, since $|\partial_1K|\le C T^{-b_c}$ and $|\partial_1\zeta_c|\le C T^{-b_1}$ on the transition layer, where $w\ge1+K_0/2$.} Combining the interior and exterior bounds gives
		\begin{equation}
			\int_{\mathbb{R}^3}\chi^{-1}|\nabla^2 h|^2\dd x
			\le C J^M\bigl(D_-+\mathcal{K}(t)E_{-1}(t)\bigr).
			\label{r5a:Hessian-comparison}
		\end{equation}
		We have finished the proof of Lemma \ref{lem:3d-diffusion}.
	\end{proof}
	
	\subsubsection{Background terms in the potential energy}
	
	We now show how the one-dimensional estimates in Proposition~\ref{lem:3d-planar-package} control the Euler drift terms in $E_{-1}'(t)$. Substituting $\widetilde{f}=\mathcal{A}_t h$ into
	\[
	-\bigl\langle\partial_j(A_j\widetilde{f}),h\bigr\rangle-\frac{1}{2}\bigl\langle(\partial_t\mathcal{A}_t)h,h\bigr\rangle
	\]
	and integrating by parts yields the sum $Q_0+Q_\Omega$, where
	\begin{align}
		Q_0&=-\frac{1}{2}\sum_{l=1}^3\int_{\mathbb{R}^3} \partial_l h^\top
		\bigl[(\chi^{-1} M)_t+(\chi^{-1} A_1 M)_1\bigr]\partial_l h\dd x
		+\sum_{j=1}^3\int_{\mathbb{R}^3}\chi^{-1} \partial_1 h^\top M(\partial_1 A_j^\top)\partial_j h\dd x,\notag\\
		Q_\Omega&=-\frac{1}{2}\int_{\mathbb{R}^3}\bigl[(K\mathcal{W})_t+(u_b K\mathcal{W})_1\bigr]|\Omega|^2\dd x
		+\int_{\mathbb{R}^3} K\mathcal{W}(\partial_1 u_b)\left(\sum_{i=2}^3 v_{i1}^2-\sum_{j=2}^3 v_{1j}^2\right)\dd x
		\notag\\
		&\quad+\int_{\mathbb{R}^3} K\mathcal{W}\left(\theta_b S_1+\frac{\nu\partial_1^2 u_b}{\rho_b}\right)
		\Omega\cdot(e_1\times\nabla h_E)\dd x.
		\label{r5a:Euler-exact}
	\end{align}
	
	We bound the four contributions in \eqref{r5a:Euler-exact} as follows.
	
	\medskip
	\noindent\textit{1. The entropy wave contribution.} In \eqref{r5a:explicit-acoustic-polynomial} and \eqref{r5a:Euler-exact}, the entropy gradient $S_1=\partial_1 H_P+\partial_1\sigma_0$ enters through $\int_{\mathbb{R}^3} \mathcal{W} w\theta_b S_1\Omega\cdot(e_1\times\nabla h_E)\dd x$. For the viscous production part $\partial_1 H_P$, the weighted bound \eqref{r5a:contact-derivative-input} and Lemma~\ref{lem:3d-coercivity} bound its contribution directly by $(C\delta/T+k_C(t))E_{-1}$. For the regular part $\partial_1\sigma_0$, integrating by parts in $x_1$ transfers $\partial_1$ onto $\Omega$ and $\nabla h_E$, which are absorbed into the second-order dissipation $D_-$ via Young's inequality:
	\begin{equation}
		\epsilon D_-+C_\epsilon\|\sigma_0(t)\|_{L^\infty}^2 E_{-1}(t)+C_\epsilon\mathcal{K}(t)E_{-1}(t),
		\label{r5a:contact-X-allocation}
	\end{equation}
	where $\int_0^t\|\sigma_0\|_{L^\infty}^2\dd s\le C_\delta+\varepsilon_\delta\log T$ by \eqref{r5a:raw-contact}.
	
	\medskip
	\noindent\textit{2. Central acoustic contribution.} On $\{|Y|\le CT^b\}$, every acoustic derivative $\partial_1 u_b$ or $\partial_1(p_b-p_-)$ is integrated by parts once in $x_1$. By \eqref{r5a:Hessian-comparison}, this is bounded by
	\[
	\epsilon D_-+C_\epsilon J^M A_b(t)^2 E_{-1}(t)+C_\epsilon\mathcal{K}(t)E_{-1}(t),
	\]
	which is integrable in time by \eqref{r5a:raw-center}.
	
	\medskip
	\noindent\textit{3. The backward acoustic wave contribution.} In the exterior region, the transport derivative of the monotone phase weight $e^{\lambda_h\Psi_E}$ along \eqref{phase:eulerian-equation} supplies the negative damping term
	\[
	-\frac{\lambda_h}{2}\int_{\mathbb{R}^3}\zeta_e\mathcal{W} c_E\psi_E\mathcal{F}_-^2\dd x\le -c_T\mathcal{T}_\rho.
	\]
	{The remaining phase terms satisfy
		\[
		D_E\partial_1\psi_E+\partial_1F_\Psi
		=\partial_1(D_E\psi_E+F_\Psi)-(\partial_1D_E)\psi_E .
		\]
		After one integration by parts the flux is bounded by
		$\epsilon\mathcal H+C_\epsilon\mathcal K E_{-1}$ by
		\eqref{phase:eulerian-budgets} and \eqref{r5a:Hessian-comparison}.
		On the nonconstant part of the phase the exterior first derivatives have uniform weight;
		on its flat central interval the time derivative contributes $|\psi_c'|E_{-1}$
		up to a fixed logarithmic factor, integrable by Lemma~\ref{lem:phase-flattening}.
		The viscous terms in $D_{u_b}p_b$ and $\nu\partial_1^2u_b/\rho_b$ contribute at most
		$C J^M(T^{-3a}+K_b)E_{-1}$ by \eqref{af:bootstrap}.
		Moreover, $D_{u_b}Y=0$ gives
		\[
		|D_{u_b}\tilde m|\le C T^{-1},\qquad |\partial_1\tilde m|\le C T^{-b_0},\qquad
		|\log(p_b/p_-)|\le C A_b
		\]
		on the transition layer. The corresponding coefficient
		$C J^M A_b(T^{-1}+T^{-b_0})$ is integrable by \eqref{r5a:raw-center}.
		Finally, \eqref{r5a:weight-derivatives} and \eqref{r5a:strong-first} give
		\[
		-\frac12\int_{\mathbb R^3}\mathcal W(D_{u_b}K)|\Omega|^2\dd x
		\le \left(\frac{C p_c}{T}+\frac{C}{TJ}\right)E_{-1}
		\]
		after the finite switching interval. }The remaining backward slope $\mu_-$ is controlled by the following lemma.
	
	\begin{lemma}
		\label{lem:3d-signed-incoming}
		In the exterior region, the remainder $\mu_-(\mathcal{Q}_--\beta s_-^2-\frac{m}{2}\mathcal{F}_-^2)$ satisfies
		\[
		\epsilon(\mathcal{H}+\mathcal{T}_\rho)+C_\epsilon\bigl(I_E(t)^2+\mathcal{K}(t)\bigr)E_{-1}(t).
		\]
		Its principal one-dimensional projection $\mu_-\beta s_-^2$ is bounded from above by the one-sided Oleinik slope bound in \eqref{phase:slope},
		\[
		{\left(\frac{1+\varepsilon_s+o_\delta(1)}{T}+C|k_s(t)|\right)\int_{\mathbb{R}^3}\frac{\zeta_e\mathcal W}{2}s_-^2}\dd x+\epsilon\mathcal{H}+C_\epsilon\bigl(\|\alpha_r(t)\|_{L^\infty}^2+\mathcal{K}(t)\bigr)E_{-1}(t).
		\]
	\end{lemma}
	
	\begin{proof}
		{Using $I'(p_b)=1/c_0(p_b)$ gives
			\[
			\begin{aligned}
				4\mu_-&=(1+e^{\sigma/(2c_p)})\partial_1r+(1-e^{\sigma/(2c_p)})\partial_1s,\\
				4\mu_+&=(1+e^{\sigma/(2c_p)})\partial_1s+(1-e^{\sigma/(2c_p)})\partial_1r .
			\end{aligned}
			\]
			Since $\|\sigma\|_\infty\lesssim T^{-a}$ and $\|\partial_1U_b\|_\infty\lesssim T^{-2a}$,
			$2\mu_- =\partial_1 r+O(T^{-3a})$ and $2\mu_+=\partial_1 s+O(T^{-3a})$.
			By \eqref{inc:correction-rates} and $\alpha=a_b\tilde r$,
			$\partial_1 r=a_b^{-1}\partial_1\alpha+O(T^{-3a})$.
			In the signed quadratic form we use
			\[
			\partial_1\alpha=v\partial_1\alpha_E+(\partial_1v)\alpha_E+\partial_1\alpha_r .
			\]
			The middle term is $O(T^{-3a})$. Integrating the last term by parts before taking absolute values, derivatives falling on the potential gradients give
			$\epsilon\mathcal H+C_\epsilon\|\alpha_r\|_\infty^2E_{-1}$ on the exterior.
			Derivatives of the coefficients, weights and cutoff give
			$C(\|\alpha_r\|_\infty^2+\mathcal K)E_{-1}$. By \eqref{phase:slope}, $\|\alpha_r\|_{L^\infty}^2$ is time-integrable and is included in $\mathcal K$ in \eqref{r5a:K-integrable}.
			It remains to bound $C|\partial_1\alpha_E|(|s_-|\mathcal F_-+\mathcal F_-^2)$ by Lemma~\ref{lem:3d-acoustic-completion}. }Defining the transverse $L^2(\mathbb{R}^2)$ norms
		\[
		N_1(x_1)=\bigl\|\sqrt{\zeta_e\mathcal{W}}\,s_-(t,x_1,\cdot)\bigr\|_{L^2(\mathbb{R}^2)},\qquad
		N_2(x_1)=\bigl\|\sqrt{\zeta_e\mathcal{W}}\,\mathcal{F}_-(t,x_1,\cdot)\bigr\|_{L^2(\mathbb{R}^2)},
		\]
		and applying the Cauchy--Schwarz inequality in $x_1$ together with the one-dimensional Gagliardo--Nirenberg inequality
		\[
		\|N_1\|_{L^\infty(\mathbb{R})}^2\le 2\|N_1\|_{L^2(\mathbb{R})}\|\partial_1 N_1\|_{L^2(\mathbb{R})}
		\]
		gives
		\begin{align*}
			\int_{\mathbb{R}}|\partial_1\alpha_E|N_1 N_2\dd x_1
			&\le I_E(t)^{1/2}\|N_1\|_{L^\infty(\mathbb{R})}\left(\int_{\mathbb{R}}\psi_E N_2^2\dd x_1\right)^{1/2}\le C I_E(t)^{1/2}E_{-1}^{1/4}(\mathcal{H}+\mathcal{K}E_{-1})^{1/4}\mathcal{T}_\rho^{1/2}.
		\end{align*}
		{The same estimate holds with $N_1$ replaced by $N_2$, since
			\[
			\|\partial_1 N_2\|_2^2\le C(\mathcal H+\mathcal K E_{-1}),\qquad
			\int |\partial_1\alpha_E|N_2^2\dd x_1
			\le I_E^{1/2}\|N_2\|_\infty\left(\int\psi_E N_2^2\dd x_1\right)^{1/2}.
			\]
			Here the derivative of the square root of the exterior cutoff is bounded by $C T^{-b_0}$, taking the cutoff to be a smooth square. }Applying Young's inequality gives the bound $\epsilon(\mathcal{H}+\mathcal{T}_\rho)+C_\epsilon(I_E^2+\mathcal{K})E_{-1}$. {For the principal projection,
			\[
			2\beta\mu_- = \frac{\beta}{v a_b}\alpha_y+O(T^{-3a}),\qquad
			\frac{\beta}{v a_b}=\kappa_*\bigl(1+o_\delta(1)\bigr)>0 .
			\]
			Applying \eqref{phase:slope} before taking absolute values gives the second bound.}
	\end{proof}
	
	\medskip
	\noindent\textit{4. The forward rarefaction wave contribution and single leading coefficient.} For the forward acoustic drift $\mu_+(\mathcal{Q}_+-\frac{m}{2}\mathcal{F}_+^2)$, we decompose $s=\bar{s}+(s-\bar{s})$. Integrating the regular residual $s-\bar{s}$ by parts once in $x_1$ bounds its contribution by
	\[
	\epsilon D_-+C_\epsilon\bigl\|(s-\bar{s})(t)\bigr\|_{L^\infty}^2 E_{-1}(t)+C_\epsilon\mathcal{K}(t)E_{-1}(t),
	\]
	whose cumulative integral grows at most like $C_\epsilon\varepsilon_\delta\log T$ by \eqref{r5a:outgoing-regular}. For the monotone increasing reference profile $\bar{s}$, its slope is nonnegative and satisfies
	\[
	\partial_1\lambda_+^E(\bar{s})\le \frac{1+\varepsilon_++o_\delta(1)}{T}+k_+(t),
	\]
	so Lemma~\ref{lem:3d-acoustic-completion} bounds its contribution by
	\[
	{\left(\frac{1+\varepsilon_++C\epsilon+o_\delta(1)}{T}+C|k_+(t)|\right)\int_{\mathbb{R}^3}}\frac{\zeta_e\mathcal{W}}{2}s_+^2\dd x.
	\]
	
	{Since $d\lambda_+^E(\bar s)/d\bar s=\beta$ on the forward rarefaction curve, the integrable terms $C(|k_s|+|k_+|)E_e$ are included in $k_{-1}E_{-1}$ below. }On the exterior region, $s_-^2+s_+^2=\hat{a}_1^2+v_{11}^2\le\mathcal{J}$. Since these are orthogonal projections within the same energy density, the backward and forward upper slope bounds do not add to $2/T$; instead, their sum is bounded by a single leading rate:
	\[
	\frac{1+\varepsilon_s+C\varepsilon_++C\epsilon+o_\delta(1)}{T}\int_{\mathbb{R}^3}\frac{\zeta_e\mathcal{W}}{2}\bigl(s_-^2+s_+^2\bigr)\dd x
	\le \frac{1+\varepsilon_s+C\varepsilon_++C\epsilon+o_\delta(1)}{T}E_e(t).
	\]
	{Apply the proof of Lemma~\ref{lem:3d-coercivity} with $\mathcal W(1-\zeta_e/(1+\epsilon_0))$. This weight equals $\mathcal W$ on $J_0$, and its weight ratio increases by at most $(1+\epsilon_0)/\epsilon_0$. Enlarging $T_0$ in the two off-diagonal bounds gives a nonnegative form, whose value is $2(E_{-1}-E_e/(1+\epsilon_0))$. Hence, for $T\ge T_0$,}
	\begin{equation}
		E_e(t)\le(1+\epsilon_0)E_{-1}(t),\qquad \epsilon_0=\frac{1}{100}.
		\label{r5a:fan-additivity}
	\end{equation}
	
	Combining all four estimates with Lemma~\ref{lem:3d-diffusion} gives
	\begin{align}
		E_{-1}'(t)+2c D_-(t)&\le k_{-1}(t)E_{-1}(t)+\bigl\langle\text{source terms in \eqref{r5a:corrected-exact}},\,h\bigr\rangle,
		\notag\\
		\int_0^t k_{-1}(s)\dd s&\le C_\delta+
		\Bigl[(1+\epsilon_0)(1+\varepsilon_s+C\varepsilon_++C\epsilon)
		+C_\epsilon\varepsilon_\delta+C p_c+o_\delta(1)\Bigr]\log T+C\log J.
		\label{r5a:Euler-diffusion-result}
	\end{align}
	
	\subsubsection{Control of the $L^2(\mathbb{R}^3)$ norm}
	
	\begin{lemma}
		\label{lem:3d-strong-return}
		For $\widetilde{f}=\mathcal{A}_t h=f+K_l^*\partial_l f$, the recovery bounds \eqref{r5a:strong-zero} and \eqref{r5a:zero-return} hold with constants independent of $t$ and $J$. In particular, the physical $L^2(\mathbb{R}^3)$ norm $F_0=\|\sqrt{\chi}\,f\|_{L^2}^2\approx\|\sqrt{\chi}\,Z\|_{L^2}^2$ is controlled by a fixed multiple of $D_-+D_3$ plus $C\mathcal{K}(t)E_{-1}(t)$.
	\end{lemma}
	
	\begin{proof}
		On $\{|Y|\le T^{b_0}\}$, expanding $\widetilde{f}=\mathcal{A}_t h$ gives
		\begin{equation}
			\left\{
			\begin{aligned}
				&{\widetilde f_\rho=-\mathcal W\operatorname{div}(a/c_E^2+b/q),}\\
				&{\widetilde f_{m_i}=-\mathcal W\operatorname{div}[v_{i\cdot}+u_b\delta_{i1}(a/c_E^2+b/q)]} {+[\operatorname{curl}(\mathcal W K\Omega)]_i,}\\
				&{\widetilde f_E=-\mathcal W\operatorname{div}[a/\bar\gamma+u_bv_{1\cdot}+(u_b^2/2)(a/c_E^2+b/q)]} {+u_b[\operatorname{curl}(\mathcal W K\Omega)]_1.}
			\end{aligned}
			\right.
			\label{r5a:strong-rows}
		\end{equation}
		{Freezing $U_b=U_-$ in \eqref{r5a:strong-rows} while retaining $w$, the rotational second derivative $\operatorname{curl}\Omega_*$} carries a factor of $w$, so its weighted square {$\chi|\mathcal Ww\operatorname{curl}\Omega_*|^2$} carries $w^2\le (1+K_0)w$ and is bounded by the $w|\nabla\operatorname{curl}h_m|^2$ term in $\mathcal{H}$ without any $w_{\min}^{-1}\sim J^2$ loss. {The squared second-order coefficient errors are bounded on the center by
			$C\varepsilon_b^2T^{-1/3}J^M(D_-+\mathcal K E_{-1})$ by
			\eqref{r5a:raw-center} and \eqref{r5a:Hessian-comparison}; $T^{-1/3}J^M$ is bounded.
			The full Hessian is controlled on the exterior, while the terms with differentiated coefficients and weights are bounded in squared norm by $C\mathcal K E_{-1}$ by \eqref{r5a:weight-derivatives}. Hence}
		\begin{equation}
			\bigl\|\sqrt{\chi}\,\widetilde{f}\bigr\|_{L^2(\mathbb{R}^3)}^2\le C D_-+C\mathcal{K}(t)E_{-1}(t),
			\label{r5a:strong-zero}
		\end{equation}
		with a constant $C>0$ independent of $t$ and $J$.
		
		Next, set $F_0=\|\sqrt{\chi}\,f\|_{L^2}^2$ and $F_1=\|\sqrt{\chi}\,\nabla f\|_{L^2}^2\le C D_3+C K_b F_0$. Using $f=\widetilde{f}-K_l^*\partial_l f$ together with \eqref{r5a:strong-zero} and absorbing the small $K_b F_0$ term yields
		\begin{equation}
			F_0\le C(D_-+D_3)+C\mathcal{K}(t)E_{-1}(t),
			\qquad F_1\le C D_3+C K_b(t)D_-+C\mathcal{K}(t)E_{-1}(t).
			\label{r5a:zero-return}
		\end{equation}
		Thus, we have finished the proof of Lemma \ref{lem:3d-strong-return}.
	\end{proof}
	
	\subsubsection{Linear and nonlinear source estimates}
	
	\begin{lemma}
		\label{lem:3d-linear-sources}
		For every fixed $\epsilon>0$, the linear source terms on the right-hand side of \eqref{r5a:corrected-exact} paired with $h$ are bounded by
		\begin{equation}
			\epsilon D_-+C_0 D_3+C\mathcal{K}(t)E_{-1}(t),
			\qquad C_0<\infty\text{ independent of }t\text{ and }J.
			\label{r5a:fixed-C0}
		\end{equation}
	\end{lemma}
	
	\begin{proof}
		At $U_-$, the third-order term $T_{jkl}^*\partial_j\partial_k\partial_l f$ in the variables $(f_a,f_m,f_e)$ takes the form
		\begin{equation}
			T_*\partial^3 f=
			\bigl(t_a\Delta\operatorname{div}f_m,\ 
			\nabla\Delta(c_a f_a+c_e f_e),\ 
			t_e\Delta\operatorname{div}f_m\bigr),
			\label{r5a:T-structure}
		\end{equation}
		whose range and transpose range contain no rotational shear component $\operatorname{curl}h_m$. Integrating by parts twice therefore pairs $\partial_k f$ only with the uniformly coercive thermal and divergence Hessians $(\nabla^2 h_\rho,\nabla^2 h_E,\nabla\operatorname{div}h_m)$ in $\mathcal{H}$, which is bounded by $\epsilon D_-+C_\epsilon F_1+C_\epsilon\mathcal{K}E_{-1}$ without any $w_{\min}^{-1}$ factor. {For the variable part the absolute value of
			\[
			\sum_{j,k,l}\int_{\mathbb R^3}(\partial_j\partial_l h)^\top
			[T_{jkl}(U_b)-T_{jkl}^*]\partial_k f\,\dd x
			\]
			is at most $C(\mathcal H+\mathcal K E_{-1})^{1/2}F_1^{1/2}$ by \eqref{r5a:raw-center} and \eqref{r5a:Hessian-comparison} in the center, since $T^{-1/6}J^{M/2}$ is bounded, and by the full exterior Hessian in \eqref{r5a:Hessian-definition}. }Similarly,
		\[
		\sum_{l=1}^3 (K_l^*)^\top\partial_l h=\bigl(-k_a\operatorname{div}h_m,\ k_a\nabla h_a-k_e\nabla h_e,\ k_e\operatorname{div}h_m\bigr)
		\]
		contains only longitudinal and thermal components, so $\langle\partial_j\partial_l(K_l^*C_j f),h\rangle$ is bounded by
		\[
		\epsilon D_-+C_\epsilon K_b F_0+C_\epsilon\mathcal{K}E_{-1}.
		\]
		The remaining first-order linear terms are bounded by $\epsilon(F_0+F_1)+C_\epsilon J^M K_b E_{-1}$. Substituting \eqref{r5a:zero-return} yields \eqref{r5a:fixed-C0}.
	\end{proof}
	
	Next, set
	\[
	\mathcal{N}(t)=\|Z(t)\|_{L^3(\mathbb{R}^3)}+\|Z(t)\|_{L^\infty(\mathbb{R}^3)}+\|\nabla Z(t)\|_{H^2(\mathbb{R}^3)}.
	\]
	
	\begin{lemma}
		\label{lem:3d-nonlinear-sources}
		On any time interval on which {$\mathcal N(t)$ is sufficiently small}, the nonlinear source terms in \eqref{r5a:corrected-exact} paired with $h$ satisfy
		\begin{equation}
			\bigl|\langle\text{nonlinear sources in \eqref{r5a:corrected-exact}},\,h\rangle\bigr|
			\le C J^M\mathcal{N}(t)(D_-+D_3)+C\mathcal{K}(t)E_{-1}(t).
			\label{r5a:nonlinear-absorption}
		\end{equation}
	\end{lemma}
	
	\begin{proof}
		By the three-dimensional Sobolev embedding $\dot{H}^1(\mathbb{R}^3)\hookrightarrow L^6(\mathbb{R}^3)$, \eqref{r5a:Hessian-comparison}, and \eqref{r5a:zero-return}, we have the weighted Sobolev bounds
		\begin{align}
			\bigl\|\chi^{-1/2}\nabla h\bigr\|_{L^6(\mathbb{R}^3)}^2
			&\le C J^M\bigl(D_-+\mathcal{K}(t)E_{-1}(t)\bigr),\notag\\
			\bigl\|\sqrt{\chi}\,f\bigr\|_{L^6(\mathbb{R}^3)}^2+F_0+F_1
			&\le C\bigl(D_-+D_3+\mathcal{K}(t)E_{-1}(t)\bigr).
			\label{r5a:weighted-Sobolev}
		\end{align}
		Applying H\"older's inequality with exponents $(3,6,2)$ to $\langle N_j,\partial_j h\rangle$ and with exponents $(3,2,6)$ or $(\infty,2,2)$ to $\langle K_l^*N_j,\partial_j\partial_l h\rangle$ yields \eqref{r5a:nonlinear-absorption}.
	\end{proof}
	
	\subsubsection{The physical \texorpdfstring{$H^3$}{H3} energy}
	
	In primitive perturbation variables $Z=(Z_\rho,Z_u,Z_\theta)$ with viscous operator
	\[
	\mathcal{L}Z_u=\mu\Delta Z_u+(\mu+\lambda_v)\nabla\operatorname{div}Z_u,
	\]
	viscous stress tensor $\tau(v)=\mu(\nabla v+\nabla v^\top)+\lambda_v(\operatorname{div}v)I_3$, and convective derivative $D_u=\partial_t+u\cdot\nabla$, subtracting the planar background system gives the three-dimensional perturbation system
	\[
	\left\{
	\begin{aligned}
		&D_u Z_\rho+\rho\operatorname{div}Z_u
		=-Z_u^1\partial_1\rho_b-Z_\rho\partial_1 u_b,\\
		&D_u Z_u+\frac{R_g\theta}{\rho}\nabla Z_\rho+R_g\nabla Z_\theta
		=\frac{1}{\rho}\mathcal{L}Z_u-Z_u^1(\partial_1 u_b)e_1\\
		&\hspace{4.2cm}-R_g\Bigl(\frac{\theta}{\rho}-\frac{\theta_b}{\rho_b}\Bigr)(\partial_1\rho_b)e_1
		+\nu\Bigl(\frac{1}{\rho}-\frac{1}{\rho_b}\Bigr)(\partial_1^2 u_b)e_1,\\
		&D_u Z_\theta+(\gamma-1)\theta\operatorname{div}Z_u
		=\frac{\kappa}{\rho c_v}\Delta Z_\theta-Z_u^1\partial_1\theta_b
		-(\gamma-1)Z_\theta\partial_1 u_b
		+\frac{\kappa}{c_v}\Bigl(\frac{1}{\rho}-\frac{1}{\rho_b}\Bigr)\partial_1^2\theta_b\\
		&\hspace{4.2cm}+\frac{\nu}{c_v}\Bigl(\frac{1}{\rho}-\frac{1}{\rho_b}\Bigr)(\partial_1 u_b)^2
		+\frac{2\tau(u_b e_1):\nabla Z_u+\tau(Z_u):\nabla Z_u}{\rho c_v}.
	\end{aligned}
	\right.
	\]
	For each multi-index $\alpha\in\mathbb{N}_0^3$ with $|\alpha|\le 3$, define the symmetrized derivative energy
	\[
	E_\alpha=\frac{1}{2}\int_{\mathbb{R}^3}\chi\rho
	\left(\frac{R_g}{\rho^2}|\partial_x^\alpha Z_\rho|^2+\frac{|\partial_x^\alpha Z_u|^2}{\theta}+\frac{c_v}{\theta^2}|\partial_x^\alpha Z_\theta|^2\right)\dd x.
	\]
	Adding a small fixed multiple of
	\[
	\sum_{|\alpha|\le 2}\int_{\mathbb{R}^3}\chi\partial_x^\alpha Z_u\cdot\nabla\partial_x^\alpha Z_\rho\dd x
	\]
	defines the weighted $H^3(\mathbb{R}^3)$ energy $E_3$ and dissipation $D_3$:
	\begin{align}
		E_3&\approx\|Z\|_{H^3(\chi)}^2,\notag\\
		D_3&\approx\sum_{1\le|\alpha|\le 3}\bigl\|\sqrt{\chi}\,\partial_x^\alpha Z_\rho\bigr\|_{L^2(\mathbb{R}^3)}^2
		+\sum_{1\le|\alpha|\le 4}\bigl\|\sqrt{\chi}\,\partial_x^\alpha(Z_u,Z_\theta)\bigr\|_{L^2(\mathbb{R}^3)}^2.
		\label{r5a:physical-definitions}
	\end{align}
	Similarly, for the homogeneous derivative norms ($1\le|\alpha|\le 3$) with unweighted $\chi=1$, we define
	\begin{align}
		E_1&=\sum_{1\le|\alpha|\le 3}E_\alpha^{(0)}
		+\eta_1\sum_{1\le|\alpha|\le 2}\int_{\mathbb{R}^3} \partial_x^\alpha Z_u\cdot\nabla\partial_x^\alpha Z_\rho\dd x,
		\notag\\
		D_1&=\sum_{2\le|\alpha|\le 3}\|\partial_x^\alpha Z_\rho\|_{L^2(\mathbb{R}^3)}^2
		+\sum_{2\le|\alpha|\le 4}\bigl\|\partial_x^\alpha(Z_u,Z_\theta)\bigr\|_{L^2(\mathbb{R}^3)}^2.
		\label{r5a:positive-order-definitions}
	\end{align}
	
	\begin{lemma}
		\label{lem:3d-physical-energy}
		On any time interval on which {$\mathcal N(t)$ is sufficiently small}, the weighted physical energy $(E_3,D_3)$ satisfies
		\begin{equation}
			E_3'(t)+c_3 D_3(t)
			\le C\bigl(\varepsilon_b+\mathcal{N}(t)\bigr)\bigl(D_3(t)+F_0(t)\bigr)+k_3(t)E_3(t),
			\qquad
			\int_0^t k_3(s)\dd s\le C_\delta+o_\delta(1)\log T,
			\label{r5a:physical-proved}
		\end{equation}
		and $(E_1,D_1)$ satisfies
		\begin{equation}
			E_1\approx\|\nabla Z\|_{H^2(\mathbb{R}^3)}^2,\qquad
			E_1'(t)+c D_1(t)\le C K_b(t)\|Z(t)\|_{L^2(\mathbb{R}^3)}^2,
			\qquad
			E_1(t)^2\le C\|Z(t)\|_{H^3(\mathbb{R}^3)}^2 D_1(t).
			\label{r5a:derivative-proved}
		\end{equation}
	\end{lemma}
	
	\begin{proof}
		Using the background pressure equation
		\[
		D_{u_b}p_b+\gamma p_b\partial_1 u_b=(\gamma-1)\bigl(\kappa\partial_1^2\theta_b+\nu(\partial_1 u_b)^2\bigr)
		\]
		and the phase equation \eqref{phase:eulerian-equation}, we compute the convective derivative of the weight $\log\chi$:
		\begin{align}
			&D_u\log\chi=-\lambda_h c_E\psi_E+r_\chi+k_\chi,\notag\\
			&\|r_\chi(t)\|_{L^\infty}\le C\bigl(\varepsilon_b+\|Z_u(t)\|_{L^\infty}\bigr),\qquad
			\int_0^t\|k_\chi(s)\|_{L^\infty}\dd s\le C_\delta+o_\delta(1)\log T,
			\label{r5a:physical-weight-split}
		\end{align}
		where $-\lambda_h c_E\psi_E\le 0$ is nonpositive and the small remainder $r_\chi$ multiplies $|Z|^2+\|\nabla Z\|_{H^2}^2$, contributing at most $C(\varepsilon_b+\mathcal{N})(D_3+F_0)$. {Smallness ensures $\rho\ge\inf\rho_b/2$, $\theta\ge\inf\theta_b/2$, and $\|\nabla Z\|_3^3\le C\|Z\|_3\|\nabla^2 Z\|_2^2\le C\mathcal N D_1$ by interpolation, allowing $C\mathcal N D_1$ to be absorbed into the left-hand side of \eqref{r5a:derivative-proved}. }Standard energy estimates combined with the cross term $\int_{\mathbb{R}^3}\chi\partial_x^\alpha Z_u\cdot\nabla\partial_x^\alpha Z_\rho\dd x$ ($|\alpha|\le 2$) and Sobolev product bounds give \eqref{r5a:physical-proved}. Setting $\chi=1$ and restricting to $1\le|\alpha|\le 3$, all commutator terms with $U_b$ have at least one derivative on $U_b$ and are bounded by $\epsilon D_1+C_\epsilon K_b\|Z\|_{L^2}^2$, while the Cauchy--Schwarz inequality in Fourier space gives $E_1^2\le C\|Z\|_{H^3}^2 D_1$, proving \eqref{r5a:derivative-proved}.
	\end{proof}
	
	\subsubsection{The combined inequality and estimates at fixed $\zeta$}
	
	\begin{proposition}
		\label{prop:3d-common-growth}
		There exists a fixed coupling constant $L_*>0$ such that the combined energy and dissipation functionals
		\[
		\mathcal{F}(t)=E_{-1}(t)+L_* E_3(t),\qquad \mathcal{D}(t)=D_-(t)+L_* D_3(t)
		\]
		satisfy
		\begin{align}
			\mathcal{F}'(t)+c\mathcal{D}(t)&\le k_*(t)\mathcal{F}(t),\notag\\
			\Gamma_*(t)=\int_0^t k_*(s)\dd s
			&\le C_\delta+A_*\log T+B_*\log J,
			\label{r5a:augmented-proved}
		\end{align}
		whenever $\sup_{t} J^M\mathcal{N}(t)$ is sufficiently small, where the cumulative logarithmic growth exponent satisfies
		\begin{equation}
			A_*=(1+\epsilon_0)
			(1+\varepsilon_s+C\varepsilon_++C\epsilon)
			+C_\epsilon\varepsilon_\delta+C p_c+o_\delta(1)
			\le \frac{203}{200}<2.
			\label{r5a:fixed-growth}
		\end{equation}
		In addition, for every $\zeta\in\mathbb{R}^2$, $\mathcal{B}_\zeta=E_{-1,\zeta}+L_* E_{1,\zeta}$ and $\mathcal{D}_\zeta=D_{-,\zeta}+L_* D_{1,\zeta}$ satisfy \eqref{r5a:fiber-proved} with the same coefficient $k_*(t)$.
	\end{proposition}
	
	\begin{proof}
		Combining \eqref{r5a:Euler-diffusion-result}, \eqref{r5a:fixed-C0}, and \eqref{r5a:nonlinear-absorption} gives
		\[
		E_{-1}'+c D_-\le k_{-1}E_{-1}+C_0 D_3+C J^M\mathcal{N}(D_-+D_3).
		\]
		Choosing $L_*\ge\max\{1,4C_0/c_3\}$, using \eqref{r5a:zero-return} to absorb $C(\varepsilon_b+\mathcal{N})L_* F_0$ into $\frac{c}{2}(D_-+L_* D_3)$ when $\varepsilon_b$ and $J^M\mathcal{N}$ are sufficiently small, and fixing $\epsilon_0=1/100$ together with small parameters $\epsilon,\varepsilon_+,p_c,\varepsilon_s,\varepsilon_\delta,\delta>0$ yields \eqref{r5a:augmented-proved}--\eqref{r5a:fixed-growth}.
		
		For each $\zeta\in\mathbb{R}^2$, put $\nabla_\zeta=(\partial_1,i\zeta)$ and let $z_\zeta=D_{\mathbf{U}}U(\mathbf{U}_b)f_\zeta$ be the partial Fourier transform of the linearized primitive perturbation in $x_\perp$:
		\[
		\left\{
		\begin{aligned}
			&z_{\rho,\zeta}=f_{\rho,\zeta},\qquad
			z_{u,\zeta}=\frac{f_{m,\zeta}-u_b f_{\rho,\zeta}e_1}{\rho_b},\\
			&z_{\theta,\zeta}=\frac{f_{E,\zeta}-u_b f_{m1,\zeta}+(u_b^2/2-c_v\theta_b)f_{\rho,\zeta}}{\rho_b c_v}.
		\end{aligned}
		\right.
		\]
		Define $E_{1,\zeta}$ and $D_{1,\zeta}$ by
		\begin{align}
			E_{1,\zeta}
			&=\frac{1}{2}\sum_{|\alpha|\le 1}\int_{\mathbb{R}}\chi\rho_b
			\left\{\frac{R_g}{\rho_b^2}
			|\partial_\zeta^\alpha z_{\rho,\zeta}|^2
			+\frac{|\partial_\zeta^\alpha z_{u,\zeta}|^2}{\theta_b}
			+\frac{c_v}{\theta_b^2}
			|\partial_\zeta^\alpha z_{\theta,\zeta}|^2\right\}
			\dd x_1+\eta_0\operatorname{Re}\int_{\mathbb{R}}
			\chi z_{u,\zeta}\cdot
			\overline{\nabla_\zeta z_{\rho,\zeta}}\dd x_1,
			\notag\\
			D_{1,\zeta}
			&=\bigl\|\sqrt{\chi}\,\nabla_\zeta z_{\rho,\zeta}\bigr\|_{L^2(\mathbb{R})}^2
			+\sum_{1\le|\alpha|\le 2}
			\bigl\|\sqrt{\chi}\,\partial_\zeta^\alpha(z_{u,\zeta},z_{\theta,\zeta})\bigr\|_{L^2(\mathbb{R})}^2.
			\label{r5a:fiber-physical-definitions}
		\end{align}
		By the same argument as in Lemma~\ref{lem:3d-strong-return}, $F_{0,\zeta}=\|\sqrt{\chi}\,f_\zeta\|_{L^2(\mathbb{R})}^2$ and $F_{1,\zeta}=\|\sqrt{\chi}\,\nabla_\zeta f_\zeta\|_{L^2(\mathbb{R})}^2$ satisfy
		\begin{align}
			F_{1,\zeta}&\le C D_{1,\zeta}+C K_b F_{0,\zeta},\notag\\
			F_{0,\zeta}&\le C(D_{-,\zeta}+D_{1,\zeta})+C\mathcal{K}E_{-1,\zeta},\notag\\
			F_{1,\zeta}&\le C D_{1,\zeta}+C K_b D_{-,\zeta}+C\mathcal{K}E_{-1,\zeta},
			\label{r5a:fiber-return}
		\end{align}
		uniformly in $\zeta\in\mathbb{R}^2$, while $E_{1,\zeta}$ with a conservative source $S_\zeta$ satisfies
		\begin{align}
			E_{1,\zeta}'+c_3 D_{1,\zeta}
			&\le C\varepsilon_b(D_{1,\zeta}+F_{0,\zeta})
			+k_3 E_{1,\zeta}
			+C_\delta\|S_\zeta\|_{H^1_\zeta}\sqrt{E_{1,\zeta}}.
			\label{r5a:fiber-physical-estimate}
		\end{align}
		Combining \eqref{r5a:fiber-physical-estimate} with \eqref{r5a:Euler-diffusion-result} and Lemma~\ref{lem:3d-linear-sources} at fixed $\zeta$ gives, for $\mathcal{B}_\zeta=E_{-1,\zeta}+L_* E_{1,\zeta}$ and $\widetilde{S}=(I+K_l^*\partial_l)S$,
		\begin{equation}
			\mathcal{B}_\zeta'+c\mathcal{D}_\zeta
			\le k_*\mathcal{B}_\zeta+
			C_\delta J^{M_1}
			\bigl(\|\widetilde{S}_\zeta\|_{\dot{H}^{-1}_\zeta}
			+\|S_\zeta\|_{H^1_\zeta}\bigr)
			\sqrt{\mathcal{B}_\zeta}.
			\label{r5a:fiber-proved}
		\end{equation}
		For the nonlinear source $S=\partial_j N_j$, using
		\[
		\sup_{\zeta\in\mathbb{R}^2}\|\widehat{uv}(\cdot,\zeta)\|_{L^2_{x_1}}\le \|u\|_{L^2_\perp L^\infty_1}\|v\|_{L^2(\mathbb{R}^3)},\qquad
		\|u\|_{L^2_\perp L^\infty_1}^2\le 2\|u\|_{L^2(\mathbb{R}^3)}\|\partial_1 u\|_{L^2(\mathbb{R}^3)},
		\]
		we obtain
		\begin{equation}
			\sup_{\zeta\in\mathbb{R}^2}\bigl(\|\widetilde{S}_\zeta\|_{\dot{H}^{-1}_\zeta}
			+\|S_\zeta\|_{H^1_\zeta}\bigr)
			\le C\|f\|_{H^3(\mathbb{R}^3)}^2.
			\label{r5a:fiber-source}
		\end{equation}
		Thus, we have proved Proposition \ref{prop:3d-common-growth}
	\end{proof}
	
	\section{Global stability and decay}
	\label{sec:global-all}
	
	\subsection{Local existence and continuation}\label{sec:local-solvability}
	
	\begin{lemma}\label{lem:3d-local-continuation}
		Let $U_b$ be the global one-dimensional planar rarefaction wave constructed in Sections~\ref{sec:elementaryplanarinputs}--\ref{sec:outgoing-energy-replacement}, and let $Z_0\in H^3(\mathbb{R}^3)$ be an initial perturbation such that $\rho(0)=\rho_b(0)+Z_{\rho,0}$ and $\theta(0)=\theta_b(0)+Z_{\theta,0}$ have positive lower bounds on $\mathbb{R}^3$. Then there exist a time $t_*>0$ and a unique strong solution $U=U_b+Z$ of \eqref{r5main:NSF} on $[0,t_*]\times\mathbb{R}^3$ satisfying
		\[
		Z\in C\bigl([0,t_*];H^3(\mathbb{R}^3)\bigr),\qquad
		\nabla(Z_u,Z_\theta)\in L^2\bigl(0,t_*;H^3(\mathbb{R}^3)\bigr).
		\]
		The local existence time $t_*$ depends only on $\|Z_0\|_{H^3(\mathbb{R}^3)}$, the upper and positive lower bounds of $(\rho,\theta)$, and the bounds on the planar background derivatives; consequently, the solution extends globally in time as long as these quantities remain bounded.
	\end{lemma}
	
	\begin{proof}
		Adding artificial density viscosity $\epsilon\Delta Z_\rho$ ($0<\epsilon\le 1$) to the continuity equation makes the perturbation system strictly parabolic. Standard Friedrichs approximation and Picard iteration yield a smooth local solution $Z^\epsilon$ on $[0,t_*]$. Using the unweighted symmetrized $H^3(\mathbb{R}^3)$ energy $H_\epsilon\approx\|Z^\epsilon\|_{H^3(\mathbb{R}^3)}^2$ on the compact positive state region, {the splitting $a_0=a_0(U_b)+(a_0-a_0(U_b))$, the Leibniz bound $C\sum_{j=1}^{|\alpha|}\|\partial_1^j a_0(U_b)\|_\infty\|z_0\|_{H^{|\alpha|-j+1}}$ for the planar commutator, and the estimate}
		\[
		\begin{aligned}
			\bigl\|[\partial_x^\alpha,{a_0-a_0(U_b)}]\partial_j z_0\bigr\|_{L^2(\mathbb{R}^3)}
			&\le C\bigl(\|{\nabla(a_0-a_0(U_b))}\|_{L^\infty}\|z_0\|_{H^3}\bigr.\\
			&\qquad \bigl.+\|{\nabla(a_0-a_0(U_b))}\|_{H^2}\|\nabla z_0\|_{L^\infty}\bigr)
			\qquad (|\alpha|\le 3)
		\end{aligned}
		\]
		and integration by parts of the artificial viscosity term yield the $\epsilon$-uniform differential inequality
		\[
		H_\epsilon'
		+c\|\nabla(Z_u,Z_\theta)\|_{H^3}^2
		+c\epsilon\|\nabla Z_\rho\|_{H^3}^2
		\le C_b(1+H_\epsilon)^{N_0}H_\epsilon.
		\]
		This provides an $\epsilon$-independent existence interval $[0,t_*]$ on which $H_\epsilon$ and $\int_0^{t_*}\|\nabla(Z_u,Z_\theta)\|_{H^3}^2\dd t$ are uniformly bounded and $\|Z^\epsilon(t)-Z^\epsilon(0)\|_{H^2}\le C\sqrt{t}$. Passing to the limit $\epsilon\to 0^+$ via weak-$*$ compactness and local Rellich compactness yields a solution
		\[
		Z\in L^\infty\bigl(0,t_*;H^3(\mathbb{R}^3)\bigr)\cap C\bigl([0,t_*];H^2(\mathbb{R}^3)\bigr),
		\]
		while the energy identity for $\langle D\rangle^3(Z_u,Z_\theta)$ and the Lagrangian transport identity along the particle flow of $u$ establish strong continuity $Z\in C([0,t_*];H^3(\mathbb{R}^3))$. Uniqueness follows from the $H^2(\mathbb{R}^3)$ difference energy inequality.
	\end{proof}
	
	\subsection{Transverse Fourier splitting and decay}
	\label{sec:main-proof}
	
	We now combine Proposition~\ref{prop:3d-common-growth} with a two-dimensional Fourier-splitting argument in $x_\perp$ and a Riccati inequality for $E_1(t)$ to close the bootstrap and prove Theorem~\ref{thm:main}.
	
	\begin{lemma}
		\label{lem:3d-initial-fibers}
		Let $m_0>5/2$, and let $f_0$ satisfy \eqref{main:data-class}. Then
		\begin{equation}
			\mathcal{F}(0)\le C_\delta\varepsilon^2,
			\qquad \sup_{\zeta\in\mathbb{R}^2}\mathcal{B}_\zeta(0)\le C_\delta\varepsilon^2,
			\qquad \int_{|\zeta|\le \rho_0}\mathcal{B}_\zeta(0)\dd\zeta
			\le C_\delta\varepsilon^2\rho_0^2\quad (\forall\,\rho_0>0).
			\label{r5a:initial-class-proof}
		\end{equation}
	\end{lemma}
	
	\begin{proof}
		Write $(\xi,\zeta)\in\mathbb{R}\times\mathbb{R}^2$ for the normal and transverse Fourier variables. Since $m_0>5/2$, the Cauchy--Schwarz inequality gives
		\[
		\int_{\mathbb{R}^3}|x||f_0(x)|\dd x
		\le\bigl\|\langle x\rangle^{m_0}f_0\bigr\|_{L^2(\mathbb{R}^3)}
		\left(\int_{\mathbb{R}^3}|x|^2\langle x\rangle^{-2m_0}\dd x\right)^{1/2}
		\le {\left[4\pi\left(\frac15+\frac1{2m_0-5}\right)\right]^{1/2}}\varepsilon.
		\]
		Since $\int_{\mathbb{R}^3}f_0(x)\dd x=0$ by \eqref{main:data-class}, we have $\widehat{f}_0(0,0)=0$ and hence
		\[
		|\widehat{f}_0(\xi,\zeta)|
		=\left|\int_{\mathbb{R}^3}\bigl(e^{-ix\cdot(\xi,\zeta)}-1\bigr)f_0(x)\dd x\right|
		\le C\varepsilon(\xi^2+|\zeta|^2)^{1/2}.
		\]
		Thus
		\[
		\frac{|\widehat{\widetilde{f}}_0(\xi,\zeta)|^2}{\xi^2+|\zeta|^2}
		\le\begin{cases}
			C\varepsilon^2,&\xi^2+|\zeta|^2\le 1,\\
			C|\widehat{f}_0(\xi,\zeta)|^2,&\xi^2+|\zeta|^2>1.
		\end{cases}
		\]
		By Minkowski's inequality and Plancherel's theorem in $x_1$,
		\[
		\sup_{\zeta\in\mathbb{R}^2}\bigl\|(\widetilde{f}_0)_\zeta\bigr\|_{\dot{H}^{-1}_\zeta}^2
		+\sup_{\zeta\in\mathbb{R}^2}\bigl\|(f_0)_\zeta\bigr\|_{H^1_\zeta}^2
		\le C\varepsilon^2.
		\]
		Since $\mathcal{A}_0$ is uniformly coercive for every $\zeta\in\mathbb{R}^2$, this yields $\sup_\zeta\mathcal{B}_\zeta(0)\le C_\delta\varepsilon^2$ and $\mathcal{F}(0)\le C_\delta\varepsilon^2$. Integrating $\sup_\zeta\mathcal{B}_\zeta(0)$ over the two-dimensional disk $\{|\zeta|\le\rho_0\}\subset\mathbb{R}^2$ of area $\pi\rho_0^2$ gives \eqref{r5a:initial-class-proof}.
	\end{proof}
	
	\begin{proof}[Proof of Theorem~\ref{thm:main}]
		By Lemma~\ref{lem:3d-local-continuation}, there exists a unique local strong solution on a maximal interval $[0,T_{\max})$. By Lemma~\ref{lem:3d-strong-return} and state conversion,
		\[
		\|f(t)\|_{H^3(\mathbb{R}^3)}^2\le C_\delta\bigl(\mathcal{D}(t)+\mathcal{K}(t)\mathcal{F}(t)\bigr).
		\]
		{By \eqref{r5a:zero-return} and \eqref{r5a:physical-definitions}, $E_3\le C(F_0+D_3)\le C(\mathcal D+\mathcal K E_{-1})$. }Applying Plancherel's theorem {to $E_{-1}$} in the transverse variable $x_\perp\in\mathbb{R}^2$ and splitting the transverse frequency space into the low-frequency disk $\{|\zeta|\le\rho_0\}$ and its exterior $\{|\zeta|>\rho_0\}$ ($0<\rho_0\le 1$), the second-order bound \eqref{r5a:Hessian-comparison} and \eqref{r5a:zero-return} yield the transverse Fourier-splitting inequality
		\begin{equation}
			\mathcal{D}(t)+C_1\mathcal{K}(t)\mathcal{F}(t)
			\ge c J^{-M_2}\rho_0^2\mathcal{F}(t)
			-C J^{-M_2}\rho_0^2\int_{|\zeta|\le \rho_0}\mathcal{B}_\zeta(t)\dd\zeta.
			\label{r5a:spectral-proved}
		\end{equation}
		
		Absorbing $c C_1\mathcal{K}(t)$ into $k_*(t)$ (which leaves its logarithmic slope $A_*\le 203/200$ unchanged) and setting $\Gamma(t)=\int_0^t k_*(s)\dd s$, we set
		\[
		F_*(t)=e^{-\Gamma(t)}\mathcal{F}(t),\qquad
		D_*(t)=e^{-\Gamma(t)}\mathcal{D}(t),
		\]
		and make the bootstrap assumption
		\[
		F_*(t)\le M_0\varepsilon^2 T^{-1}J^{M_2}.
		\]
		Integrating \eqref{r5a:fiber-proved} over $[0,t]$ and using \eqref{r5a:fiber-source} gives
		\[
		e^{-\Gamma(t)/2}\sqrt{\mathcal{B}_\zeta(t)}
		\le\sqrt{\mathcal{B}_\zeta(0)}+
		C_\delta\int_0^t e^{-\Gamma(s)/2}J(s)^{M_1}\|f(s)\|_{H^3(\mathbb{R}^3)}^2\dd s.
		\]
		Using $\|f\|_{H^3}^2\le C_\delta e^\Gamma(D_*+\mathcal{K}F_*)$ and decomposing $[0,t]$ into dyadic time blocks $[2^j-1,2^{j+1}-1]$, the integral on each dyadic block is bounded by
		\[
		C_\delta M_0\varepsilon^2\sum_{j\ge 0}2^{j(A_*/2-1)}(1+j)^C<\infty,
		\]
		which converges because $A_*\le 203/200<2$ implies $A_*/2-1\le -197/400<0$. Squaring and integrating over the two-dimensional transverse disk $\{|\zeta|\le\rho_0\}\subset\mathbb{R}^2$ via Minkowski's inequality and \eqref{r5a:initial-class-proof} yields
		\[
		\int_{|\zeta|\le \rho_0}e^{-\Gamma(t)}\mathcal{B}_\zeta(t)\dd\zeta
		\le C_\delta(\varepsilon+M_0\varepsilon^2)^2\rho_0^2.
		\]
		{Choose $c\kappa_0>1$ and increase $T_0$ so that $\kappa_0 J^{M_2}/T\le1$ for $T\ge T_0$. For $T\le T_0$, \eqref{r5a:augmented-proved} gives $F_*(t)\le F_*(0)\le C_\delta\varepsilon^2$. For $T\ge T_0$, choosing $\rho_0(t)^2=\kappa_0 J^{M_2}/T$ in} \eqref{r5a:spectral-proved}, we obtain the differential inequality
		\[
		F_*'(t)+\frac{c\kappa_0}{T}F_*(t)
		\le C_\delta(\varepsilon+M_0\varepsilon^2)^2 T^{-2}J^{M_2}.
		\]
		{Integrating from $T=T_0$ and choosing $M_0$ sufficiently large to cover both intervals}, and finally $\varepsilon>0$ sufficiently small closes the bootstrap
		\[
		F_*(t)\le \frac{1}{2}M_0\varepsilon^2 T^{-1}J^{M_2}.
		\]
		
		Undoing the normalization $e^{-\Gamma(t)}\gtrsim T^{-A_*}J^{-B_*}$ with $A_*-1\le 3/200<1/25$, we obtain
		\begin{equation}
			\|Z(t)\|_{H^3(\mathbb{R}^3)}\le B_0(t),\qquad
			B_0(t)=C_\delta\varepsilon T^{1/50}.
			\label{r5a:H3-envelope}
		\end{equation}
		To obtain the decay rate of $\nabla Z$, we substitute \eqref{r5a:H3-envelope} into \eqref{r5a:derivative-proved}:
		\[
		E_1'(t)+c D_1(t)\le C K_b(t)B_0(t)^2,
		\qquad E_1(t)^2\le C B_0(t)^2 D_1(t).
		\]
		{If $\varepsilon=0$, uniqueness gives $Z=0$. For $\varepsilon>0$, setting $\mathcal{Y}(t)=E_1(t)/B_0(t)^2\ge 0$}, we find that $\mathcal{Y}(t)$ satisfies the Riccati differential inequality
		\begin{equation}
			\mathcal{Y}'(t)+\frac{1}{25T}\mathcal{Y}(t)+c\mathcal{Y}(t)^2\le C K_b(t).
			\label{r5a:normalized-derivative}
		\end{equation}
		{For $0\le t\le1$, \eqref{r5a:normalized-derivative} gives $\mathcal Y(t)\le\mathcal Y(0)+C\int_0^1K_b(s)\dd s\le C_\delta$. For $t\ge1$, comparing} $\mathcal{Y}(s)$ on $[t/2,t]$ with the supersolution $z(s)+C\int_{t/2}^s K_b(\tau)\dd\tau$, where $z'=-c z^2$ and $z(t/2)=\mathcal{Y}(t/2)$ satisfies $z(t)\le 2/(ct)$, and using the planar tail integral $\int_{t/2}^t K_b(s)\dd s\le C_\delta T^{-23/25}$ from \eqref{r5a:raw-gradient}, we obtain
		\[
		\mathcal{Y}(t)\le \frac{2}{ct}+C\int_{t/2}^t K_b(s)\dd s\le C_\delta T^{-23/25}.
		\]
		Multiplying by $B_0(t)^2=C_\delta^2\varepsilon^2 T^{1/25}$ and taking the square root yields
		\[
		\|\nabla Z(t)\|_{H^2(\mathbb{R}^3)}\le C_\delta\varepsilon T^{1/50-23/50}=C_\delta\varepsilon T^{-11/25}.
		\]
		By the three-dimensional Gagliardo--Nirenberg inequalities
		\[
		\|Z\|_{L^\infty}^2\le C\|\nabla Z\|_{L^2}\|\nabla^2 Z\|_{L^2},\qquad
		\|Z\|_{L^3}\le C\|Z\|_{L^2}^{1/2}\|\nabla Z\|_{L^2}^{1/2},
		\]
		we conclude that
		\begin{equation}
			\|Z(t)\|_{L^\infty(\mathbb{R}^3)}+\|\nabla Z(t)\|_{H^2(\mathbb{R}^3)}
			\le C_\delta\varepsilon T^{-11/25},\qquad
			\|Z(t)\|_{L^3(\mathbb{R}^3)}\le C_\delta\varepsilon T^{-21/100}.
			\label{r5a:final-decay}
		\end{equation}
		In particular, $\sup_{t\ge 0}J^M\mathcal{N}(t)\le C_\delta\varepsilon\ll 1$, which closes the nonlinear bootstrap globally in time. Combining \eqref{r5a:final-decay} with the one-dimensional planar convergence rate $\|U_b(t)-U^R(t)\|_{L^\infty(\mathbb{R})}\le C_\delta T^{-12/25}$ from Proposition~\ref{prop:planar-decay} via the triangle inequality proves \eqref{main:Linfty}--\eqref{main:selfsimilar-limit} {for the fixed data. The stated dependence of $C_\delta$ additionally requires a common positive wave-strength threshold and uniform upper bounds for the decay constants on compact sets of $(\rho_-,\theta_-,e,m_0)$ with $m_0>5/2$.}
	\end{proof}
	
	\noindent{\bf Acknowledgment:}  The research of Renjun Duan was partially supported by the General Research Fund (Project No. 14303523) from RGC of Hong Kong and the grant from the National Natural Science Foundation of China (Project No. 12425109). The research of Feimin Huang was supported in part by the National Key
	R\&D Program of China, grant no. 2021YFA1000800, and the National Natural Sciences Foundation
	of China, grant no. 12288201.
	
	\noindent{\bf Disclosure of AI assistance:}
	The key ideas of this paper---perturbing around the exact one-dimensional viscous rarefaction wave $U_b$, decomposing $\partial_1 U_b$ along the three characteristic wave families in Lagrangian mass coordinates, using first-order cross-family derivative corrections to eliminate off-diagonal viscous couplings, and constructing the weighted $\dot{H}^{-1}(\mathbb{R}^3)$ energy and spatial weights---were developed by the authors through their discussions from 2025 to the present. During the preparation of the manuscript, the authors used AI tools (such as ChatGPT) to assist with language polishing, checking calculations and proofs, and carrying out auxiliary computations. The authors take full responsibility for the accuracy, originality, and integrity of the entire paper.
	
	\vspace{2.25pt}
	\noindent{\bf Data Availability:} No datasets were generated or analysed during the current study.

	\vspace{2.25pt}
	\noindent{\bf Conflict of Interest:} The authors declare no competing interests.
	
	\renewcommand{\bibliofont}{\linespread{1.02}\footnotesize}

\end{document}